\documentclass[12pt, uplatex]{amsart}
\usepackage{amsmath,amsfonts,amsthm,amscd,amssymb,mathrsfs,amssymb,bm,pb-diagram, here, extsizes, midpage}
\usepackage{setspace}
\usepackage[all,cmtip]{xy}

\usepackage{blkarray}
\usepackage{multirow}
\usepackage{mathtools}
\usepackage{graphicx}

\usepackage{mathrsfs, arydshln}
\usepackage{graphicx}
\usepackage[all]{xy}
\newtheorem{theorem}{Theorem}

\newtheorem{proposition}[theorem]{Proposition}
\newtheorem{lemma}[theorem]{Lemma}
\newtheorem{definition}[theorem]{Definition}

\newtheorem{corollary}[theorem]{Corollary}

\newtheorem{remark}[theorem]{Remark}
\newcommand{\aaa}{\alpha}
\newcommand{\bbb}{\beta}

\newcommand{\CCC}{\Gamma}
\newcommand{\ddd}{\delta}
\newcommand{\DDD}{\Delta}

\newcommand{\id}{{\rm{id}}}
\newcommand{\lmd}{\lambda}
\newcommand{\Lmd}{\Lambda}

\newcommand{\CP}{\mathbb{CP}}
\newcommand{\CC}{\mathbb{C}}

\newcommand{\RR}{\mathbb{R}}
\newcommand{\ZZ}{\mathbb{Z}}

\newcommand{\RP}{\mathbb{RP}}

\newcommand{\PGL}{\mathrm{PGL}}

\newcommand{\SU}{{\rm{SU}}}

\newcommand{\Sing}{{\rm{Sing}\,}}
\newcommand{\Supp}{{\rm{Supp}\,}}
\newcommand{\qdr}{\CP_1\times\CP_1}

\newcommand{\upd}{^{(d)}}
\newcommand{\upg}{^{(g)}}

\newcommand{\uptwo}{^{(2)}}

\renewcommand{\hat}{\widehat}
\renewcommand{\tilde}{\widetilde}

\newcommand{\mf}{\mathfrak}

\newcommand{\ol}{\overline}

\newcommand{\lras}{\,\longrightarrow\,}

\newcommand{\set}{\,|\,}

\newcommand{\proofend}{\hfill$\square$}

\newcommand{\inv}{^{-1}}

\newcommand{\Aut}{{\rm{Aut}}}
\newcommand{\Bs}{{\rm{Bs}}}

\newcommand{\ms}{\mathscr}
\newcommand{\minus}{\backslash}
\newcommand{\ptl}{\partial}

\newcommand{\qandq}{\quad{\text{and}}\quad}

\newcommand{\pr}{{\rm{pr}}}

\newcommand{\us}{^{\sigma}}
\newcommand{\uc}{^{\circ}}
\newcommand{\ug}{^{(g)}}
\newcommand{\ut}{^{(2)}}

\DeclareFontFamily{U}{mathx}{}
\DeclareFontShape{U}{mathx}{m}{n}{<-> mathx10}{}
\DeclareSymbolFont{mathx}{U}{mathx}{m}{n}
\DeclareMathAccent{\widehat}{0}{mathx}{"70}
\DeclareMathAccent{\widecheck}{0}{mathx}{"71}

\renewcommand{\check}{\widecheck}

\numberwithin{equation}{section}
\numberwithin{theorem}{section}

\begin{document}
\bibliographystyle{alpha} 
\title[]
{Hyperelliptic curves, minitwistors, and spacelike Zoll spaces}
\author{Nobuhiro Honda}
\address{Department of Mathematics, Institute
of Science Tokyo, 2-12-1, O-okayama, Meguro, 152-8551, JAPAN}
\email{honda@math.titech.ac.jp}

\thanks{The author was partially supported by JSPS KAKENHI Grant 22K03308.
}
\begin{abstract}

We construct a compact minitwistor space from a hyperelliptic curve with real structure and show that it yields a lot of new Lorentzian Einstein-Weyl spaces all of which are diffeomorphic to the 3-dimensional deSitter space.
These structures are real analytic, admit a circle symmetry and moreover, all their spacelike geodesics are closed and simple. 
The number of the nodes of minitwistor lines on the minitwistor space is equal to the genus of the hyperelliptic curve and is taken arbitrarily. 
These Einstein-Weyl structures deform as the hyperelliptic curves deform, and so have $(2g-1)$-dimensional moduli space, where $g$ is the genus of the hyperelliptic curve.
A relationship between the minitwistor spaces recently obtained by Hitchin from ALE gravitational instantons is also given for A$_{\rm odd}$-type.
\end{abstract}

\maketitle

\setcounter{tocdepth}{1}

\section{Introduction}
In this paper, we aim to uncover the intriguing differential geometry hidden in hyperelliptic curves with real structure.
The main theorems can be regarded as an extension of the results of our recent paper \cite{HN22} to the case of arbitrary genus.

Jones-Tod \cite{JT85} defined a minitwistor space as a complex surface that contains smooth rational curves with self-intersection number two. These curves are known as minitwistor lines and by a theorem of Kodaira \cite{K62}, they form a 3-dimensional manifold. Hitchin \cite{Hi82} demonstrated that this 3-manifold has a natural geometric structure, called the Einstein-Weyl structure (abbreviated as EW structure in the sequel) and conversely every 3-dimensional EW manifold locally arises in this way.
The EW structure is a pair of a conformal structure and a compatible affine connection that satisfies an Einstein equation. 

In our older work with Nakata \cite{HN11}, we showed that the space of rational curves on a complex surface admits an EW structure even if the rational curves have several ordinary nodes as their only singularities provided that the maximal family of such rational curves is 3-dimensional. We still refer to such rational curves as minitwistor lines and a complex surface having such curves as a minitwistor space. 
We define the genus of a minitwistor space as the number of nodes on its minitwistor lines. If such a minitwistor space is compact, then it is a rational surface and the complete linear system generated by minitwistor lines always defines a birational morphism onto its image \cite{HN22}. 
Hence, the minitwistor lines are basically hyperplane sections.
If $g$ denotes the genus of the minitwistor space, then the hyperplanes cutting out minitwistor lines must be tangent to the minitwistor space at precisely $g$ points. The closure of the space of these hyperplanes is an algebraic variety and is called the Severi variety of a projective surface in a general context. 
So the result in \cite{HN11} can be roughly stated that the Severi variety of nodal rational curves has EW structure on its regular locus if the variety is 3-dimensional.

From the general construction of the affine connection of the EW structure, every complex geodesic on it is a section of the complex EW space by a codimension-two linear subspace.
As a result, all complex geodesics on the complex EW space are essentially algebraic curves. 
So we can expect the closedness of geodesics on the real EW space formed by real minitwistor lines when the minitwistor space is equipped with a real structure.
With this background, in \cite{HN22}, we studied Lorentzian EW spaces that are associated to compact minitwistor spaces of genus one equipped with a real structure and obtained many nice properties of the EW spaces. 
In particular, we showed that the Lorentzian EW spaces are diffeomorphic to the deSitter space $\mathbb S^2\times I$ where $I$ is an open interval, their EW structures are real analytic, and that they admit an $\mathbb S^1$-action.
Moreover, we proved that all spacelike geodesics on them are simple and closed. 

In this article, extending these results, we construct compact minitwistor spaces of any given genus and investigate the associated EW spaces. The main results are as follows.


\begin{theorem}\label{t:main}
For any $g>1$, there exist compact minitwistor spaces of genus $g$ equipped with real structures
whose real EW spaces satisfy the following properties:
\begin{itemize}
\setlength{\itemsep}{-.5mm}
\item[\em (i)]
the EW spaces are diffeomorphic to $\mathbb S^2\times I$ and
the EW structures on them are Lorentzian and real analytic (Proposition \ref{p:W0} and Theorem \ref{t:extend}),
\item[\em (ii)]
the identity component of the automorphism group of each EW structure is the circle,
\item[\em(iii)]
all spacelike geodesics on them are closed and simple (Theorem \ref{t:Zoll}),
\item[\em(iv)]
the moduli spaces of the EW structures are of $(2g-1)$-dimensional (Corollary \ref{c:mod}).
\end{itemize}
\end{theorem}

From  (ii) of the theorem, none of the EW structures are isomorphic to the standard deSitter structure, and from (iv), since different points of the moduli space correspond to strictly different EW structures, they are not included in the ones obtained in \cite{HN22} because they have a 1-dimensional moduli space.

In \cite{LM09}, LeBrun and Mason proved that there is a one-to-one correspondence between orientation-reversing diffeomorphisms from $\CP_1$ to $\CP_1$ and the EW spaces that satisfy some geometric conditions. 
It is not clear as in the case $g=1$ whether the EW spaces in Theorem \ref{t:main} satisfy their geometric conditions.

We mention that the minitwistor spaces in Theorem \ref{t:main} are (rather accidentally) isomorphic to the ones obtained from Joyce's self-dual metrics with torus action \cite{J95} by a 1-dimensional reduction, which were first constructed in \cite{Hon10}. Furthermore, it turns out that their resolutions are isomorphic to some of the minitwistor spaces recently obtained by Hitchin \cite{Hi24} from ALE gravitational instantons.
On the other hand, the minitwistor lines we will consider are different from the images of the real twistor lines on the 3-dimensional twistor spaces under the quotient map. Of course, this is because the EW structures obtained from the above two structures are positive definite, while the present EW structures are indefinite.
However, the minitwistor lines for all cases belong to the common linear system on the minitwistor spaces. This means that the space of real minitwistor lines has many connected components in the linear system and different choices of the components give definite/indefinite EW structures.
We will prove that there exist precisely $2^{g-1}(g+1)$ connected components that yield the Lorentzian EW spaces and all of them satisfy the properties as in Theorem \ref{t:main}.

We briefly describe the proof of Theorem \ref{t:main}.
In Section \ref{s:Jac}, we first realize a hyperelliptic curve $\Sigma$ of genus $g$ equipped with a real structure in the cone ${\rm C}(\Lmd)$ over a non-degenerate rational normal curve $\Lmd$ in $\CP_{g+1}$.
The real locus of $\Sigma$ will consist of $(g+1)$ circles.
Utilizing the Jacobian variety and the Abel-Jacobi map for $\Sigma$, we shall classify all real hyperplanes in $\CP_{g+2}$ which are tangent to $\Sigma$ at $g$ points belonging to each of these real circles except one component and which intersect $\Sigma$ 
transversally at two points on the exceptional component.
We show that these real hyperplanes are parameterized by $2^g$ real smooth surfaces in the Jacobian variety, each of which is diffeomorphic to a disk (Proposition \ref{p:Sei1}). 
These disks constitute $2^{g-1}$ pairs in terms of their adjacent property in the Jacobian variety.
From their origin, we will call the pairs Seifert surfaces.

In Section \ref{s:mt}, we construct the compact minitwistor spaces of genus $g$ as the double covers of the cone ${\rm C}(\Lmd)$ branched along the hyperelliptic curve $\Sigma$.
The minitwistor lines will be the inverse images of the hyperplane sections of the cone classified in Section \ref{s:Jac} and have real nodes over the tangent points.
From the intersection property with the exceptional component, the minitwistor lines have real circles.
Moving these minitwistor lines by an $\mathbb S^1$-action on the minitwistor space, we obtain 3-dimensional families of minitwistor lines.
The parameter spaces of these families are all diffeomorphic to $\mathbb S^2\times I$ minus two intervals of the form $\{\text{a point}\}\times I$.
The circle action in Theorem \ref{t:main} has the union of these intervals as the rotation axis.
These punctured 3-manifolds admit Lorentzian EW structures
(Proposition \ref{p:W0}).
We see in Sections \ref{ss:imtl} and \ref{ss:cplt} that these two intervals naturally parameterize reducible hyperplane sections which consist of one irreducible conic and $2g$ straight lines. 
As in \cite{HN22}, we call these highly reducible curves irregular minitwistor lines.

In Section \ref{s:EW}, we show that the EW structures on the punctured 3-manifolds extend smoothly across the rotation axis, thereby obtaining the EW structures on the whole of $\mathbb S^2\times I$ (Theorem \ref{t:extend}). We will show it by introducing a kind of modification to an open subset of the minitwistor space and showing that all minitwistor lines contained in the subset, including the irregular ones, will be transformed into ordinary smooth minitwistor lines. So the modification transforms the open subset into the minitwistor space in the original sense. This modification is a natural extension of the one we gave in \cite{HN22} in the case $g=1$.
In particular, the transformed minitwistor spaces give examples of non-compactifiable complex surfaces (Proposition \ref{p:nc}).
Next, by using the deformation theory of singular curves, we investigate spacelike geodesics on the EW space $\mathbb S^2\times I$ and prove that all the geodesics are closed and simple (Theorem \ref{t:Zoll}).

In Section \ref{s:mod}, we investigate the automorphism groups and the moduli spaces of the EW spaces. 
Since each EW space is the parameter space of minitwistor lines, there is a natural map from the automorphism group and the moduli space of the minitwistor spaces to those of EW structures.
But because the general construction of the minitwistor space from a real EW space is not established in the global setting, and also due to the presence of the reducible minitwistor lines, we cannot immediately conclude that these mappings are bijective.
As in the case $g=1$, we show this is indeed the case in two steps.
The first one is the invariance of the punctured space under any automorphisms of the EW spaces.
We show it in Section \ref{ss:invI} by using 
the modification in Section \ref{ss:modify} and spacelike geodesics of a special kind.

The second step is to select a suitable Zariski open subset from the minitwistor space and show that the punctured EW space can be regarded as the parameter space of closed disks in the open subset, whose boundaries belong to a real sphere.
These are done in Sections \ref{ss:limtl} and \ref{ss:disk}.
For this purpose, we make use of the boundary of the EW space which consists of two spheres. We find hyperplane sections of the minitwistor space that correspond to points of these spheres through a realization of the minitwistor space as a ramified coverings of degree $(g+1)$ over quadric surfaces (Propositions \ref{p:cov1} and \ref{p:limtl}).
Using these, in Section \ref{ss:autmod}, we show that the identity component of the automorphism group of the EW space is the circle, and that if the hyperelliptic curve has no non-trivial automorphism, then the whole automorphism group is generated by the circle and the involution induced by the covering transformation of the minitwistor space over the cone (Theorem \ref{t:aut}).
The involution is a reflection with respect to an annulus $\mathbb S^1\times I$ in $\mathbb S^2\times I$.
Finally, we prove that the EW structures strictly depend on the complex structure of the initial hyperelliptic curves (Theorem \ref{t:mod}).

%
%
%

In the Appendix, we will show that including the real structure, the present minitwistor spaces are isomorphic to the minitwistor spaces that are recently obtained by Hitchin \cite{Hi24} from the twistor spaces of ALE gravitational instantons of type A$_{2k-1}$ that admit the scalar circle action, and also that the minitwistor lines in the ALE case are linearly equivalent to the ones in the Lorentzian case.

\section{Hyperelliptic curves and their tangent hyperplanes
}\label{s:Jac}
\subsection{The Jacobian varieties of hyperelliptic curves with real structures}\label{ss:J}
In this subsection, we first investigate real forms of the
Jacobian varieties of 
hyperelliptic curves equipped with a real structure
and also 2-torsion subgroups of the Jacobian varieties. 
Next, we prepare some properties of hyperelliptic curves that will be needed in the rest of this paper.

Let $g> 1$ be any integer, 
$
\lmd_1,\lmd_1',
\lmd_2,\lmd_2',
\dots,
\lmd_{g+1},\lmd_{g+1}'
$
be any distinct $(2g+2)$ real numbers arranging in this order,
and let $\Sigma$ be the hyperelliptic curve whose branches are these  points.
We denote
$$
\pi:\Sigma\lras \CP_1\qandq\tau:\Sigma\lras \Sigma,
$$
for the double covering map and the hyperelliptic involution of $\Sigma$ respectively.
We denote $r_i$ and $r'_i$ for the ramification points over $\lmd_i$ and $\lmd'_i$ respectively.
For the canonical line bundle and the canonical divisors on $\Sigma$, we have
\begin{align}\label{can1}
K_{\Sigma}\simeq \pi^*\ms O(g-1) \qandq
H^0(K_{\Sigma}) \simeq \pi^*H^0\big(\ms O(g-1)\big).
\end{align}

Let $z$ be a non-homogeneous coordinate on $\CP_1$,
and $w$ a fiber coordinate on the line bundle $\ms O(g+1)$ over $\CP_1$ that is valid on $\CC$ on which the coordinate $z$ makes sense.
Then the hyperelliptic curve $\Sigma$ can be realized in the total space of $\ms O(g+1)$ by the equation
\begin{align}\label{Sigma}
w^2 = - 
\prod_{i=1}^{g+1}\big(z-\lmd_i\big)\big(z-\lmd'_i\big).
\end{align}
Let $\Lmd\subset\CP_{g+1}$ be a rational normal curve of degree $(g+1)$ and ${\rm C}(\Lmd)\subset\CP_{g+2}$ the cone over $\Lmd$.
Then thinking the coordinate $w$ as the one in the direction of generating lines of the cone ${\rm C}(\Lmd)$, the hyperelliptic curve $\Sigma$ in \eqref{Sigma} may be regarded as a curve embedded in the cone as a double covering over $\Lmd$. 
This realization of $\Sigma\subset {\rm C}(\Lmd)$ is induced by the complete linear system $|\sum_{i=1}^{g+1}(r_i+r'_i)|$ and 
will be often used later.
Note that this means that $\Sigma$ is a non-degenerate curve of degree $(2g+2)$ in $\CP_{g+2}$.

Under the realization \eqref{Sigma} of $\Sigma$, as a basis of $H^0(K_{\Sigma})\simeq\CC^g$, we can take
\begin{align}\label{basis1}
\vec{\omega}:=\Big(\frac{dz}{w},\,\,z\frac{dz}{w},\,\,z^2\frac{dz}{w},\dots, 
z^{g-1}\frac{dz}{w}\Big).
\end{align}
So $\tau$ acts on $H^0(K_{\Sigma})\simeq\CC^g$ as $(-\id)$.

\begin{figure}
\includegraphics[height=45mm]{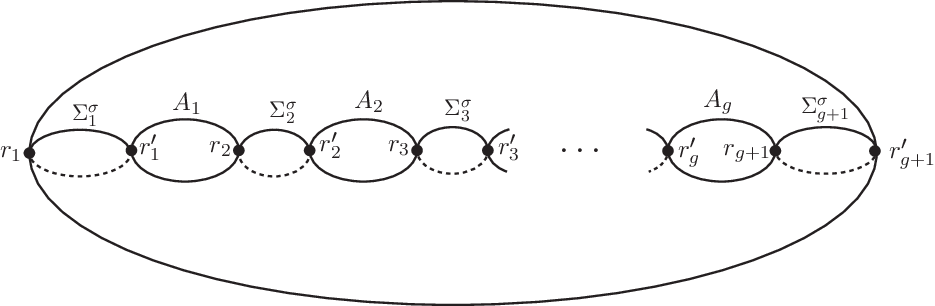}
\caption{The 1-cycles in $\Sigma$ 
}
\label{fig:Rsurface}
\end{figure}

Since all branch points of $\pi$ are real, the real structure $z\longmapsto\ol z$ lifts to the curve $\Sigma$ as $(z,w)\longmapsto (\ol z,\ol w)$ in the above coordinates. We denote $\sigma:\Sigma\lras \Sigma$ for this real structure.
If we define disjoint open intervals in $\RP_1=\RR\cup\{\infty\}$ by 
\begin{align}\label{intervals1}
K_i= 
\big(\lmd_i,\lmd'_i\big),\quad 1\le i\le g+1,
\end{align}
then the right-hand side of \eqref{Sigma} is positive on every $K_i$, and denoting $\ol K_i$ for the closed interval $[\lmd_i,\lmd'_i]$ in $\RR$,
the real locus $\Sigma^{\sigma}$ of $\Sigma$ is a double cover of
the disjoint intervals $\ol K_1\sqcup \ol K_2\sqcup\dots \sqcup \ol K_{g+1}$ branched at the ends $\lmd_i$ and $\lmd'_{i}$ of $\ol K_i$.
Let $\Sigma^{\sigma}_i\simeq \mathbb S^1$ be the connected component of $\Sigma^{\sigma}$ lying over $\ol K_i$. 
Every component of the vector-valued 1-form $\vec\omega$ in \eqref{basis1} is real-valued on these circles.
Further, for $1\le i\le g$, let $A_i$ be the cycle on $\Sigma$ defined by
\begin{align*}
A_i:=\pi\inv\big([\lmd'_i,\lmd_{i+1}]\big)\simeq \mathbb S^1.
\end{align*}
On every $A_i$, all 1-forms in \eqref{basis1} are pure imaginary. 
See Figure \ref{fig:Rsurface} for the circles $\Sigma\us_i$ and $A_i$.
The involution $\tau$ is the rotation with angle $\pi$ around the line passing through all branch points $r_1,r'_1,\dots,r_{g+1},r'_{g+1}$, and the real structure $\sigma$ is the reflection with respect to the vertical plane containing the real circles $\Sigma\us_1,\dots,\Sigma\us_{g+1}$.

The 1-cycles $\Sigma_i^{\sigma}$ and $A_i$ $(1\le i\le g)$ in $\Sigma$ constitute a basis of the homology group $H_1(\Sigma,\ZZ)\simeq\ZZ^{2g}$.
So as periods of the hyperelliptic curve $\Sigma$, we may choose
\begin{align}\label{per1}
\vec\Pi_i:=\int_{\Sigma^{\sigma}_i}\vec{\omega},\quad
\vec\Pi_{g+i}:=\int_{A_i}\vec{\omega},\quad
1\le i\le g,
\end{align}
where we are putting any one of the two possible orientations on these $2g$ circles.
We then have $\vec\Pi_i\in \RR^g$ and $\vec\Pi_{g+i}\in (\sqrt{-1}\RR)^g$ for any $1\le i\le g$.
The period lattice $\Pi$ of $\Sigma$ is defined as the lattice in $\CC^g$ generated by the periods $\vec{\Pi}_1,\dots,\vec{\Pi}_{2g}$,
and the Jacobian variety ${\rm J}_{\Sigma}$ of $\Sigma$ is defined by 
$$
{\rm J}_{\Sigma} = \CC^g/ \Pi.
$$
Since $\vec\Pi_i\in\RR^g$ and 
$\vec\Pi_{g+i}\in (\sqrt{-1}\RR)^g$ as above,
there is a decomposition 
$$
{\rm J}_{\Sigma} = \big(\RR^g/\langle\vec\Pi_1,\dots,\vec\Pi_g\rangle\big) 
\oplus \big((\sqrt{-1}\RR)^g/\langle\vec\Pi_{g+1},\dots,\vec\Pi_{2g}\rangle\big). 
$$
We call the former and the latter factors the real part and the imaginary part of ${\rm J}_{\Sigma}$ respectively.
We also note that, thanks to the hyperelliptic involution $\tau$, 
the half periods of $\Sigma$ are given by
\begin{align}\label{half1}
\frac 12 \vec\Pi_i=\int_{r_i}^{r'_i}\vec\omega, \quad
\frac 12 \vec\Pi_{g+i}=\int_{r'_i}^{r_{i+1}}\vec\omega, \quad 1\le i\le g,
\end{align}
where the integrations are done along $\Sigma\us_i$ and $A_i$ respectively in the same directions as those in \eqref{per1}.

Let $\mf a:\Sigma\lras {\rm J}_{\Sigma}$ be the Abel-Jacobi map whose base point is $r_1$. So $\mf a(r_1) = \mf o$, where $\mf o$ denotes the origin of ${\rm J}_{\Sigma}$.
Since $\tau(r_1) = r_1$ and $\tau^*\omega = -\omega$ for any differential $\omega\in H^0(K_{\Sigma})$, we have, for any point $p\in \Sigma$,
\begin{align}\label{sym0}
\mf a\big(\tau(p)\big) = -\mf a(p).
\end{align}
As usual, we extend the mapping $\mf a$ to the one from the symmetric product $\Sigma\upd$ to ${\rm J}_{\Sigma}$ for any $d>0$. 
Abel's theorem states that for any $d>0$
and any $D_1,D_2\in \Sigma^{(d)}$, the divisors $D_1$ and $D_2$
are linearly equivalent if and only if $\mf a(D_1) = \mf a(D_2)$, and the Jacobi inversion states that 
the map $\mf a:\Sigma^{(g)}\lras  {\rm J}_{\Sigma}$ is surjective.
Sometimes we write $\mf a^{(d)}$ when we want to indicate the domain of $\mf a$.
We denote
$$
W_d:=\mf a\big(\Sigma^{(d)}\big).
$$ 
This is a $d$-dimensional subvariety of $ {\rm J}_{\Sigma}$.
The symmetric product $\Sigma\upd$ is a non-singular variety and the natural mapping 
$\Sigma^d\lras \Sigma\upd$ from the $d$-fold product and the Abel-Jacobi mapping $\mf a\upd:\Sigma\upd \lras {\rm J}_{\Sigma}$ are morphisms of varieties. 
From \eqref{sym0}, 
the subvariety $W_d$ is symmetric for any $d$. Namely,
$W_d = -W_d.$
By abuse of notation, we use $W_g$ to denote the Jacobian itself.

From \eqref{half1}, we readily obtain that, for any $1\le i\le g$,
\begin{align}\label{arr}
\mf a(r_i) \equiv \frac12\sum_{j=1}^{i-1}\vec\Pi_j + \frac12\sum_{j=1}^{i-1}\vec\Pi_{g+j}\qandq
\mf a(r'_i) \equiv \frac12\sum_{j=1}^{i}\vec\Pi_j + \frac12\sum_{j=1}^{i-1}\vec\Pi_{g+j} \mod\Pi.
\end{align}

The real structure on ${\rm J}_{\Sigma}$ induced by $\sigma$ is just the complex conjugation on $\CC^g/\Pi$.
Hence, the real locus $({\rm J}_{\Sigma})^{\sigma}$ of ${\rm J}_{\Sigma}$ consists of $2^g$ number of real $g$-dimensional tori which are concretely given by 
\begin{align}\label{ccJ}
\big(\RR^g/\langle\vec\Pi_1,\dots,\vec\Pi_g\rangle\big)  + \frac{\ddd_1\vec\Pi_{g+1} + \dots + \ddd_g\vec\Pi_{2g}}2,
\quad
\ddd_1,\dots,\ddd_g\in\{0,1\}.
\end{align}
We denote $({\rm J}_{\Sigma})\us_{\mf o}$ for the identity component $\RR^g/\langle\vec\Pi_1,\dots,\vec\Pi_g\rangle$. 
Translations by elements of ${\rm J}_{\Sigma}[2]$, the subgroup of 2-torsion elements, preserve $({\rm J}_{\Sigma})\us$ but not each of the components \eqref{ccJ}.
We denote ${\rm J}_{\Sigma}[2]_{\mf o}$ for the subgroup of ${\rm J}_{\Sigma}[2]$ consisting of elements which preserve each component of $({\rm J}_{\Sigma})\us$ under translations.
From \eqref{ccJ}, we have
\begin{align}\label{ccJ2}
{\rm J}_{\Sigma}[2]_{\mf o}={\rm J}_{\Sigma}[2]\cap ({\rm J}_{\Sigma})\us_{\mf o}=
\Big\langle \frac 12\vec\Pi_1,\dots,\frac 12\vec\Pi_g\Big\rangle
\simeq(\ZZ/2\ZZ)^g.
\end{align}

Also, since we are taking the real point $r_1$ as a base point for the Abel-Jacobi map, for any $d>0$, the mapping $\mf a\upd:\Sigma\upd\lras {\rm J}_{\Sigma}$ commutes with the real structure and the subvariety $W_d$ is invariant under the real structures.

In the rest of this subsection, we give three propositions that hold for any hyperelliptic curve $\Sigma$ without assuming the existence of the real structure. 
For the first one, if an effective divisor $D$ on a hyperelliptic curve $\Sigma$ of genus $g$ contains a divisor of 
the form $p + \tau(p)$ for some point $p\in \Sigma$, 
then $\dim |D|>0$. Conversely,

\begin{proposition}\label{prop:he1}
If an effective divisor $D$ on $\Sigma$ with $\deg D\le g$ satisfies
 $\dim |D|>0$, then 
$D\ge p + \tau (p)$ for some point $p\in \Sigma$.
\end{proposition}

We omit a proof.
The second proposition is about Abel-Jacobi map $\mf a\uptwo:\Sigma\uptwo\lras {\rm J}_{\Sigma}$ and the subvariety $W_2$ for hyperelliptic curve $\Sigma$.

\begin{proposition}\label{prop:W2}
If $g>2$, the subvariety $W_2\subset {\rm J}_{\Sigma}$ has a singularity only at the origin
and it is isomorphic to the cyclic quotient singularity 
obtained by $\ZZ/(g-1)\ZZ$-action on $\CC^2$ generated by
$(z,w)\longmapsto (\zeta z,\zeta w),\,\zeta=e^{2\pi i/(g-1)}$.
\end{proposition}

Again we omit a proof. See \cite[Chapter 2, Section 7]{GrHr} or \cite[Lecture III]{Mum74}.
Since the symmetric product $\Sigma\uptwo$ is smooth
and the exceptional locus of $\mf a\uptwo:\Sigma\uptwo\lras W_2$ is the rational curve
$\big\{p+\tau(p)\in \Sigma^{(2)}\set p\in \Sigma\big\}\simeq\CP_1$, the proposition means that the map $\mf a\uptwo$ 
is exactly the minimal resolution of the quotient singularity of $W_2$ at the origin (if $g>2$).
We will use this description in Section \ref{ss:cplt} when we discuss minitwistor lines of a particular kind.

As the third proposition, we prove a property about linear systems on the hyperelliptic curve $\Sigma$ (of genus $g$) that will be important for our purpose.
To state it, we temporarily reset the notation to the effect that instead of $r_1,r'_1,\dots,r_{g+1},r'_{g+1}$ for the ramification points of $\pi:\Sigma\lras \CP_1$, we use $r_1,r_2,\dots,r_{2g+2}$.
Then we have the linear equivalence
\begin{align}\label{totalsum}
\sum_{i=1}^{2g+2} r_i \sim 2(g+1)r_1,
\end{align}
as is readily seen by thinking of the hyperelliptic curve $\Sigma$ as embedded in the cone ${\rm C}(\Lmd)$ and taking the osculating hyperplane to the rational normal curve $\Lmd$ at the branch point $\lmd_1=\pi(r_1)$.
Further, for any $(g+1)$ indices 
$1\le i_1,\dots,i_{g+1}\le 2g+2$ that are {\em mutually distinct}, 
we define an effective divisor of degree $(g+1)$ on $\Sigma$ by
\begin{align}\label{Lab}
D(i_1,\dots, i_{g+1}):=r_{i_1}+ \dots+ r_{i_{g+1}}.
\end{align}

\begin{proposition}\label{prop:L12}
The complete linear system $|D(i_1, \dots, i_{g+1})|$ is always a base point free pencil. Further, if $j_1,\dots,j_{g+1}$ are $(g+1)$ indices that are also mutually distinct,
then we have the linear equivalence
$D({i_1}, \dots, {i_{g+1}})\sim D({j_1}, \dots, {j_{g+1}})$ if and only if $\{i_1,\dots, i_{g+1}\} = \{j_1,\dots, j_{g+1}\}$
or $\{i_1,\dots, i_{g+1}\} \cap \{j_1,\dots, j_{g+1}\}=\emptyset$
(namely, complimentary to each other).
\end{proposition}

\proof
Put $F:= [D({i_1}, \dots, {i_{g+1}})]$ for simplicity.
By Riemann-Roch, $h^0(F) - h^0(K_{\Sigma}-F) = 2$. 
If $H^0(K_{\Sigma}-F)\neq 0$, then taking $2\sum_{k=1}^{g-1} r_{i_k}$ as a canonical divisor on $\Sigma$ from \eqref{can1}, we obtain that the linear system 
$|(r_{i_1}+ \dots+ r_{i_{g-1}}) - (r_{i_g} + r_{i_{g+1}})|$ is not empty.
But by Proposition \ref{prop:he1}, the linear system $|r_{i_1}+ \dots+ r_{i_{g-1}}|$ already consists of a single member and therefore $|(r_{i_1}+ \dots+ r_{i_{g-1}}) - (r_{i_g} + r_{i_{g+1}})|=\emptyset$.
Therefore, $H^0(K_{\Sigma}-F)=0$. So $h^0(F) = 2$.
If the pencil $|F|$ would have a base point, then we obtain a sub-divisor $D'$ of $D({i_1}, \dots, {i_{g+1}})$ whose degree is at most $g$ but which satisfies $\dim |D'|>0$.
From \eqref{Lab}, this again contradicts Proposition \ref{prop:he1}. Therefore $\Bs\,|F|=\emptyset$.

For the latter assertion, put
$D =  D({i_1}, \dots, {i_{g+1}})$ and 
$E =  D({j_1}, \dots, {j_{g+1}})$.
By Abel's theorem, $D\sim E$ if and only if 
$\mf a(D) =\mf a(E) $.
Suppose $\{i_1,\dots, i_{g+1}\} \cap \{j_1,\dots, j_{g+1}\}=\emptyset$.
Then from $D+E=
\sum_{i=1}^{2g+2} r_i$ and  $\sum_{i=1}^{2g+2} r_i\sim 2(g+1)r_1$ 
as in \eqref{totalsum}, $D+E\sim 2(g+1)r_1$.
Hence $\mf a(D) + \mf a(E) = \mf o$.
So $\mf a(D) =- \mf a(E)$. 
Further, since $2E = 2 r_{j_1} + \dots + 2 r_{j_{g+1}}$,
$\mf a(2E) = \mf o$.
Therefore, $\mf a(E) = -\mf a (E)$, and we obtain $\mf a(D) = \mf a (E)$. 
Conversely, suppose $\mf a(D) =\mf a(E)$.
If $(\Supp D)\cap (\Supp E)\neq \emptyset$, then we may assume $i_{g+1}= j_{g+1}$. From the relation
$D \sim  E$, this means $r_{i_1} + \dots +  r_{i_{g}} \sim  r_{j_1} + \dots +  r_{j_{g}}$. By Proposition \ref{prop:he1}, this means $r_{i_1} + \dots +  r_{i_{g}} =  r_{j_1} + \dots +  r_{j_{g}}$.
Hence $D=E$.
If $(\Supp D)\cap (\Supp E)= \emptyset$, then $D$ and $E$ are complementary to each other in the set $\{r_1,r_2,\dots, r_{2g+2}\}$, as required.
\proofend

\subsection{Some subsets in the real locus of the Jacobian variety}\label{ss:ssJ}
As in the previous subsection, let $\Sigma$ be the hyperelliptic curve of genus $g$ equipped with the real structure $\sigma$.
So far, we have treated the real circles $\Sigma\us_1,\dots,\Sigma\us_{g+1}$ equally, but in the rest of this paper, any component $\Sigma\us_i$
($1\le i\le g+1$) chosen among these components plays a special role. 
Without loss of generality, we may suppose that this component is $\Sigma\us_{g+1}$
by applying a real M\"obius transformation of $\Lmd\simeq\CP_1$ which is orientation preserving.

\begin{definition}{\em
For any $1\le i\le g+1$, we define subsets of the real locus $({\rm J}_{\Sigma})\us$ by
$$
W_{i,i}:=
\big\{ \mf a (\xi + \eta) \set \xi,\eta\in \Sigma^{\sigma}_i\big\} 
$$
and 
$$
W_{1,2,\dots,g}=\Big\{\sum_{i=1}^g\mf a(p_i)
\,\Big|\, p_1\in \Sigma\us_1, \dots,
p_g\in \Sigma\us_g \Big\}.
$$
}
\end{definition}

It is easy to see that these are closed and connected subsets of $({\rm J}_{\Sigma})\us$.
Further, using \eqref{arr}, we see that $W_{i,i}$ is included in the identity component $({\rm J}_{\Sigma})\us_{\mf o}$ and 
$W_{1,2,\dots,g}$ is contained in the following particular connected component of $({\rm J}_{\Sigma})\us$:
\begin{align}\label{compo1}
({\rm J}_{\Sigma})\us_{\mf o} + \frac12 \sum_{i=1}^{\frac g2} \vec\Pi_{g+2i-1}  & \quad {\text{if $g$ is even,}}
\end{align}
or
\begin{align}\label{compo2}
({\rm J}_{\Sigma})\us_{\mf o} + \frac12 \sum_{i=1}^{\frac{g-1}2} \vec\Pi_{g+2i}  & \quad {\text{if $g$ is odd.}}
\end{align}
We denote $({\rm J}_{\Sigma})\us_1$ for this connected component of $({\rm J}_{\Sigma})\us_{\mf o}$. 

For any integer $k>0$, we write $T^k$ for a $k$-dimensional real torus.

\begin{proposition}\label{p:123}
The subset $W_{1,2,\dots,g}$ is diffeomorphic to $\Sigma\us_1\times\dots\times \Sigma\us_g\simeq T^g$ and $W_{1,2,\dots,g} = ({\rm J}_{\Sigma})\us_1$ holds.
\end{proposition}

\proof 
Since $\Sigma\us_i\cap \Sigma\us_j=\emptyset$ if $i\neq j$, the natural mapping from $\Sigma\times \dots \times \Sigma$ ($g$-fold product) to the symmetric product $\Sigma^{(g)}$ gives an embedding from $\Sigma\us_1\times \dots\times \Sigma\us_g\simeq T^g$ to $\Sigma\upg$. In the following, we identify $\Sigma\us_1\times\dots\times \Sigma\us_g$ with the image of this embedding.
Then the subset $W_{1,2,\dots,g}$ is exactly the image under $\mf a^{(g)}$ of $\Sigma\us_1\times\dots\times \Sigma\us_g$.
We show that the restriction of $\mf a^{(g)}$ to $\Sigma\us_1\times\dots\times \Sigma\us_g$ is a diffeomorphism onto $W_{1,2,\dots,g}$.
Suppose that $\mf a^{(g)}$ is not injective on $\Sigma\us_1\times\dots\times \Sigma\us_g$, and let $D=\sum_{i=1}^g p_i$ and $E=\sum_{i=1}^g q_i$
($p_i,q_i\in \Sigma\us_i$) be points on $\Sigma\us_1\times\dots\times \Sigma\us_g\subset \Sigma^{(g)}$ that satisfies $\mf a^{(g)}(D) = \mf a^{(g)}(E)$.
By Abel's theorem, this means a linear equivalence $D\sim E$.
Since $\Sigma\us_i\cap \Sigma\us_j=\emptyset$ if $i\neq j$, 
the divisor $D$ does not have a subdivisor of the form $\pi^*\lmd$ for $\lmd\in\CP_1$, and the same thing holds for $E$. By Proposition \ref{prop:he1}, this means the coincidence
$D=E$.
Hence, $\mf a^{(g)}$ is injective on $\Sigma\us_1\times\dots\times \Sigma\us_g\subset \Sigma^{(g)}$.

Next, we show that the mapping $\mf a^{(g)}:\Sigma^{(g)}\lras {\rm J}_{\Sigma}$ is of maximal rank at any point of $\Sigma\us_1\times\dots\times \Sigma\us_g\subset \Sigma^{(g)}$.
Let $\phi:\Sigma\lras \CP_{g-1}$ be the canonical map of $\Sigma$. Since $\Sigma$ is hyperellptic, the image $\phi(\Sigma)$ is a rational normal curve in $\CP_{g-1}$ and $\phi$ is identified with the double covering map $\pi:\Sigma\lras \CP_1$.
Take any $D=\sum_{i=1}^g p_i\in \Sigma\us_1\times\dots\times \Sigma\us_g\subset \Sigma\upg$.
From the definition of the Abel-Jacobi mapping $\mf a^{(g)}$, 
if we choose local lifts $\tilde{\phi}_i$ of $\phi$ to a map from a neighborhood of each $p_i$ to $\CC^g\minus\{0\}$,
then the Jacobian matrix of $\mf a^{(g)}$ at the point $D\in \Sigma^{(g)}$ is given by the $g\times g$ matrix (\cite[p.\,342]{GrHr})
$$
\big(\tilde\phi_1(p_1),\tilde\phi_2(p_2)\dots,\tilde\phi_g(p_g)\big),
$$
where we are thinking each $\tilde\phi_i(p_i)\in \CC^g$ as a column vector.
This matrix is invertible since the image $\phi(\Sigma)$ is a rational normal curve in $\CP_{g-1}$ as above and since any distinct $g$ points on a rational normal curve in $\CP_{g-1}$ is linearly independent and since the points $\phi(p_1),\phi(p_2)\dots,\phi(p_g)$ are indeed distinct
because $\pi(p_i)\neq\pi(p_j)$ if $i\neq j$ as $\ol K_i\cap \ol K_j =\emptyset$ if $i\neq j$.
Therefore, $\mf a^{(g)}$ is of maximal rank at any point of $\Sigma\us_1\times\dots\times \Sigma\us_g$.
Hence, the image $\mf a\upg(\Sigma\us_1\times\dots\times \Sigma\us_g) = W_{1,2,\dots,g}$ is a (real) submanifold of ${\rm J}_{\Sigma}$ and it is diffeomorphic to $\Sigma\us_1\times\dots\times \Sigma\us_g\simeq T^g$.
Since $W_{1,2,\dots,g}$ is included in the connected component $({\rm J}_{\Sigma})\us_1\simeq T^g$ which is compact, we obtain the coincidence $W_{1,2,\dots,g}=({\rm J}_{\Sigma})\us_1\simeq T^g$.
\proofend

\medskip 
In the following, by $\mf t$, we mean the homomorphism from ${\rm J}_{\Sigma}$ to itself defined by $\mf a(D)\longmapsto \mf a(2D) $ where $\deg D=g$.
This mapping $\mf t:{\rm J}_{\Sigma}\lras {\rm J}_{\Sigma}$ can be identified with the quotient map for the action of ${\rm J}_{\Sigma}[2]$ on ${\rm J}_{\Sigma}$ by translations, and commutes with the real structure. 
The following property follows immediately from \eqref{ccJ} and \eqref{ccJ2}, and will be used to classify hyperplanes in $\CP_{g+2}$ that are tangent to the hyperelliptic curve $\Sigma$ in a certain proper way.

\begin{proposition}\label{prop:t1}
We have the inclusion 
$\mf t\big( ({\rm J}_{\Sigma})\us) \subset ({\rm J}_{\Sigma})\us_{\mf o}.$
Further, the restriction of $\mf t$ onto each connected component of $({\rm J}_{\Sigma})\us$ is identified with the quotient map of the component by the action of the subgroup ${\rm J}_{\Sigma}[2]_{\mf o}$ in \eqref{ccJ2} by translations.
\end{proposition}


\medskip
Composing $\mf a:\Sigma\lras {\rm J}_{\Sigma}$ with $\mf t$, we obtain a holomorphic mapping $\mf t\circ \mf a:\Sigma\lras {\rm J}_{\Sigma}$. We often denote $2\mf a$ for this mapping.

\begin{proposition}\label{prop:bdry}
For any $1\le i\le g+1$, the restriction of the mapping $2\mf a=\mf t\circ \mf a:\Sigma\lras {\rm J}_{\Sigma}$ to the circle $\Sigma\us_i\subset \Sigma$ is an immersion into $({\rm J}_{\Sigma})\us_{\mf o}\simeq T^g$, and the image $2\mf a(\Sigma\us_i)=\mf t\circ\mf a(\Sigma\us_i)$ is a real curve which has the origin $\mf o\in {\rm J}_{\Sigma}$ as the unique singular point.
Further, the curve $2\mf a(\Sigma\us_i)$ has exactly two branches at $\mf o$,
and the intersection of these branches is transversal.
\end{proposition}

\proof
The restriction $2\mf a|_{\Sigma\us_i}$ is an immersion since $\mf a:\Sigma\lras {\rm J}_{\Sigma}$ is an embedding and the mapping $\mf t:{\rm J}_{\Sigma}\lras {\rm J}_{\Sigma}$ is an unramified covering (of degree $4^g$).
Suppose that $p$ and $q$ are distinct points of $\Sigma\us_i$ which satisfy $2\mf a(p) = 2\mf a(q)$. By Abel's theorem, this means $2p\sim 2q$. 
From Proposition \ref{prop:he1}, this can happen only when $p,q\in\{r_i,r'_i\}$. 
We may suppose $p = r_i$ and $q = r'_i$, and it follows that as a map to its image, the restriction $2\mf a|_{\Sigma\us_i}$ is a map which identifies the two points $r_i$ and $r'_i$. Further, the image of these two points is the origin $\mf o$.
Then it remains to show the transversality of the intersection of the two branches of $2\mf a(\Sigma\us_i)$ at $\mf o$.
But this can be shown in a similar way to the final part of the proof of Proposition \ref{p:123}, by using $\pi(r_i)=\lmd_i\neq \pi(r'_i)=\lmd'_i$. \proofend

\medskip
Therefore, for any value of $1\le i\le g+1$, the real curve $2\mf a(\Sigma\us_i)$ is diffeomorphic to the shape of figure 8. However, according to the next proposition, this closed curve should be understood as having a winding number of two when viewed as a curve on the real surface $(W_2)\us$.

\begin{figure}
\includegraphics[width=45mm]{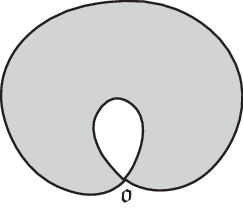}
\caption{ The region
$W_{g+1,\,g+1}$ 
}
\label{fig:Seifert1}
\end{figure}

\begin{proposition}\label{prop:bdry2}
The closed subset $W_{i,i}$ of $(W_2)\us$ itself is diffeomorphic to a closed, non-simply connected domain in $\RR^2$ bounded by a closed curve as illustrated in Figure \ref{fig:Seifert1},
and the boundary curve is naturally identified with the closed curve $2\mf a(\Sigma\us_i)$.
\end{proposition}

\proof From the definition of $W_{i,i}$ and a property of the Abel-Jacobi map $\mf a^{(2)}$ given in the previous subsection, the set $W_{i,i}$ is obtained from $\Sigma\us_i\times \Sigma\us_i\simeq T^2$ by first identifying points $(p,q)$ with $(q,p)$
($p,q\in \Sigma\us_i$ are arbitrary) and next identifying all points of the form $p+\tau (p)$ (where $p\in \Sigma\us_i$) of the symmetric product of $\Sigma\us_i$. The first identification gives a M\"obius strip whose boundary is the image of the diagonal of $\Sigma\us_i\times \Sigma\us_i$.
The boundary of the M\"obius strip exactly consists of the divisors on $\Sigma$ of the form $2p$, $p\in \Sigma\us_i$.
The second identification contracts the image of the off-diagonal in $\Sigma\us_i\times \Sigma\us_i$ into the M\"obius strip to the origin $\mf o\in {\rm J}_{\Sigma}$.
The image of the off-diagonal can be regarded as a fiber of the M\"obius strip if we view it as a non-trivial interval bundle over a circle. Then the resulting surface $W_{i,i}$ is exactly as stated in the proposition.
Moreover, the boundary is naturally identified with the curve $2\mf a(\Sigma\us_i)$ from the above description of the map $\Sigma\us_i\times \Sigma\us_i \lras W_{i,i}$.
\proofend

\medskip
Applying Proposition \ref{prop:L12} to the situation where the real structure exists,
we immediately obtain the following proposition. 
This will be important for our purpose.

\begin{proposition}\label{prop:L123}
To each $i=1,2,\dots, g+1$, take any element $\rho_i$ from $\{r_i,\,r'_i \}$.
Then the complete linear system $|\sum_{i=1}^{g+1}\rho_i|$ on $\Sigma$ is 
a fixed point free pencil.
If $\nu_i$ ($1\le i \le g+1$) are also elements of $\{r_i,\,r'_i \}$, 
then the linear equivalence $\sum_{i=1}^{g+1}\rho_i \sim \sum_{i=1}^{g+1}\nu_i$ holds
if and only if either
\begin{align}\label{cplm}
\big\{\rho_1,\dots, \rho_{g+1}\big\} = \big\{\nu_1,\dots, \nu_{g+1}\big\}
\quad{\text{
or}}\quad  
\big\{\rho_1,\dots, \rho_{g+1}\big\} \cap \big\{\nu_1,\dots, \nu_{g+1}\big\}=\emptyset
\end{align}
holds.
\end{proposition}

For any choices 
$\rho_i\in\{r_i,r'_i\}$ as in the proposition, we denote $\ms P(\rho_1,\dots,\rho_{g+1})$
for the real locus of the pencil $|\rho_{1} + \dots + \rho_{{g+1}}|$.
If $\nu_i\in\{r_i,r'_i\}$ ($1\le i\le g+1)$ are other choices, 
then by Proposition \ref{prop:L123}, the coincidence $\ms P(\rho_1,\dots,\rho_{g+1})= \ms P(\nu_1,\dots,\nu_{g+1})$ holds
if and only if they satisfy one of the conditions in \eqref{cplm}.
So we can slightly normalize the choices of $\rho_i\in\{r_i,r'_i\}$ ($1\le i\le g+1)$ by letting $r_1$ for $\rho_1$.
Therefore, the number of pencils in Proposition \ref{prop:L123} is exactly $2^{g+1}/2=2^g$.
In the rest of this paper, we will always take $\rho_1=r_1$.
Any member of the real pencil $\ms P(\rho_1,\dots,\rho_{g+1})$ is of the form
$p_1 + \dots + p_{g+1}$ with $p_i\in \Sigma\us_i$, so 
$\ms P(\rho_1,\dots,\rho_{g+1})$ can be regarded as a circle embedded in the torus $\Sigma\us_1\times\dots\times \Sigma\us_{g+1}\simeq T^{g+1}$.
We consider the restriction of the composition  
\begin{align}\label{S1}
\Sigma\us_1\times\dots\times \Sigma\us_{g+1} \stackrel{\pr}\lras 
\Sigma\us_1\times\dots\times \Sigma\us_{g} \stackrel{\mf a}\lras {\rm J}_{\Sigma} 
\end{align}
to the circle $\ms P(\rho_1,\dots,\rho_{g+1})$, where $\pr$ is the projection that just drops the last component and $\Sigma\us_1\times\dots\times \Sigma\us_{g}\simeq T^g$ is considered as a subset of the symmetric product $\Sigma^{(g)}$ as in Section \ref{ss:ssJ} so that the composition with $\mf a=\mf a^{(g)}$ makes sense.

\begin{proposition}\label{prop:diff1}
For any choices of $\rho_i$ ($1\le i\le g+1$) as in Proposition \ref{prop:L123}, this mapping (from $\ms P(\rho_1,\dots,\rho_{g+1})\simeq \mathbb S^1$ to ${\rm J}_{\Sigma}$) is a smooth embedding.
\end{proposition}

\proof
Since the restriction of $\mf a^{(g)}$ to the torus $\Sigma\us_1\times\dots \times \Sigma\us_g$ is an embedding as shown in the proof of Proposition \ref{p:123}, it suffices to show that the restriction of the map $\pr$ to the circle $\ms P(\rho_1,\dots,\rho_{g+1})$ is a diffeomorphism onto the image. 
For this, since $\ms P(\rho_1,\dots,\rho_{g+1})$ is a base point pencil, for any point $p_1\in \Sigma\us_1$, there exist unique elements $p_2\in \Sigma\us_2,\dots,p_{g+1}\in \Sigma\us_{g+1}$ that satisfy $p_1 + p_2 + \dots + p_{g+1} \in \ms P(\rho_1,\dots,\rho_{g+1})$. Hence, the circle $\ms P(\rho_1,\dots,\rho_{g+1})$ is the graph of a mapping $p_1\longmapsto (p_2,\dots,p_{g+1})$ from $\Sigma\us_1$ to $\Sigma\us_2\times\dots\times \Sigma\us_{g+1}$.
Since this mapping is obviously differentiable, the mapping $i_{g+1}:p_1\longmapsto (p_1,p_2,\dots,p_{g+1})$ from $\Sigma\us_1$ to  $\Sigma\us_1\times\dots\times \Sigma\us_{g+1}$ is an embedding.
For the same reason, the mapping $i_g:p_1\longmapsto (p_1,p_2,\dots,p_g)$ from $\Sigma\us_1$ to  $\Sigma\us_1\times\dots\times \Sigma\us_{g}$ is also an embedding.
Since $\pr\circ i_{g+1} = i_g$, we conclude that $\pr|_{\ms P(\rho_1,\dots,\rho_{g+1})}$ is a diffeomorphism onto the image.
\proofend

\subsection{Seifert surfaces in the real form of the Jacobian varieties}\label{ss:Sei}
Let $\rho_i\in\{r_i,r'_i\}$ ($1\le i\le g+1$) be as in Proposition \ref{prop:diff1}. In the following, we denote
\begin{align}\label{ms1}
\ms S^1(\rho_1,\dots,\rho_{g+1}):=(\mf a\circ\pr)\big(\ms P(\rho_1,\dots,\rho_{g+1})\big).
\end{align}
From the previous proposition, these are real smooth circles in ${\rm J}_{\Sigma}$
and from Proposition \ref{p:123}, they are all contained in the connected component $({\rm J}_{\Sigma})\us_1 = W_{1,2,\dots,g}$ of $({\rm J}_{\Sigma})\us$.
Further, any $\ms S^1(\rho_1,\dots,\rho_{g+1})$ is invariant under $(-\id)$ on ${\rm J}_{\Sigma}$ since
$p_1+\dots+p_{g+1}\sim \rho_1+\dots+\rho_{g+1}$ means
$\tau(p_1)+\dots+\tau(p_{g+1})\sim \rho_1+\dots+\rho_{g+1}$ as $\tau(\rho_i) = \rho_i$ and hence 
$-(\mf a(p_1) + \dots + \mf a(p_g))\in \ms S^1(\rho_1,\dots,\rho_{g+1})$
as $\tau$ acts as $(-\id)$ on ${\rm J}_{\Sigma}$ as in \eqref{sym0}.

\begin{proposition}
\label{p:inter3} 
Under the normalization $\rho_1=\nu_1$,
if $\ms S^1(\rho_1,\dots,\rho_{g+1})$ and $\ms S^1(\nu_1,\dots,\nu_{g+1})$ are distinct circles among the above $2^g$ ones,
$\ms S^1(\rho_1,\dots,\rho_{g+1})\cap \ms S^1(\nu_1,\dots,\nu_{g+1})\neq\emptyset$ if and only if $\nu_i=\rho_i$ for any $1< i\le g$ 
(and so $\rho_{g+1}\neq \nu_{g+1}$).
\end{proposition}

\proof If $\nu_i=\rho_i$ for any $1\le i\le g$ and $\rho_{g+1}\neq \nu_{g+1}$, then letting $\rho'_i$ be the unique element of $\{r_i,r'_i\}\minus\{\rho_i\}$, we readily have
\begin{align}\label{2pt}
\ms S^1(\rho_1,\dots,\rho_{g+1})\cap\ms S^1(\nu_1,\dots,\nu_{g+1})
= \Big\{
\sum_{i=1}^g \mf a(\rho_i), \sum_{i=1}^g \mf a(\rho'_i)
\Big\}.
\end{align}

To prove the converse, suppose
$\ms S^1(\rho_1,\dots,\rho_{g+1})\cap \ms S^1(\nu_1,\dots,\nu_{g+1})\neq\emptyset$ and assume the normalization $\rho_1=\nu_1$.
Then from the definition of these circles, there exist divisors $\sum_{i=1}^{g+1} p_i \in \ms P(\rho_1,\dots,\rho_{g+1})$
and $\sum_{i=1}^{g+1} q_i \in \ms P(\nu_1,\dots,\nu_{g+1})$ that
satisfy $\sum_{i=1}^{g} \mf a(p_i) = \sum_{i=1}^{g} \mf a(q_i)$.
By Abel's theorem, this means $\sum_{i=1}^{g} p_i \sim \sum_{i=1}^{g} q_i$.
As $\Sigma\us_i\cap \Sigma\us_j = \emptyset$ if $i\neq j$, 
using Proposition \ref{prop:he1}, this means that 
$p_i = q_i$ for any $1\le i\le g$.
As $\sum_{i=1}^{g+1} p_i\sim \sum_{i=1}^{g+1} \rho_i$ and 
$\sum_{i=1}^{g+1} q_i\sim \sum_{i=1}^{g+1} \nu_i$,
we obtain from this that 
\begin{align}\label{li15}
\sum_{i=1}^{g} p_i + p_{g+1}\sim \sum_{i=1}^{g+1} \rho_i\qandq
\sum_{i=1}^{g} p_i + q_{g+1}\sim \sum_{i=1}^{g+1} \nu_i.
\end{align}
Since both $\sum_{i=1}^{g+1} \mf a(\rho_i)$ and $\sum_{i=1}^{g+1} \mf a(\nu_i)$ belong to ${\rm J}_{\Sigma}[2]$, this readily implies $\mf a(p_{g+1} - q_{g+1})\in {\rm J}_{\Sigma}[2]$.
Hence, $2p_{g+1}\sim 2q_{g+1}$.
Again by Proposition \ref{prop:he1}, this means that 
both $p_{g+1}$ and $q_{g+1}$ belong to $\{\rho_{g+1},\rho'_{g+1}\}
\,(=\{r_{g+1},r'_{g+1}\}$).

If $p_{g+1} = q_{g+1}$ would hold, then 
$\sum_{i=1}^{g+1} p_i = \sum_{i=1}^{g+1} q_i$ and hence
$\ms S^1(\rho_1,\dots,\rho_{g+1})=\ms S^1(\nu_1,\dots,\nu_{g+1})$, contradicting our assumption.
Hence $p_{g+1} \neq q_{g+1}$ and therefore  
$(p_{g+1}, q_{g+1})$ is equal to either $(\rho_{g+1}, \rho'_{g+1})$ or 
 $(\rho'_{g+1}, \rho_{g+1})$.
If the former is the case, then from the coincidence of the last component (i.e., the point belonging to $\Sigma\us_{g+1}$) and base point freeness of the pencil $\ms P(\rho_1,\dots,\rho_{g+1})$ as in Proposition \ref{prop:L123},
the former of \eqref{li15} means $p_i = \rho_i$ for any $1\le i\le g$.
From the latter of \eqref{li15}, this means $\sum_{i=1}^g\rho_i + \rho'_{g+1} \sim \sum_{i=1}^{g+1}\nu_i$. 
By the normalization $\rho_1 = \nu_1$, this implies 
$\sum_{i=1}^g\rho_i + \rho'_{g+1} = \sum_{i=1}^{g+1}\nu_i$ and we are done. 
Similarly, if $(p_{g+1}, q_{g+1})=(\rho'_{g+1}, \rho_{g+1})$, then we obtain $\sum_{i=1}^g\rho'_i + \rho_{g+1} \sim \sum_{i=1}^{g+1}\nu_i$. 
By Proposition \ref{prop:L123}, we can replace LHS with the complementary divisor and we obtain 
$\sum_{i=1}^g\rho_i + \rho'_{g+1} \sim \sum_{i=1}^{g+1}\nu_i$.
By the normalization $\rho_1 = \nu_1$, this again implies 
$\sum_{i=1}^g\rho_i + \rho'_{g+1} = \sum_{i=1}^{g+1}\nu_i$ and we are done. 
\proofend

\medskip From this proposition, we introduce the following 

\begin{definition}
\label{def:ac}
{\em
To any circle $\ms S^1=\ms S^1(\rho_1,\dots,\rho_{g+1})$ among the $2^g$ ones in $({\rm J}_{\Sigma})\us_1=W_{1,2,\dots,g}$, we denote $\hat{\ms S}^1$ for the unique circle that intersects $\ms S^1$ as in Proposition \ref{p:inter3}, and we call $\ms S^1$ and $\hat{\ms S^1}$ the {\em intersecting pair}.
}
\end{definition}

The number of the intersecting pairs is obviously $2^{g}/2=2^{g-1}$.

We keep any choices of $\rho_i\in\{r_i,r'_i\}$, $i=1,\dots,g+1$.
The mapping $\mf t:{\rm J}_{\Sigma}\lras {\rm J}_{\Sigma}$ was defined by $x\longmapsto 2x$ for  $x\in {\rm J}_{\Sigma}$, and by Proposition \ref{prop:t1} it maps $({\rm J}_{\Sigma})\us$ to $({\rm J}_{\Sigma})\us_{\mf o}$.
In the sequel, for simplicity of notation, we write $\mf p$ to mean the image $\mf a(p)$ for $p\in {\Sigma}$.
So from \eqref{ms1}, any element of $\ms S^1(\rho_1,\dots,\rho_{g+1})$ may be written 
$\mf p_1 + \dots + \mf p_{g}$ for some points $p_1,\dots,p_g,p_{g+1}$ satisfying $\sum_{i=1}^{g+1}p_i\sim \sum_{i=1}^{g+1}\rho_i$, where $p_i\in {\Sigma}\us_i$ for any $i$.
For such a point of $\ms S^1(\rho_1,\dots,\rho_{g+1})$,
since $\mf t(\mf p_1 + \dots + \mf p_{g+1}) = 
2\mf a(\rho_1) + \dots + 2\mf a(\rho_{g+1})=\mf o$ as $\rho_i\in\{r_i,r'_i\}$, we obtain
\begin{align*}
-\mf t(\mf p_1 + \dots + \mf p_g) &= \mf t(\mf p_{g+1})\\
&=2\mf p_{g+1}.
\end{align*}
Together with the above $(-\id)$-invariance of the circle $\ms S^1(\rho_1,\dots,\rho_{g+1})
\subset {\rm J}_{\Sigma}$,
this implies 
$\mf t(\ms S^1(\rho_1,\dots,\rho_{g+1})) \subset 2\mf a({\Sigma}\us_{g+1}).$
Furthermore, using that every point $p_{g+1}\in {\Sigma}\us_{g+1}$ admits points $p_1\in {\Sigma}\us_1,\dots,p_g\in {\Sigma}\us_g$ that satisfy
$\sum_{i=1}^{g+1}p_i \in \ms P(\rho_1,\dots,\rho_{g+1})$,
we readily obtain the reverse inclusion. Hence we obtain 
\begin{align}\label{image3}
\mf t\big(\ms S^1(\rho_1,\dots,\rho_{g+1})\big) = 2\mf a({\Sigma}\us_{g+1}).
\end{align}
As before, we mean by $\rho'_i$ the unique element of $\{r_i,r'_i\}\minus \{\rho_i\}$.
Then for the two points appearing in the intersection \eqref{2pt},
\begin{align}\label{t2}
\mf t\Big(\sum_{i=1}^{g}\mf a(\rho_i)\Big) = \mf t\Big(\sum_{i=1}^{g}\mf a(\rho'_i)\Big) = \mf o.
\end{align}

\begin{proposition}\label{prop:circles0}
The surjective mapping
$$\mf t|_{\ms S^1(\rho_1,\dots,\rho_{g+1})}:\ms S^1(\rho_1,\dots,\rho_{g+1})\lras 2\mf a(\Sigma\us_{g+1})$$ is exactly the one that identifies the two points $\sum_{i=1}^g \mf a(\rho_i)$ and 
$\sum_{i=1}^g \mf a(\rho'_i)$ of the circle $\ms S^1(\rho_1,\dots,\rho_{g+1})$.
\end{proposition}

\proof By \eqref{t2},
it is enough to show that if $D=\sum_{i=1}^{g+1} p_i$ and $E=\sum_{i=1}^{g+1} q_i$ are distinct members of the pencil $ \ms P (\rho_1,\dots,\rho_{g+1})$ that satisfy
$\mf t(\sum_{i=1}^g \mf p_i) = \mf t(\sum_{i=1}^g \mf q_i)$ then 
either $D= \sum_{i=1}^{g+1} \rho_i$ and $E= \sum_{i=1}^{g+1} \rho'_i$
or $D= \sum_{i=1}^{g+1} \rho'_i$ and $E= \sum_{i=1}^{g+1} \rho_i$.
The assumptions mean $2\mf p_{g+1} = 2\mf q_{g+1}$, which implies the linear equivalence $2p_{g+1}\sim 2q_{g+1}$ by Abel's theorem.
If $p_{g+1}= q_{g+1}$, then from the base point freeness of the pencil $ \ms P (\rho_1,\dots,\rho_{g+1})$, we obtain that $\sum_{i=1}^{g+1} p_i =\sum_{i=1}^{g+1} q_i$, which contradicts $D\neq E$.
If  $p_{g+1}\neq  q_{g+1}$, then $2p_{g+1}\sim 2q_{g+1}$ and  Proposition \ref{prop:he1} mean that either $(p_{g+1},q_{g+1}) = (\rho_{g+1}, \rho'_{g+1})$ or $(p_{g+1},q_{g+1}) = (\rho'_{g+1}, \rho_{g+1})$.
Again from the base point freeness of the same pencil, the former 
(resp.\,the latter) implies 
$(p_i,q_i) = (\rho_i,\rho'_i)$ 
(resp.\,$(p_i,q_i) = (\rho'_i,\rho_i)$)
for any $1\le i\le g$.
This means that $D= \sum_{i=1}^{g+1} \rho_i$ and $E= \sum_{i=1}^{g+1} \rho'_i$
(resp.\,$D= \sum_{i=1}^{g+1} \rho'_i$ and $E= \sum_{i=1}^{g+1} \rho_i$). 
\proofend

\medskip
From Proposition \ref{prop:t1}, the mapping 
\begin{align}\label{t1}
\mf t_1:=\mf t|_{({\rm J}_{\Sigma})\us_1}:({\rm J}_{\Sigma})\us_1=W_{1,2,\dots,g}\lras ({\rm J}_{\Sigma})\us_{\mf o}
\end{align}
can be regarded as a quotient map for the translations by ${\rm J}_{\Sigma}[2]_{\mf o}$, so $\mf t_1$
realizes $({\rm J}_{\Sigma})\us_1$ as a $2^g$-fold covering over $({\rm J}_{\Sigma})\us_{\mf o}$. For the inverse image of the boundary curve $2\mf a(\Sigma\us_{g+1})$ of the region $W_{g+1,\,g+1}$ under $\mf t_1$, we readily have:

\begin{proposition}\label{prop:circles}
$\mf t_1\inv\big(2\mf a(\Sigma\us_{g+1})\big) = \displaystyle\bigcup_{\rho_i\in\{r_i,r'_i\}}
\ms S^1(\rho_1,\dots,\rho_{g+1}).
$
\end{proposition}

\proof
The inclusion `$\supset$' follows immediately from \eqref{image3}.
The reverse inclusion follows from the facts that $\mf t_1:({\rm J}_{\Sigma})\us_1\lras ({\rm J}_{\Sigma})\us_{\mf o}$ is a $2^g$-fold unramified covering, 
that $\mf t_1$ is generically one-to-one on each $\ms S^1(\rho_1,\dots,\rho_{g+1})$ by Proposition \ref{prop:circles0}, as well as the fact that there are exactly $2^g$ circles of the form $\ms S^1(\rho_1,\dots,\rho_{g+1})$.
\proofend

\medskip
As above,
there are $2^g$ circles of the form $\ms S^1(\rho_1,\dots,\rho_{g+1})$ and they are grouped into $2^{g-1}$ intersecting pairs as in  Definition \ref{def:ac}.
To simplify the notation, in the following, we choose any one from each intersecting pair and write them 
\begin{align}\label{bdry1}
\ms S^1_1,\ms S^1_2,\dots, \ms S^1_{2^{g-1}}.
\end{align}
Here no rule is put to order these circles.
In the notation of Definition \ref{def:ac}, the remaining circles are accordingly arranged as
\begin{align}\label{bdry2}
\hat{\ms S}^1_1,\hat{\ms S}^1_2,\dots, \hat{\ms S}^1_{2^{g-1}}.
\end{align}
By Proposition \ref{p:inter3}, two distinct circles among the $2g$ circles \eqref{bdry1} and \eqref{bdry2} intersect only for $\ms S^1_k$ and $\hat {\ms S}^1_k$, and the intersection 
$\ms S^1_k\cap \hat{\ms S}^1_k$ consists of the two points \eqref{2pt} for any $1\le k\le 2^{g-1}$.

Recall that the region $W_{g+1,\,g+1}$ in the real surface $(W_2)\us$ is bounded by the curve $2\mf a(\Sigma\us_{g+1})$ with a self-intersection point $\mf o$ (Propositions \ref{prop:bdry} and \ref{prop:bdry2}). We denote the interior of $W_{g+1,\,g+1}$ by 
$W^{\circ}_{g+1,\,g+1}$. Namely, 
$$
W^{\circ}_{g+1,\,g+1}= W_{g+1,\,g+1} \minus\, 2\mf a(\Sigma\us_{g+1}).
$$

\begin{proposition}
Let $\mf t_1:({\rm J}_{\Sigma})\us_1\lras ({\rm J}_{\Sigma})\us_{\mf o}$ be as in \eqref{t1}.
The inverse image $\mf t_1\inv \big(W_{g+1,\,g+1}\uc\big)\subset ({\rm J}_{\Sigma})\us_1=W_{1,2,\dots,g}$ consists of $2^g$ connected components and each component is diffeomorphic to a 2-dimensional open disk. These components are bounded by the $2^g$ circles 
$\ms S^1_k$ and $\hat{\ms S}^1_k$, $1\le k\le 2^{g-1}$.
\end{proposition}

\proof
As in \eqref{2pt}, 
the intersection $\ms S^1_k\cap \hat{\ms S}^1_k$ consists of the two points  
$\sum_{i=1}^g \mf a(\rho_i)$ and $\sum_{i=1}^g \mf a(\rho'_i)$, where $\rho_1,\dots,\rho_{g+1}$ are the points that determine the circle $\ms S_k^1$.
We denote $\mf o_k$ and ${\mf o}'_k$ for these two points respectively.
Using \eqref{arr}, the element $\mf c:=\mf a(\rho_{g+1} +\rho'_{g+1})\in {\rm J}_{\Sigma}$ belongs to the subgroup ${\rm J}_{\Sigma}[2]_{\mf o}$, so it defines a covering transformation of 
$\mf t_1:({\rm J}_{\Sigma})\us_1\lras ({\rm J}_{\Sigma})\us_{\mf o}$.
Further, using Proposition \ref{prop:L123}, 
for any $1\le k\le 2^{g-1}$,
\begin{align*}
\mf o_k + \mf c &= \sum_{i=1}^g \mf a(\rho_i) + \big(\mf a(\rho_{g+1}) + \mf a(\rho'_{g+1})\big)\\
&= \sum_{i=1}^{g+1} \mf a(\rho'_i) + \mf a(\rho'_{g+1}) \\
&=\mf o'_k.
\end{align*}
Since only $\ms S^1_k$ and $\hat{\ms S}^1_k$ pass $\mf o_k$ and $\mf o'_k$ among the $2^g$ circles by Proposition \ref{p:inter3}, 
this means that 
only $\ms S_1^k + \mf c = \ms S_1^k$ and $\ms S_1^k + \mf c = \hat{\ms S}_1^k$ can happen.
But if the former would hold, then the mapping $\mf t_1|_{\ms S_1^k}:\ms S_1^k\lras 2\mf a(\Sigma\us_{g+1})$ would be two-to-one which contradicts Proposition \ref{prop:circles0}.
Therefore, $\ms S_1^k + \mf c = \hat{\ms S}_1^k$.

\begin{figure}
\includegraphics[width=70mm]{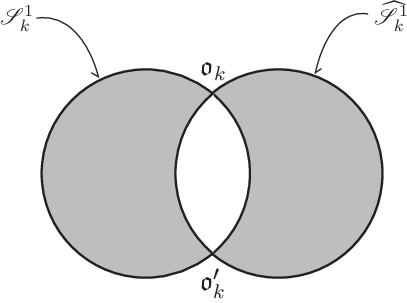}
\caption{
The Seifert surface $\mf S_k$ and its boundary $\ms S^1_k \cup \hat{\ms S}^1_k$
}
\label{fig:Seifert2}
\end{figure}

This means that each of the lifts of the open region $W_{g+1,\,g+1}\uc$ to the covering $({\rm J}_{\Sigma})_1\us$ separates the origin $\mf o$ into two points, and the two points have to be $\mf o_k$ and $\mf o'_k$ for some $1\le k\le 2^{g-1}$. (See Figure \ref{fig:Seifert2}.)
Further, the lift of $W_{g+1,\,g+1}\uc$ consists of two connected components both of which are mapped bijectively and differentiably onto the original region $W_{g+1,\,g+1}\uc$.
Clearly, each of these connected components is diffeomorphic to a disk. Further, since there are exactly $2^{k-1}$ intersection pairs of the circles, the number of the connected components of $\mf t_1\inv (W_{g+1,\,g+1}\uc)$ is $2^g$.
\proofend

\medskip
As in the proof, all disks that are connected components of $\mf t_1\inv \big(W_{g+1,\,g+1}\uc\big)$ are mapped differentiably and bijectively onto $W^{\circ}_{g+1,\,g+1}$ by  $\mf t_1$, but these maps are not homeomorphisms because the two points $\mf o_k$ and $\mf o'_k$ are identified.
Also, the proof means that the closed inverse image $\mf t_1\inv \big(W_{g+1,\,g+1}\big)$ consists of $2^{g-1}$ connected components and each component
consists of two closed `disks' that share the two points 
$\mf o_k = \sum_{i=1}^g \mf a(\rho_i)$ and $\mf o'_k = \sum_{i=1}^g \mf a(\rho'_i)$. (See Figure \ref{fig:Seifert2}.)
Associating the boundary, the connected components of $\mf t_1\inv(W_{g+1,\,g+1})$ are in one-to-one correspondence with
the set of intersecting circles $\ms S^1_k\cup\hat{\ms S}^1_k$ in the sense of Definition \ref{def:ac}.

\begin{definition}\label{def:Sei}{\em
By abuse of terminology, we call the connected components of the closed set $\mf t_1\inv \big(W_{g+1,\,g+1}\big)$ the {\em Seifert surfaces} and denote $\mf S_k$ for the component bounded by the curve $\ms S^1_k\cup\hat{\ms S}^1_k$.
We call the number $k$ the {\em boundary data} of $\mf S_k$ and the two points $\mf o_k$ and $\mf o'_k$ the {\em singularities} of $\mf S_k$. (See Figure \ref{fig:Seifert2}.)}
\proofend
\end{definition}

The boundary data is in one-to-one correspondence with the intersection pairs and it is determined from the choice of $g$ elements $\rho_1,\dots, \rho_g$ among $\{r_1,r'_1\},\dots, \{r_g,r'_g\}$ respectively. In fact, these $g$ points determine the intersecting pair $\ms S^1(\rho_1,\dots, \rho_{g}, r_{g+1})$ and $\ms S^1(\rho_1,\dots, \rho_{g}, r'_{g+1})$. 
As above, the complementary $g$ points $\rho'_1,\dots, \rho'_g$ give the same intersecting pair.
Therefore, we can always assume that $\rho_1 = r_1$ and the number $2^{g-1}$ corresponds to the choices of $\rho_2,\dots,\rho_g$.

The distinction between the singularities $\mf o_k$ and $\mf o'_k$ is unimportant.
From the definition, we have
\begin{align}\label{disj1}
\mf t_1\inv \big(W_{g+1,\,g+1}\big)
= \bigsqcup_{1\le k\le 2^{g-1}} \mf S_k 
\end{align}
and 
\begin{align}\label{disj2}
\ptl \mf S_k = \ms S_k^1\cup\hat{\ms S}^1_k, \quad 1\le  k\le 2^{g-1}.
\end{align}
We denote the interior of the Seifert surface $\mf S_k$ by $\mf S_k\uc$.
This consists of two connected components and 
\begin{align}\label{disj3}
\mf t_1\inv \big(W_{g+1,\,g+1}\uc\big)
= \bigsqcup_{1\le k\le 2^{g-1}} \mf S_k\uc.
\end{align}

The next proposition means that the Seifert surface $\mf S_k$ parameterizes real hyperplanes in $\CP_{g+2}$ that are tangent to $\Sigma$ in a certain proper way. It will be a key for describing the space of minitwistor lines.
Recall that the component $\mf S_k$ is determined from the choice of the pencil of the form $\ms P(\rho_1,\dots,\rho_{g+1})$ on $\Sigma$, where $\rho_i\in \{r_i,r'_i\}$.

\begin{proposition}\label{p:Sei1} Let $1\le k\le 2^{g-1}$ be any boundary data determining the Seifert surface $\mf S_k$ and $\mf s$ any point of $\mf S_k$.
Then there exist a unique point $(p_1,\dots,p_g)\in \Sigma\us_1\times\dots\times \Sigma\us_g$ satisfying the equation $\sum_{i=1}^g \mf a( p_i) = \mf s$ and 
a unique real hyperplane $h\subset\CP_{g+2}$ that satisfy
\begin{align}\label{hB1}
h|_{\Sigma} = \sum_{i=1}^g 2p_i + \xi + \eta
\end{align}
for some points $\xi,\eta\in \Sigma\us_{g+1}$.
Further, the tangent points $p_1,\dots,p_g$ and the two points $\xi,\eta$ satisfy the following:
\begin{itemize}
\setlength{\itemsep}{0cm}
\item[\em (i)] 
if $\mf s\not\in \ptl \mf S_k$, then none of the points $p_1,\dots, p_g$ are branch points of $\pi:\Sigma\lras \CP_1$, the divisor $\xi+\eta$ is unique, and further, $\eta\not\in\{\xi, \tau(\xi)\}$,
\item[\em (ii)] 
if $\mf s\in \ptl \mf S_k\minus\{\mf o_k,\mf o'_k\}$, then again none of the points $p_1,\dots, p_g$ are branch points of $\pi$, 
$\xi=\eta\not\in\{r_{g+1},r'_{g+1}\}$, and 
the point $\xi\,(=\eta)$ is unique,
\item[\em (iii)] 
if $\mf s=\mf o_k$ (resp.\,$\mf s = \mf o'_k$), then $p_i=\rho_i$ 
(resp.\,$p_i = \rho'_i$) for any $1\le i\le g$ and further,
$\eta=\tau(\xi)$ where $\xi\in \Sigma\us_{g+1}$ is arbitrary.
\end{itemize}
Conversely, any real hyperplane $h\subset\CP_{g+2}$ that satisfies \eqref{hB1} for some $(p_1,\dots,p_g)\in \Sigma\us_1\times\dots\times \Sigma\us_g$ and $\xi,\eta\in \Sigma\us_{g+1}$ is obtained this way.
\end{proposition}

\proof
Since $\mf S_k$ is included in $({\rm J}_{\Sigma})\us_1$ and $({\rm J}_{\Sigma})\us_1=W_{1,2,\dots,g}$ is diffeomorphic to $\Sigma\us_1\times \dots\times \Sigma\us_g$ through $\mf a$ by Proposition \ref{p:123} (ii),
any point $\mf s\in \mf S_k$ uniquely determines a point $(p_1,\dots,p_g)\in \Sigma\us_1\times \dots\times \Sigma\us_g$ that satisfies $\sum_{i=1}^g \mf a(p_i) = \mf s$.
Next, since $\mf t(\mf s)\in W_{g+1,\,g+1}$ as $\mf s\in \mf S_k$  and since any point of $W_{g+1,\,g+1}$ is written as $\mf a(\xi+\eta)$ for some points $\xi,\eta\in \Sigma\us_{g+1}$, noting also 
$W_{g+1,\,g+1}=-W_{g+1,\,g+1}$ by the hyperelliptic involution $\tau$,
we can write $\mf t(\mf s) = - \mf a(\xi+\eta)$ using such $\xi$ and $\eta$.
This is equivalent to $\mf a(\sum_{i=1}^g 2p_i + \xi+\eta) = \mf o$.
Since $\mf a(h|_{\Sigma}) = \mf o$ for the hyperplane class $h|_{\Sigma}$, by Abel's theorem, this implies the linear equivalence
$h|_{\Sigma}\sim \sum_{i=1}^g 2p_i + \xi+\eta$. 
Then the uniqueness of the hyperplane $h$ such that \eqref{hB1} holds follows immediately from the non-degeneracy of the curve $\Sigma$ in $\CP_{g+2}$.
This hyperplane $h$ is real since the divisor in the right-hand side of \eqref{hB1} is real and the embedding $\Sigma\subset\CP_{g+2}$ preserves the real structure.

To show (i) and (ii), if for some $i$ the point $p_i$ would be one of the branch points $\rho_i$ and $\rho'_i$ of $\pi$, then 
from the tangency condition, $h$ contains the generating line through $p_i$.
Again from the tangency, this means that all points $p_1,\dots, p_g$ are branch points of $\pi$. Namely, $p_i\in\{\rho_i,\rho'_i\}$ for any $i\le g$.
From the constraint $\sum_{1\le i\le g} \mf a(p_i) = \mf s\in \mf S_k$, this implies either $p_i = \rho_i$ for any $i\le g$ or $p_i=\rho'_i$ for any $i\le g$. In the former case $\mf s = \mf o_k$ and in the latter case $\mf s = \mf o'_k$. Thus, none of $p_i$ can be a branch point of $\pi$ if $\mf s\neq \mf o_k,\mf o'_k$. 
In a similar way, $\eta=\tau(\xi)$ implies $\mf s\in\{ \mf o_k, \mf o'_k\}$ and therefore $\eta\neq\tau(\xi)$ if $\mf s\neq \mf o_k,\mf o'_k$.
Moreover, if $\mf s\not\in\{\mf o_k,\mf o'_k\}$, then $\mf t(\mf s) \in W_{g+1,\,g+1}\minus \{\mf o\}$ as $\mf t:\mf S_k\to W_{g+1,\,g+1}$ is a double covering.  
Further, any element of $W_{g+1,\,g+1}\minus\{\mf o\}$ may be written $\mf a(\xi+\eta)$ with $\xi,\eta\in \Sigma\us_{g+1}$ in a unique way by Proposition \ref{prop:he1}.
Hence, we obtain (i) and (ii). 

To show (iii), if $\mf s = \mf o_k$, then $p_i = \rho_i$ for any $i\le g$ from the isomorphism $W_{1,\dots,g}\simeq \Sigma\us_1\times\dots\times \Sigma\us_g$ under $\mf a^{(g)}$, so $\mf t(\mf s) = \mf o$. From Proposition \ref{prop:he1}, this implies $\eta=\tau(\xi)$, with $\xi\in \Sigma\us_{g+1}$ being arbitrary. 
The case $\mf s = \mf o'_k$ can be shown similarly, 
and we obtain (iii).

For the converse, applying $\mf a$ to \eqref{hB1}, we obtain $2\sum_{1\le i\le g} \mf a(p_i) = - \mf a(\xi+\eta)$.
As $W_{g+1,\,g+1} = - W_{g+1,\,g+1}$, this means $\mf t(\sum_{1\le i\le g} \mf a(p_i)) \in W_{g+1,\,g+1}$.
So $ \sum_{1\le i\le g} \mf a(p_i) \in \mf S_k$ for some $1\le k\le 2^{g-1}$.
The rest is easy and we omit the detail.
\proofend

\medskip
As in the last item of the proposition, 
the singularities $\mf s=\mf o_k,\mf o'_k$ do not determine the hyperplane $h$ uniquely.
This non-uniqueness can be seen directly as follows.
Let $\lmd\in K_{g+1}$ be an arbitrary point.
Because $\Lmd$ was a rational normal curve in $\CP_{g+1}$,
the divisor $\sum_{1\le i\le g} l_i + l(\lmd)$ on the cone ${\rm C}(\Lmd)$, where $l_i$ and $l(\lmd)$ are the generating lines of the cone ${\rm C}(\Lmd)$ over the points $\pi(\rho_i)$ and $\lmd\in K_{g+1}$ respectively, is a hyperplane section of ${\rm C}(\Lmd)$ and its restriction to $\Sigma$ is of the form $\sum_{1\le i\le g} 2\rho_i + \xi + \tau(\xi)$, where $\{\xi,\tau(\xi)\} = l(\lmd)\cap \Sigma\us_{g+1}$. Namely, there is a freedom for choosing $\lmd\in K_{g+1}$.
These highly reducible hyperplane sections of ${\rm C}(\Lmd)$ will be important in the rest of this paper.

\section{The minitwistor spaces and their minitwistor lines}\label{s:mt}

\subsection{The minitwistor spaces}\label{ss:mt}
First, we construct projective algebraic surfaces which will be minitwistor spaces of genus $g$ equipped with real structures.
As in the beginning of Section \ref{s:Jac}, let $\Sigma$ be the hyperelliptic curve of genus $g$ that is realized as a double covering of $\CP_1$ branched at the real points $\lmd_1<\lmd'_1<\dots<\lmd_{g+1}<\lmd'_{g+1}$. 
This was equipped with a natural real structure $\sigma$.
Also, as before, $\Lmd\subset\CP_{g+1}$ and ${\rm C}(\Lmd)\subset\CP_{g+2}$ denote the rational normal curve of degree $(g+1)$ and the cone over $\Lmd$ respectively, and let 
${\bf v}\in\CP_{g+2}$ denote the vertex of the cone ${\rm C}(\Lmd)$.
Then from the equation \eqref{Sigma} of $\Sigma$, the hyperelliptic curve $\Sigma$ is written ${\rm C}(\Lmd)\cap Q$ where $Q\subset\CP_{g+2}$ is a quadric that does not pass through $\bf v$.

Since $Q$ is of even degree, we can take a double covering of ${\rm C}(\Lmd)$ branched along $\Sigma$, and it is realized in $\CP_{g+3}$.
We define $\ms T$ to be this double covering.
This is a projective algebraic surface that is determined from the real numbers $\{\lmd_i,\lmd'_i\set 1\le i\le g+1\}$,
and singularities of $\ms T$ are the two cone singularities over the vertex $\bf v$ of ${\rm C}(\Lmd)$.
As in the beginning of Section \ref{s:Jac}, let $z$ be a coordinate on $\CP_1\simeq\Lmd$ and $w$ the coordinate in the direction of generating lines of the cone.
Then the surface $\ms T$ is explicitly given by the real equation
\begin{align}\label{T}
y^2 = - \prod_{i=1}^{g+1}(z-\lmd_i)(z-\lmd'_i) - w^2,
\end{align}
where $y$ is a fiber coordinate of the projection $\Pi:\CP_{g+3}\lras\CP_{g+2}$ from a point.
We will use the same letter $\Pi$ to mean the double covering map $\ms T\lras {\rm C}(\Lmd)$.
$\ms T$ has a natural real structure that is concretely defined by $(z,w,y)\longmapsto (\ol z,\ol w,\ol y)$ in the above coordinates.
We denote this real structure by the same letter $\sigma$. 
On the other hand, we denote the projection $\CP_{g+2}\lras\CP_{g+1}$ from the vertex $\bf v$ or its restriction to ${\rm C}(\Lmd)$ by $\pi$.
Recall that the letter $\pi$ has been used to mean the double covering map $\Sigma\lras \Lmd$ but we will use the same letter since the projection from $\bf v$ restricts to the double covering map.

The composition $\CP_{g+3}\stackrel{\Pi}\lras\CP_{g+2}\stackrel{\pi}\lras \CP_{g+1}$ is a projection from a line, and from \eqref{T}, the restriction of $\pi\circ\Pi$ to $\ms T$ realizes $\ms T$ as a rational conic bundle over $\Lmd$.
We denote $l_{\infty}\subset\CP_{g+3}$ for the last line.
This is real.
The surface $\ms T$ intersects $l_{\infty}$ at precisely two points
and they are exactly the two singular points of $\ms T$.
It is easy to show that these two points are exchanged by $\sigma$.
So in the following, we denote the two singularities of $\ms T$
by $\bm p_{\infty}$ and $\ol{\bm p}_{\infty}$.
These are the points of indeterminacy of the rational conic bundle structure $\pi\circ\Pi:\ms T\lras \Lmd$.

By the coordinate change $u:= y+\sqrt{-1}w$ and $v:= y-\sqrt{-1}w$,
the equation \eqref{T} of $\ms T$ can be rewritten in a very simple form as
\begin{align}\label{T1}
uv = - \prod_{i=1}^{g+1}(z-\lmd_i)(z-\lmd'_i).
\end{align}
If we identify the branch divisor $\Sigma\subset {\rm C}(\Lmd)$ of $\Pi:\ms T\lras {\rm C}(\Lmd)$ and the ramification divisor on $\ms T$,
then $\Sigma\subset\ms T$ is defined by $u+v=0$ and the ramification points of $\pi:\Sigma\lras \Lmd$ are
concretely given by $r_i = \{u=v=0,z=\lmd_i\}$ and $r'_i = \{u=v=0,z=\lmd'_i\}$.
The real structure $\sigma$ is given by 
\begin{align}\label{rs}
(z,u,v)\longmapsto (\ol z, \ol v,\ol u).
\end{align}

The surface $\ms T$ is invariant under the effective $\mathbb S^1$-action defined by 
\begin{align}\label{S1act}
(z,u,v)\longmapsto \big(z,e^{\sqrt{-1}\theta}u,e^{-\sqrt{-1}\theta} v\big),
\quad e^{\sqrt{-1}\theta}\in U(1)=\mathbb S^1,
\end{align}
or equivalently, 
\begin{align}\label{S1act2}
(z,w,y)\longmapsto \big(z,\cos\theta\, w + \sin\theta\, y,\, -\sin\theta\, w + \cos\theta\, y\big),
\quad \theta\in\RR.
\end{align}
This $\mathbb S^1$-action commutes with the above real structure $\sigma$ on $\ms T$.
Although this $\mathbb S^1$-action on $\ms T$ is not a lift of any $\mathbb S^1$-action on the cone ${\rm C}(\Lmd)$,
it preserves each fiber of the rational conic bundle map $\pi\circ\Pi:\ms T\lras\Lmd$.
This can also be regarded as a quotient map by the $\CC^*$-action on $\ms T$ obtained by allowing $e^{\sqrt{-1}\theta}$ to move in $\CC^*$. 

The real locus of $\ms T$ will be important, and it can be readily obtained from the above explicit construction as follows.
Put $f(z):=-\prod_{i=1}^{g+1}(z-\lmd_i)(z-\lmd'_i)$ for simplicity.
As in Section \ref{ss:J}, let $K_i$ be the interval $\{\lmd_i<z<\lmd'_i\}$ in $\Lmd\us\simeq \mathbb S^1$.
These were exactly the locus in $\Lmd^{\sigma}$ where $f(z)>0$.
We define subsets $\mathbb D_1,\dots, \mathbb D_{g+1}$ of the cone ${\rm C}(\Lmd)$ by
\begin{align}\label{Di}
\mathbb D_i = \big\{(z,w)\in {\rm C}(\Lmd)\set z\in \ol K_i,\,\,w\in \RR, \,\, |w|
\le 
\sqrt{f(z)}\big\},\quad i = 1,\dots, g+1.
\end{align}
These are the locus in ${\rm C}(\Lmd)^{\sigma}\minus\{\text{\bf v}\}$ where the right-hand side of the equation \eqref{T} is non-negative.
Each $\mathbb D_i$ is diffeomorphic to a closed disk and in terms of the real circles $\Sigma_i^{\sigma}$ of $\Sigma$,
we have 
\begin{align}\label{D2}
\ptl \mathbb D_i= \Sigma^{\sigma}_i,\quad 1\le i\le g+1.
\end{align}
(See Figure \ref{fig:cone1}.)
\begin{midpage}
\begin{figure}
\includegraphics[height=100mm]{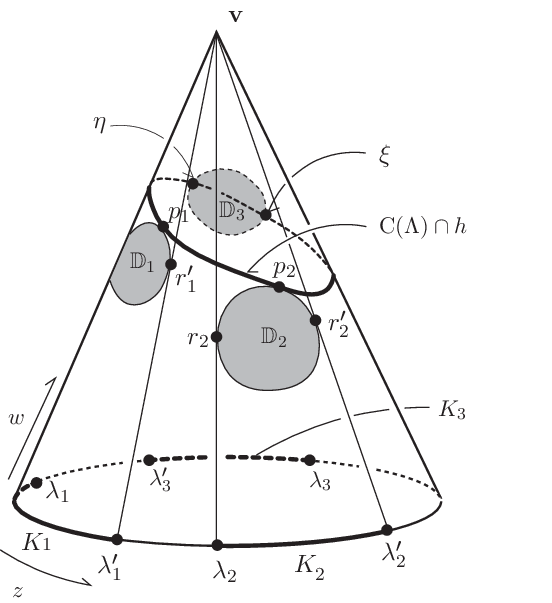}
\caption{The disks in the cone ${\rm C}(\Lmd)$ in the case $g=2$.
}
\label{fig:cone1}
\end{figure}
\end{midpage}

\begin{proposition}\label{p:rl1} 
The real locus $\ms T^{\sigma}$ of $\ms T$ consists of $(g+1)$ connected components and all of them are diffeomorphic to a smooth 2-sphere. Further, they do not pass through the singularities $\bm p_{\infty}$ nor $\ol{\bm p}_{\infty}$ of $\ms T$.
\end{proposition}

\proof 
From the equation \eqref{T} of $\ms T$, the real locus $\ms T^{\sigma}$ is a double cover of $\mathbb D_1\cup\dots\cup \mathbb D_{g+1}$ branched along the boundary circles $\ptl \mathbb D_i = \Sigma_i^{\sigma}$.
Again noting that the real locus $\ms T^{\sigma}$ is automtically smooth except for $\bf v$, $\ms T^{\sigma}$ consists of $(g+1)$ spheres smoothly embedded in $\ms T$. From \eqref{Di} these do not pass through the singularities of $\ms T$.
\proofend

\medskip
We denote the sphere that is the double cover of the disk $\mathbb D_i$
($1\le i\le g+1)$ in \eqref{Di} by $\ms T\us_i$.
Then $\ms T\us = \ms T\us_1\sqcup\dots\sqcup\ms T\us_{g+1}$,
and in the above coordinates $(z,u,v)$, they are given by
$$
\ms T\us_i = \Big\{ z\in \ol K_i,\, v=\ol u,\, |u|^2 = 
- \prod_{i=1}^{g+1}(z-\lmd_i)(z-\lmd'_i)\Big\}.
$$
From this, these spheres are invariant under the $\mathbb S^1$-action \eqref{S1act} and the ramification points
$r_i$ and $r'_i$ are exactly the $\mathbb S^1$-fixed points on $\ms T\us_i$, where we are still identifying the branch divisor $\Sigma\subset {\rm C}(\Lmd)$ with the ramification divisor in $\ms T$.

%
%

\subsection{Regular minitwistor lines}\label{ss:rmtl}
Let $S$ be a complex surface and $g>0$ an integer.
We call a divisor $C$ on $S$ a {\em regular minitwistor line with $g$ nodes}  if $C$ is reduced, irreducible, rational, have exactly $g$ ordinary nodes as its only singularities, $C\cap \Sing(S)=\emptyset$, and moreover $C$ satisfies $C^2 = 2g+2$. 
If $S$ has such a minitwistor line, then we call $S$ a {\em minitwistor space of genus $g$}.

As in the previous subsection, let $\Pi:\CP_{g+3}\lras\CP_{g+2}$ be the projection from a point that restricts to the double covering map
$\ms T\lras {\rm C}(\Lmd)$.
The latter has the hyperelliptic curve $\Sigma$ as the branch divisor.
Pick any boundary data $1\le k\le 2^{g-1}$ (see Definition \ref{def:Sei}) and let $\mf S_k$ be the corresponding Seifert surface.
Take any point $\mf s\in \mf S_k\minus\{\mf o_k,\mf o'_k\}$.
By Proposition \ref{p:Sei1}, it uniquely determines $g$ points $p_i\in \Sigma\us_i$ ($1\le i\le g$), points $\xi,\eta\in \Sigma\us_{g+1}$ and a real hyperplane $h\subset\CP_{g+2}$ such that 
\begin{align}\label{hB}
h|_{\Sigma} =\sum_{i=1}^g 2p_i + \xi+\eta.
\end{align}
The two points $\xi$ and $\eta$ were distinct if and only if $\mf s\in \mf S\uc_k$ (the interior).


By the assignment $\mf s\longmapsto h$, 
the Seifert surface $\mf S_k\minus\{\mf o_k,\mf o'_k\}$ can be regarded as a subset of the dual projective space $\CP_{g+2}^*$, and in the following, we do so (rather than thinking it as a subset of the torus $({\rm J}_{\Sigma})\us_1$).
Recall that the whole surface $\mf S_k$ does not admit a natural embedding into $\CP_{g+2}^*$ 
 because the two singular points $\mf o_k$ and $\mf o'_k$
do not uniquely determine the hyperplane as explained at the end of Section \ref{ss:Sei}.
This will be significant in Section \ref{ss:cplt} when we make the completion of an `incomplete' EW space.

\begin{proposition}\label{p:mtl1}
Let $\mf s$ be a point of the interior $\mf S_k\uc$ and
$h\subset\CP_{g+2}$ the real hyperplane determined from $\mf s$ satisfying \eqref{hB} as above.
Then the hyperplane section $\ms T\cap \Pi\inv(h)$ is a real regular minitwistor line with $g$ nodes.
Further, the intersection $\mathbb D_{g+1}\cap h$ is always an arc bounded by $\xi$ and $\eta$ and 
the real locus of the minitwistor line consists of the $g$ nodes $p_1,\dots, p_g$ and a circle in $\ms T\us_{g+1}$ which is over the  arc.
\end{proposition}

\proof
Write $C:=\ms T\cap \Pi\inv(h)$ for simplicity.
Since $h$ is real and the projection $\Pi:\ms T\lras {\rm C}(\Lmd)$ preserves the real structures, $C$ is also real.
Further, since $\Pi|_C:C\lras {\rm C}(\Lmd)\cap h$ is a double covering and $\deg \Lmd = g+1$,
\begin{align}\label{degT}
C^2 = 2\deg {\rm C}(\Lmd) = 2\deg\Lmd = 2g+2.
\end{align}
Next, since we are assuming $\mf s\in \mf S_k\uc$, ${\bf v}\not\in h$ because otherwise ${\rm C}(\Lmd)\cap h$ would consist of generating lines of the cone, which means $\eta=\tau(\xi)$, contradicting Proposition \ref{p:Sei1} (iii).
Therefore, the cut ${\rm C}(\Lmd)\cap h$ is a smooth rational curve.
Further, from \eqref{hB}, the double covering $C$ of ${\rm C}(\Lmd)\cap h$ has ordinary nodes over the tangent points $p_1,\dots,p_g\in \Sigma$
and has exactly two simple branch points $\xi$ and $\eta\in \Sigma\us_{g+1}$.
Passing to the normalization $\hat C$ of $C$, 
the composition $\hat C\lras C\lras {\rm C}(\Lmd)\cap h$ has $\xi$ and $\eta$ as the only branch points.
Since ${\rm C}(\Lmd)\cap h\simeq\Lmd\simeq\CP_1$, by Riemann-Hurwitz, this implies that $C$ is a rational curve.
In summary, $C$ is a real regular minitwistor line with $g$ nodes.

For the assertion about the real locus of $C$,
since the intersections between ${\rm C}(\Lmd)\cap h$ and $\Sigma$ at $\xi$ and $\eta$ are transversal, 
the intersections between the real locus ${\rm C}(\Lmd)\cap h\us$ and $\Sigma\us_{g+1}$ at $\xi$ and $\eta$ are also transversal.
This means that the intersection $\mathbb D_{g+1}\cap h\us$ has a (real) 1-dimensional component. So the real locus $C\us$ also has a (real) 1-dimensional component.
Therefore, the real structure induced on the normalization $\hat C\simeq\CP_1$ is identified with the complex conjugation and the real circle has to be over the intersection $\mathbb D_{g+1}\cap h\us$.
Then obviously this intersection has to be an arc whose ends are $\xi$ and $\eta$, and
$C$ cannot have a real point other than this circle and $p_1,\dots, p_g$.
\proofend

\medskip
Let $\Pi^*:\CP_{g+2}^*\hookrightarrow\CP_{g+3}^*$ be the dual map of the projection $\Pi:\CP_{g+3}\lras\CP_{g+2}$.
Then from Proposition \ref{p:mtl1}, every element of $\Pi^*(\mf S_k\uc)$ cuts out from $\ms T$ a real regular minitwistor line with $g$ nodes.

\begin{definition}\label{d:mt}{\em
For any hyperelliptic curve $\Sigma$ of genus $g$ as given in Section \ref{ss:J},
we call the compact complex surface $\ms T\subset\CP_{g+3}$ equipped with the real structure $\sigma$ {\em 
the minitwistor space associated to the real hyperelliptic curve $\Sigma$.}
Further, we call the component $\ms T\us_{g+1}\simeq\mathbb S^2$ {\em 
the real minitwistor space.}
}\proofend
\end{definition}

The reason for the last terminology is that it will be the space of (real) null surfaces in the EW spaces.

As in the previous subsection, $\ms T$ has exactly two cone singularities $\bm p_{\infty}$ and $\ol{\bm p}_{\infty}$ over the vertex $\bf v$ of the cone ${\rm C}(\Lmd)$, and is smooth away from these points.
Since $\mf S_k$ is (real) 2-dimensional, we have so far only a 2-dimensional family of real minitwistor lines for each boundary data $k$. We will soon extend this to a 3-dimensional family of minitwistor lines by using the $\mathbb S^1$-action.

\begin{remark}\label{r:rmtl}
{\em
(i)
As in Proposition \ref{p:rl1}, the real locus $\ms T\us$ consists of $(g+1)$ number of spheres $\ms T\us_1,\dots,\ms T\us_{g+1}$ and from 
Proposition \ref{p:mtl1} the last sphere $\ms T\us_{g+1}$ plays the special role. This is a consequence of the M\"obius transformation we have applied at the beginning of Section \ref{ss:ssJ} and
any $\ms T\us_i$ can be the real minitwistor space.
(ii) Under the identification of the minitwistor space $\ms T$ with the minitwistor space obtained from ALE gravitational instanton which will be given in the Appendix, the real minitwistor space is exactly the central sphere in the sense of Hitchin \cite{Hi21}.\proofend
}
\end{remark}


We will need the next property on $\ms T$.

\begin{proposition}\label{p:pn}
The surface $\ms T$ is a rational surface.
The restriction homomorphism $H^0\big(\ms O_{\CP_{g+3}}(1)\big)\lras H^0\big(\ms O_{\CP_{g+3}}(1)|_{\ms T}\big)$ is isomorphic.
Namely, $\ms T$ is non-degenerate and projectively normal in $\CP_{g+3}$.
\end{proposition}

\proof 
By \cite[Proposition 2.6]{HN11}, any compact minitwistor space is a rational surface. Hence, from Proposition \ref{p:mtl1}, the surface $\ms T$ is a rational surface. 
Since $\ms T\subset \CP_{g+3}$ is realized as a double covering of the cone ${\rm C}(\Lmd)\subset\CP_{g+2}$ and the cone is non-degenerate, $\ms T$ is also non-degenerate. So the restriction map in the proposition is injective. 
To show that it is isomorphic, it suffices to show that 
$h^0\big(\ms O_{\CP_{g+3}}(1)\big)= h^0\big(\ms O_{\CP_{g+3}}(1)|_{\ms T}\big)$.
Take a hyperplane $h\subset\CP_{g+2}$ that intersects the branch curve $\Sigma$ transversally at any point and put $C=\Pi\inv(h)\cap\ms T$. This is a smooth curve of genus $g$ as $\Sigma.\,h =\deg \Sigma = 2g+2$.
By Riemann-Roch, we obtain $h^0(\ms O_{\CP_{g+3}}(1)|_C) = g+3$. Hence, using the rationality of $\ms T$, $h^0(\ms O_{\CP_{g+3}}(1)|_{\ms T}) = g+4$.
This implies the desired equality.
\proofend

\medskip
We have obtained families of real regular minitwistor lines with $g$ nodes on $\ms T$, parameterized by the real surfaces $\Pi^*(\mf S\uc_k)$, $1\le k\le 2^{g-1}$.
For any $e^{\sqrt{-1}\theta}\in \mathbb S^1$, 
the images of these minitwistor lines under the isomorphism $e^{\sqrt{-1}\theta}:\ms T\lras \ms T$ defined in \eqref{S1act} are again minitwistor lines since $e^{\sqrt{-1}\theta}$ is a holomorphic automorphism preserving the real structure.
Hence, moving the minitwistor lines $C\in \Pi^*(\mf S_k\uc)$ by the $\mathbb S^1$-action \eqref{S1act}, 
we obtain a 3-dimensional family of real regular minitwistor lines with $g$ nodes.
One might expect from this that each space of minitwistor lines is the product manifold $\Pi^*(\mf S_k\uc) \times \mathbb S^1$.
But there might be a non-identity element $e^{\sqrt{-1}\theta}\in \mathbb S^1$ 
that preserves $\Pi^*(\mf S_k\uc)$ and if this happens, the expectation is incorrect.
The next proposition means that this happens exactly for $e^{\sqrt{-1}\theta}=-1$. It will be used to show that the EW spaces arising from $\ms T$ are diffeomorphic to $\mathbb S^2\times I$.

\begin{proposition}\label{prop:interval4}
For any $C=\ms T\cap H\in \Pi^*(\mf S_k\uc)$, we have $e^{\sqrt{-1}\theta} (C)\in \Pi^*(\mf S_k\uc)$ if and only if $e^{\sqrt{-1}\theta}=\pm 1$. 
Further, $C$ and $(-1)C$ belong to distinct connected components of the surface $\Pi^*(\mf S\uc_k)\subset\CP_{g+3}^*$. 
\end{proposition}

\proof From the explicit form \eqref{S1act2} of the $\mathbb S^1$-action on $\ms T$, 
the automorphism of $\ms T$ given by $e^{\sqrt{-1}\theta}\in \mathbb S^1$ is a lift of an automorphism of the cone ${\rm C}(\Lmd)$ if and only if $e^{\sqrt{-1}\theta} = \pm1$. 
Therefore, $e^{\sqrt{-1}\theta} (C)\in \Pi^*(\mf S_k\uc)$ can hold only when 
$e^{\sqrt{-1}\theta} = \pm1$. 
Conversely, let $\iota$ be the involution on ${\rm C}(\Lmd)$ whose lift on $\ms T$ is the one given by $e^{\sqrt{-1}\theta} = -1$. 
Again from the explicit form \eqref{S1act2}, we have
\begin{align}\label{iota}
\iota(z,w)= (z,-w).
\end{align}
From the equation \eqref{Sigma} of $\Sigma$, this means that $\iota$ preserves $\Sigma$, and on $\Sigma$, $\iota$ is equal to the hyperelliptic involution $\tau$.
In particular, $\iota$ preserves each of the circles $\Sigma\us_1,\dots,\Sigma\us_{g+1}$.
So if $h\subset\CP_{g+2}$ is the hyperplane such that $\ms T\cap\Pi\inv(h)= C$, then since the cut ${\rm C}(\Lmd)\cap h$ is tangent to $\Sigma$ at points on $\Sigma\us_1,\dots,\Sigma\us_g$ and intersects $\Sigma\us_{g+1}$ at two distinct points, so is the hyperplane $\iota(h)$.
Since the union $\cup_{k=1}^{2^{g-1}}\mf S_k\uc$ exhausts all hyperplanes in $\CP_{g+2}$ having these properties by Proposition \ref{prop:charW}, this means that
$\iota(h)\in \mf S_{l}\uc$ for some $1\le l\le 2^{g-1}$.
Since $\iota$ is continuous, this implies $\iota(\mf S_k\uc) = \mf S_l\uc$ for some $1\le l\le 2^{g-1}$.
Taking the closure, we obtain $\iota(\mf S_k) = \mf S_l$.

To show $k=l$, we make use of the boundaries of the Seifert surfaces.
As before, let $\ms S_k^1\cup\hat{\ms S}^1_k$ be the boundary of $\mf S_k$ and regard the circle $\ms S_k^1$ as a subset of $({\rm J}_{\Sigma})\us_1=W_{1,2,\dots,g}$ as it was.
Since $\iota(\Sigma)=\Sigma$ as above, $\iota$ naturally induces an involution on ${\rm J}_{\Sigma}$, and since $\iota|_{\Sigma} = \tau$ as above, the latter involution is equal to $(-\id)$.
Therefore, as shown at the beginning of Section \ref{ss:Sei}, $\iota(\ms S^1_k) = \ms S^1_k$.
Similarly, $\iota(\hat{\ms S}^1_k) = \hat{\ms S}^1_k$.
Since $\mf S_k\cap\mf S_{l}\neq\emptyset$ if and only if $k=l$ as in \eqref{disj1}, these imply that $\iota(\mf S_k) = \mf S_k$ for any $1\le k\le 2^{g-1}$. 
Thus, we have obtained, for any $1\le k\le2^{g-1}$,
\begin{align}\label{iota3}
\iota(\mf S_k) = \mf S_k,\quad
\iota(\mf S_k\uc) = \mf S_k\uc,\quad
\iota(\ms S^1_k) = \ms S^1_k,\quad
\iota(\hat{\ms S}^1_k) = \iota(\hat{\ms S}^1_k).
\end{align}

Finally, we show that the above involution $\iota$ on $\CP_{g+2}$ exchanges the two connected components of $\mf S_k\uc$ for any $1\le k\le 2^{g-1}$, again by making use of the boundary circles.
From \eqref{iota3}, either $\iota$ preserves each connected component of $\mf S_k\uc$ or $\iota$ flips the two connected components of $\mf S_k\uc$. 
To show that the latter happens,
let $\rho_1,\dots\rho_{g+1}$
be choices such that $\ms S^1_k=\ms S^1_k(\rho_1,\dots,\rho_{g+1})$
and $\rho'_i$ the unique element of $\{r_i,r'_i\}\minus\{\rho_i\}$.
Since $\rho_i$ and $\rho'_i$ are branch points of $\pi:\Sigma\lras \CP_1$, the two intersection points $\mf o_k=\sum_{i=1}^g\rho_i$ and $\mf o'_k=\sum_{i=1}^g\rho'_i$ of $\ms S^1_k$ and $\hat{\ms S}^1_k$ are fixed by the involution $\iota|_{\Sigma} = \tau$, and all other points of $\ms S^1_k\cup\hat{\ms S}^1_k$ are not fixed by the involution. Hence, on each of $\ms S^1_k$ and $\hat{\ms S}^1_k$, $\iota$ is a reflection that fixes the two points $\mf o_k$ and $\mf o'_k$.
On the other hand, from the construction of the Seifert surfaces in Section \ref{ss:Sei}, the two connected components of $\ms S^1_k\minus\{\mf o_k,\mf o'_k\}$ are parts of the boundaries of distinct connected components of $\mf S_k\uc$ (see Figure \ref{fig:Seifert2}).
Hence, $\iota$ has to flip the two connected components of $\mf S_k\uc$.
\proofend

\begin{figure}
\includegraphics[height=50mm]{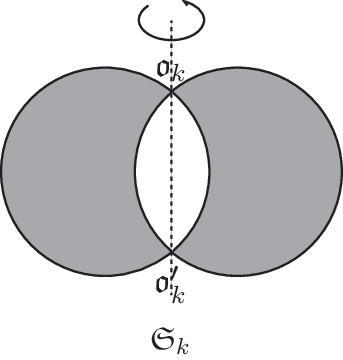}
\caption{
Obtaining $M_k\uc$ from the Seifert surface $\mf S_k$
}
\label{fig:Seifert3}
\end{figure}

\medskip
As above, for each $e^{i\theta}\in \mathbb S^1$ and $C\in \Pi^*(\mf S_k\uc)$, the image $e^{i\theta}(C)$ is a regular minitwistor line with $g$ real nodes.
Let 
$$
M_k\uc:=\big\{ e^{i\theta}(C) \set e^{i\theta}\in \mathbb S^1,\,
C\in \Pi^*(\mf S_k\uc)\big\}\subset\CP_{g+3}^*
$$ 
be the space of minitwistor lines obtained this way.
By Proposition \ref{prop:interval4}, 
$M_k\uc$ can be understood as the space obtained from one of the two connected components of $\mf S_k\uc $ by rotating around a line through $\mf o_k$ and $\mf o'_k$ (see Figure \ref{fig:Seifert3}). So we obtain a diffeomorphism
\begin{align}\label{W0}
M_k\uc
\simeq
\big(\mathbb S^2\minus\{\text{2 points}\}\big)\times I,
\end{align}
where $I$ is an open interval in $\RR$ and the two removed points correspond to the two singularties $\mf o_k$ and $\mf o'_k$ of $\mf S_k$.
Since $\ms T\us_1,\dots,\ms T\us_{g+1}$ are $\mathbb S^1$-invariant as we have already shown,
the tangential and intersection properties of $C\in \Pi^*(\mf S_k)$ with the spheres $\ms T_1\us,\dots,\ms T_{g+1}\us$ is preserved by the $\mathbb S^1$-action.
Hence, any $e^{i\theta}(C)\in M_k\uc$ is tangent to $\ms T$ at a point of each $\ms T_1\us,\dots,\ms T_{g}\us$ and intersects the real minitwistor $\ms T_{g+1}\us$ in a circle.
Hence, we have obtained a 3-dimensional family of real regular minitwistor lines with $g$ nodes on $\ms T$ for any $1\le k\le 2^{g-1}$.

\begin{proposition}\label{p:W0}
For any boundary data $1\le k\le 2^{g-1}$, the 3-manifold $M_k\uc$ admits a real analytic Lorentzian Einstein-Weyl structure.
\end{proposition}

\proof This can be shown in the same way to \cite[Proposition 3.16]{HN22}.
Indeed, since $M_k\uc$ parameterizes real regular minitwistor lines, by the result of \cite{HN11}, $M_k\uc$ admits a natural EW structure. 
Since the minitwistor lines parameterized by $M_k\uc$ have real circles as above, the EW structure is Lorentzian.
The real analyticity of the structure is a consequence of the fact that the EW space $M_k\uc$ is obtained as a real slice of holomorphic EW structure on a 3-dimensional complex manifold, and in turn this fact can be shown (without using the deformation theory of singular curves) as follows.

Take any minitwistor line $C\in M_k\uc$. Then taking a real tubular neighborhood $T$ of $C$ in $\ms T$, as in \cite{HN11}, we can take a `normalization' $\hat T$ of $T$, and then the inverse image of $C$ to $\hat T$ has a component which is a smooth rational curve with self-intersection number two.
Hence, $\hat T$ is a minitwistor space in the original sense and it is equipped with a real structure since $T$ is chosen to be real. 
So by \cite{Hi82} we obtain a holomorphic EW space as well as a real EW space as a real slice. The latter is automatically real analytic. 
Further, descending the twistor lines to $T$, we obtain equisingular displacements of $C$ in $T$ and
all such displacement of $C$ in $T$ can be obtained this way.
Thus, the EW structure on $M_k\uc$ is a real slice of a complex EW structure.
\proofend

\medskip
From the construction, the EW space $M_k\uc$ admits an $\mathbb S^1$-action induced from that on $\ms T$.

As in \eqref{W0}, the 3-manifold $M_k\uc$ has two `holes', which can be regarded as a rotation axis.
In the next subsection, we will see that these holes can be filled 
by the space of hyperplane sections that have an appropriate degeneracy and prove that the EW structure on $M_k\uc$ extends real analytically on the whole space $\mathbb S^2\times I$.

\subsection{Irreguglar minitwistor lines}\label{ss:imtl}
As in Section \ref{ss:mt}, the minitwistor space $\ms T$ has a (rational) projection $\pi\circ\Pi:\ms T\lras\Lmd$ and this gives $\ms T$ the structure of (rational) conic bundles over $\Lmd$, whose discriminant locus is the $(2g+2)$ points $\lmd_1,\lmd'_1,\dots, \lmd_{g+1}, \lmd'_{g+1}$.
The fibers over these points consist of two lines.
Using the equation \eqref{T1} of $\ms T$, 
they are concretely given by
\begin{align}\label{li}
l(r_i):=\big\{(z,u,v)\set v = 0,\,z=\lmd_i\big\},
\quad
\ol{l(r_i)}:=\big\{(z,u,v)\set u = 0,\,z=\lmd_i\big\},
\end{align}
and
\begin{align}\label{lip}
l(r'_i):=\big\{(z,u,v)\set v = 0,\,z=\lmd'_i\big\},
\quad
\ol{l(r'_i)}:=\big\{(z,u,v)\set u = 0,\,z=\lmd'_i\big\},
\end{align}
where $1\le i\le g+1$. 
We have $l(r_i)\cap \ol{l(r_i)} = \{r_i\}$ and $l(r'_i)\cap \ol{l(r'_i)} = \{r'_i\}$. 
The two lines in \eqref{li} and \eqref{lip} are real pairs.
All of $l(r_i)$ and $l(r'_i)$ pass through the same singularity of $\ms T$ and all of $\ol {l(r_i)}$ and $\ol {l(r'_i)}$ pass through another singularity of $\ms T$.
We may suppose that $\bm p_{\infty}\in l(r_i), l(r'_i)$ and
$\ol{\bm p}_{\infty}\in \ol l(r_i), \ol l(r'_i)$.
It is easy to show that the lines \eqref{li} and \eqref{lip} are all straight lines on the minitwistor space $\ms T\subset\CP_{g+3}$.

%

\begin{figure}
\includegraphics[height=60mm]{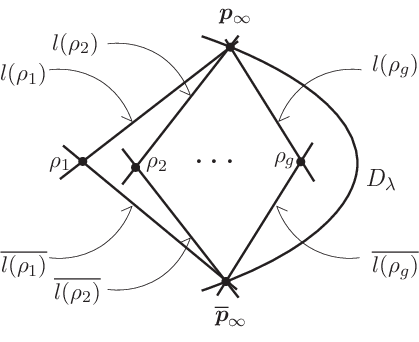}
\caption{
The rregular minitwistor line
}
\label{fig:imtl0}
\end{figure}

As before, let $1\le k\le 2^{g-1}$ be any boundary data and 
$\rho_1,\dots,\rho_{g}$ (where $\rho_i\in\{r_i,r'_i\}$ with the normalization $\rho_1 = r_1$) be the choices that corresponds to $k$.
For any real point $\lmd\in\Lmd$ which belongs to the open interval $K_{g+1}$ (bounded by $\lmd_{g+1}$ and $\lmd'_{g+1}$), let $D_{\lmd}$ be the fiber conic over the point $\lmd$, which is real and irreducible. 
The $(g+1)$ points $\lmd$ and $\rho_1,\dots,\rho_g$ span a real hyperplane in $\CP_{g+1}$ and if $H\subset \CP_{g+3}$ is the real hyperplane obtained as the pullback of this hyperplane by the projection $\pi\circ\Pi:\CP_{g+3}\lras\CP_{g+1}$, we have (see Figure \ref{fig:imtl0})
\begin{align}\label{imtl1}
H|_{\ms T} = \sum_{i=1}^g \big\{l(\rho_i) + \ol {l(\rho_i)}\big\} + D_{\lmd}.
\end{align}
In a similar way, letting $\rho'_i$ be the unique element of $\{r_i,r'_i\}\minus\{\rho_i\}$ as before and using the hyperplane spanned by the points $\rho'_1,\dots,\rho'_g$ and $\lmd$ instead, we obtain the real hyperplane $H'\subset\CP_{g+3}$ that satisfies
\begin{align}\label{imtl1'}
H'|_{\ms T} = \sum_{i=1}^g \big\{l(\rho'_i) + \ol {l(\rho'_i)}\big\} + D_{\lmd}.
\end{align}
We call these reducible curves on $\ms T$
{\em irregular minitwistor lines}.
As above, they are determined from the boundary data $k$ and a point $\lmd\in K_{g+1}$ and therefore for any fixed $k$ they are parameterized by two open intervals, both of which are identified with $K_{g+1}$.
In the following, we denote
\begin{align}\label{}
\bm I_k\qandq \bm I'_k
\end{align}
for the families of irregular minitwistor lines that are of the form \eqref{imtl1} and \eqref{imtl1'} respectively.
%
%
Since they are hyperplane sections, we have a natural inclusion
$\bm I_k\sqcup\bm I'_k\subset\CP_{g+3}^*$.
These pairs of intervals play an important role in proving the closedness of all spacelike geodesics on the EW spaces.

\subsection{The completion of $M_k\uc$}\label{ss:cplt}
In this subsection, we show that the two intervals $\bm I_k$ and 
${\bm I}'_k$
are naturally attached to the 3-manifold $M_k\uc\simeq (\mathbb S^2 \minus\{\text{2 points}\})\times I$ in \eqref{W0} to fill the two `holes', namely $\{\text{2 points}\}\times I$.

For this purpose, we define a subvariety $V_2$ of the product variety ${\rm J}_{\Sigma}\times W_2$ by
$$
V_2:=\big\{(\mf u,\mf w)\in {\rm J}_{\Sigma}\times W_2\set 2\mf u = \mf w \big\}.
$$
Obviously, the projection from $V_2$ to the second factor $W_2$ realizes $V_2$ as a $4^g$-fold unramified covering over $W_2$. In particular, by Proposition \ref{prop:W2}, $V_2$ has exactly $4^g$ singular points and they are cyclic quotient singularities of order $(g-1)$.
They are mapped to ${\rm J}_{\Sigma}[2]$ bijectively by the projection to the first factor ${\rm J}_{\Sigma}$. It might be convenient to regard $V_2$ as the graph of the multi-valued mapping $\mf w\longmapsto \frac 12 \mf w$ from $W_2$ to ${\rm J}_{\Sigma}$.

Next, let 
$$
\aaa:\Sigma^{(g)}\times \Sigma\ut \lras {\rm J}_{\Sigma}\times W_2
$$
be the product map $\mf a\ug\times\mf a\ut$.
Since both $\mf a\ug$ and $\mf a\ut$ are birational,
so is $\aaa$.
Let $\tilde V_2 \subset \Sigma\ug\times \Sigma\ut$ be the strict transform of $V_2$ into $\Sigma^{(g)}\times \Sigma\ut$.
Let 
$$
\bbb:\Sigma\ug \times \Sigma\ut \lras \Sigma^{(2g+2)}
$$
be the mapping defined by $(D_1,D_2) \longmapsto 2D_1 + D_2$.
This is a morphism of varieties.
Using $\aaa$, we have the following commutative diagram:
\begin{equation}\label{cd1}
   \vcenter{
   \xymatrix@=18pt{
\tilde V_2 
\ar @{^{(}->}[r]
\ar[d] 
& \Sigma\ug\times \Sigma\ut 
\ar[d]_{\aaa}
\ar[r]^-{\pr_2} 
& \Sigma\ut
\ar[d]^{\mf a^{(2)}} \\
V_2
\ar @{^{(}->}[r] 
& {\rm J}_{\Sigma}\times W_2
\ar[r] 
& W_2
}
}
\end{equation}
Since $\mf a\ut:\Sigma\ut\lras W_2$ is the minimal resolution of the quotient singularity $\mf o$ of $W_2$ as in Proposition \ref{prop:W2}, this diagram \eqref{cd1} implies that the mapping $\aaa|_{\tilde V_2}:\tilde V_2\lras W_2$ is exactly the minimal resolutions of all the $4^g$ singularities of $V_2$.
The exceptional curve of $\mf a\ut:\Sigma\ut\lras W_2$ is naturally identified with the rational normal curve $\Lmd$ (over which the hyperelliptic curve $\Sigma$ lies as a double covering).

For any $1\le k\le 2^{g-1}$, the Seifert surface $\mf S_k$ was a subset of the connected component $({\rm J}_{\Sigma})_1\us=W_{1,2,\dots,g}$ of $({\rm J}_{\Sigma})\us$ (see Definition \ref{def:Sei}), so any point $\mf s\in \mf S_k$ can be uniquely written as $\sum_{i=1}^g \mf a(p_i)$ with $p_i\in \Sigma\us_i$ that satisfies $\mf t(\sum_{i=1}^g \mf a(p_i)) = \mf a(\xi+\eta)$ for some $\xi,\eta\in \Sigma\us_{g+1}$.
Moreover, $\xi$ and $\eta$ were unique if $\mf s\not\in\{\mf o_k,\mf o'_k\}$.
Hence, the point $(\sum_{i=1}^g \mf a(p_i),\mf a(\xi+\eta))$ belongs to $V_2$, and the assignment $\mf S_k\minus\{\mf o_k,\mf o'_k\}\ni\sum_{i=1}^g \mf a(p_i)\longmapsto (\sum_{i=1}^g \mf a(p_i),\mf a(\xi+\eta))$ gives an embedding of $\mf S_k\minus\{\mf o_k,\mf o'_k\}$ into $V_2$.
If we regard $\mf S_k\minus\{\mf o_k,\mf o'_k\}$ as a subset of $V_2$ this way, then it lies over the subset $W_{g+1,\,g+1}$ of $W_2$ under the $4^g$-fold covering $V_2\lras W_2$.

\begin{definition}\label{d:ann}{\em
By the {\em strict transform of $\mf S_k$ into $\tilde V_2$}, we mean 
the closure of $\aaa\inv(\mf S_k\minus \{\mf o_k,\mf o'_k\})$ in $\tilde V_2$, and we write $\tilde{\mf S}_k$ for it.
Further, we mean the interior of $\tilde{\mf S}_k$ by $\tilde{\mf S}_k\uc$.
}
\end{definition}

\begin{figure}
\includegraphics[height=40mm]{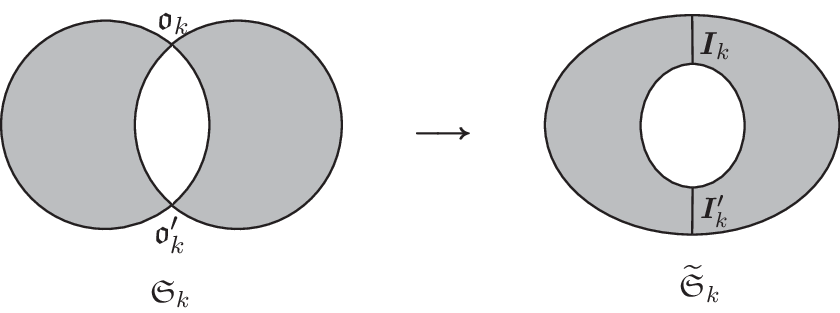}
\caption{ From $\mf S_k$ to the resolution $\tilde{\mf S}_k$ 
}
\label{fig:Seifert4}
\end{figure}

From the above identification of the exceptional curves of the minimal resolution $\tilde V_2\lras V_2$ with $\Lmd$, 
the inverse images for the mapping $\tilde{\mf S}_k\lras \mf S_k$ of
the two singular points $\mf o_k$ and $\mf o'_k$ are both identified with the interval $K_{g+1}\subset\Lmd\us$ because if 
$\sum_{i=1}^g \mf a(p_i)$ is a generic point of $\mf S_k$, then
the two points $\xi,\eta$ on $\Sigma$ that satisfy $\mf t_1(\sum_{i=1}^g \mf a(p_i)) = \mf a(\xi+\eta)$ belong to $\Sigma\us_{g+1}$ and $\pi(\Sigma\us_{g+1}) = \ol K_{g+1}$.
Hence, the strict transform $\tilde{\mf S}_k$ of any Seifert surface $\mf S_k$ is obtained from $\mf S_k$ by replacing the two singularities $\mf o_k$ and $\mf o'_k$ with the intervals $\bm I_k$ and ${\bm I}'_k$ given in the previous subsection.
Therefore, $\tilde{\mf S}_k$ has a structure of a closed interval bundle over $\mathbb S^1$. So it is diffeomorphic to either an annulus or a M\"obius strip. 
But since the boundary of $\tilde{\mf S}_k$ consists of (the strict transform of) the two connected components $\ms S^1_k$ and $\hat{\ms S}^1_k$,  $\tilde{\mf S}_k$ is diffeomorphic to a closed annulus
$\mathbb S^1\times \ol I$, where $\ol I$ is an closed interval in $\RR$. 
The interior $\tilde{\mf S}_k\uc$ is diffeomorphic to an open annulus
$\mathbb S^1\times I$.
(See Figure \ref{fig:Seifert4}.)

Similarly to \eqref{cd1}, there is another obvious commutative diagram
 \begin{equation}\label{cd2}
   \vcenter{
\xymatrix@=18pt{
\tilde V_2 
\ar @{^{(}->}[r]
\ar[d] 
& \Sigma\ug\times \Sigma\ut 
\ar[d]_{\aaa}
\ar[r]^-{\bbb} 
& \Sigma^{(2g+2)} 
\ar[d]^{\mf a^{(2g+2)}} 
\\
V_2  
\ar @{^{(}->}[r]
& {\rm J}_{\Sigma}\times W_2 
\ar[r]^-{\ol\bbb} 
& {\rm J}_{\Sigma} 
}
}
\end{equation}
where $\ol\bbb:{\rm J}_{\Sigma}\times W_2\lras {\rm J}_{\Sigma}$ is defined as
the assignment $(\mf a(D_1),\mf a(D_2))\longmapsto\mf a(2D_1+ D_2)$.
Although the mapping $\bbb$ is clearly not injective, the restriction of it to the subvariety $\tilde V_2$ is injective as $\Sigma\us_i\cap \Sigma\us_{g+1} = \emptyset$ if $i<g+1$.
Since the hyperelliptic curve $\Sigma$ is realized in $\CP_{g+2}$ as a non-degenerete curve of degree $2g+2$, there is a natural injection from the dual projective space $\CP_{g+2}^*$ to  $\Sigma^{(2g+2)}$ which assigns the divisor $h|_{\Sigma}$ to each hyperplane $h\in \CP_{g+2}^*$.
The image of this map is exactly the inverse image of the origin $\mf o\in {\rm J}_{\Sigma}$ under the Abel-Jacobi map $\mf a^{(2g+2)}:\Sigma^{(2g+2)}\lras {\rm J}_{\Sigma}$. 
From these, the subvariety $\tilde V_2$ is naturally embedded in $\CP_{g+2}^*$.
Therefore, (unlike $\mf S_k$ itself) the strict transform $\tilde{\mf S}_k$ of each Seifert surface $\mf S_k$ can be realized as a space of hyperplanes in $\CP_{g+2}$ in which the hyperelliptic curve $\Sigma$ is embedded.

As before, let $\Pi:\CP_{g+3}\lras\CP_{g+2}$ be the projection from a point that realizes $\ms T$ as a double covering over the cone ${\rm C}(\Lmd)$, and $\Pi^*:\CP_{g+2}^*\hookrightarrow\CP_{g+3}^*$ the dual embedding. 
Then since $\Pi^*$ is an embedding, the image $\Pi^*(\tilde{\mf S}_k)$ 
is diffeomorphic to $\tilde{\mf S}_k$ 
and as above the latter is diffeomorphic to an annulus $\mathbb S^1\times \ol I$.
Similarly to the definition of $M_k\uc$, using the $\mathbb S^1$-action \eqref{S1act}, we define 
\begin{align}\label{Wk}
\ol M_k:=\big\{ e^{i\theta}(C) \set e^{i\theta}\in \mathbb S^1,\,
C\in \Pi^*( \tilde{\mf S}_k)\big\}\\
M_k:=\big\{ e^{i\theta}(C) \set e^{i\theta}\in \mathbb S^1,\,
C\in \Pi^*( \tilde{\mf S}_k\uc)\big\},
\end{align}
where $\tilde{\mf S}_k\uc$ means the interior of $\tilde{\mf S}_k$.
Since the induced $\mathbb S^1$-action on $\CP_{g+3}^*$ fixes points of the form $\Pi^*(h)$, $h\in \bm I_k\cup{\bm I}'_k$,
moving the hyperplanes belonging to $\Pi^*(\tilde{\mf S}_k)$ by the $\mathbb S^1$-action amounts to 
rotating the annulus $\Pi^*(\tilde{\mf S}_k)$ around the intervals $\bm I_k\cup{\bm I}'_k$. (See Figure \ref{fig:Seifert4}, right.)
Therefore, the subset $\ol M_k$ (resp.\,$M_k$) of $\CP_{g+3}^*$ is diffeomorphic to $\mathbb S^2\times \ol I$ (resp.\,$\mathbb S^2\times I$). 
The space $M_k$ naturally contains the EW space $M\uc_k$ given in \eqref{W0} and its complement in $M_k$ is identified with $\bm I_k\cup {\bm I}'_k$ through the pullback $\Pi^*$.
Thus, we have obtained the natural `completion' of $M\uc_k$ for any $1\le k\le 2^{g-1}$.
In the next section, we will show that the EW structure on $M\uc_k$ given in Proposition \ref{p:W0} extends real analytically to the 3-manifold $M_k\simeq \mathbb S^2\times I$.

%

\begin{definition}\label{d:infty}{\em
We call the two connected components of the boundary $\ol M_k\minus M_k$ the {\em future infinity} and the {\em past infinity} of the 3-manifold $M_k$. These are 
diffeomorphic to the 2-sphere.
We denote $\ms H_k$ and $\hat{\ms H_k}$ for the future and past infinity of $M_k$.}
\end{definition}

These boundaries will play some role in the next section in proving the basic properties of the EW structure on $M_k$.

Finally in this section, 
we give a characterization of elements of the 3-manifolds $M_k$ and its closure $\ol M_k$, viewed as hyperplanes in $\CP_{g+3}$ in which $\ms T$ is embedded.
It is an immediate consequence of Proposition \ref{p:Sei1}.
Recall that each boundary data $k$ ($1\le k\le 2^{g-1}$) determines 
an element $\rho_i\in\{r_i,r'_i\}$ for any $1\le i\le g$
under the normalization $\rho_1=r_1$
and we have written $\rho'_i$ for the unique element of $\{r_i,r'_i\}\minus\{\rho_i\}$. Also, for any point $\lmd\in\Lmd$, we denote $D_{\lmd}$ for $(\pi\circ\Pi)\inv(\lmd)$ which is the fiber conic.

\begin{proposition}\label{prop:charW}
If $H\subset\CP_{g+3}$ is a hyperplane belonging to $M_k$ for some $1\le k\le 2^{g-1}$, then $H$ satisfies the following properties:
\begin{itemize}
\setlength{\parskip}{0cm}
\item[\em (i)]for any $1\le i\le g$, $H$ is tangent to $\ms T$ at a single point of the sphere $\ms T_i\us$, 
\item[\em (ii)] the intersection $\ms T\us_{g+1}\cap H$ is a (smooth) circle.
\end{itemize}
If in addition $H\in \bm I_k$ or $H\in\bm {\bm I}'_k$, then $\ms T\cap H$ is respectively of the form
\begin{align}\label{irrmtl}
\sum_{1\le i\le g} \Big\{l(\rho_i) + \ol{l(\rho_i)}\Big\} + D_{\lmd} \,{\text{ or }}\,
\sum_{1\le i\le g} \left\{l(\rho'_i) + \ol{l(\rho'_i)}\right\} + D_{\lmd}
\end{align}
for some $\lmd\in K_{g+1}$.
If $H$ belongs to the boundary $\ptl M_k=\ms H_k\cup\hat{\ms H}_k$, then $H$ satisfies (i) and the property
\begin{itemize}
\item[\em (ii)'] $\ms T\us_{g+1}\cap H$ is a single point and $H$ is tangent to $\ms T$ at the point.
\end{itemize}
Conversely, these properties characterize real hyperplanes in $\CP_{g+3}$ belonging to $M_k$, $\bm I_k\cup \bm I'_k$ and $\ptl M_k$ respectively.
\end{proposition}

Note that the properties (i), (ii) and (ii)' are independent of $1\le k\le 2^{g-1}$ and the dependence on $k$ appears only on the property \eqref{irrmtl}.

\begin{remark}\label{rmk:charW}
{\em
As in the proposition, any hyperplane $H\in \ptl M_k$ intersects the real minitwistor space $\ms T\us_{g+1}$ at a single point.
Therefore, the proposition remains true even if we replace the property (ii) with the one that ``the intersection $\ms T\us_{g+1}\cap H$ consists of {\em at least two points}''.
This stronger version will be used at the final step in the proof of the closedness of all spacelike geodesics in $M_k$.
}
\end{remark}

The following property will also be used in the same proof:
\begin{proposition}\label{prop:disj1}
If $1\le k\neq l\le 2^{g-1}$, then $M_k\cap M_l=\emptyset$ in $\CP_{g+3}^*$.
\end{proposition}

\proof
Obviously, any distinct intervals of the form 
$\bm I_k$ and $\bm I'_k$
are disjoint.
So to prove $M_k\cap M_l=\emptyset$ if $k\neq l$,  
it suffices to show $M_k\uc\cap M_l\uc=\emptyset$ if $k\neq l$.
Since the Seifert surfaces are connected components of $\mf t_1\inv(W_{g+1,\,g+1})$ as in \eqref{disj1}, we have $\mf S_k\cap\mf S_l=\emptyset$ if $k\neq l$.
Since $M_k\uc$ and $M_l\uc$ are obtained from $\mf S_k$ and $\mf S_l$ by using the $\mathbb S^1$-action, this means $M_k\uc\cap M_l\uc=\emptyset$ if $k\neq l$.
\proofend

\section{Properties of the Einstein-Weyl structures}
\label{s:EW}
\subsection{Extension of the EW structures}
\label{ss:modify}
As before, let $\ms T$ be the minitwistor space associated to the real hyperelliptic curve $\Sigma$ of genus $g$ (see Definition \ref{d:mt}).
In this subsection, we first select an open subset $U$ of $\ms T$ and next apply a modification to it, obtaining an open smooth complex surface $\check U$, which will be a minitwistor space in the original sense.
The construction is a natural generalization of those given in \cite[Section 4.1]{HN22} to the case of arbitrary genus.
In the construction, the role of lines $l_1$ and $l_2$ in \cite{HN22} is played by the $2g$ lines $l(r_1),l(r'_1),\dots, l(r_g),l(r'_g)$, and that of $l_3$ and $l_4$ is played by $l(r_{g+1})$ and $l(r'_{g+1})$. 

As before, let $k$ be any boundary data and for any $1\le i\le g$, let $\rho_i$ be one of the points $r_i,r'_i$ determined from a boundary data $k$.
We have chosen $\rho_1=r_1$, but this assumption will not be needed. 
These determine $2g$ lines $\{l(\rho_i),\ol{l(\rho_i)}\set 1\le i\le g\}$ on $\ms T$, see \eqref{li} and \eqref{lip} for the definition of these lines. 
For simplicity, in the sequel, we write $l_i = l(\rho_i)$ and $\ol l_i=\ol{l(\rho_i)}$.
Let 
$$
\mu:\tilde{\ms T}\lras \ms T
$$ 
be the minimal resolution of the two singularities $\bm p_{\infty}$ and $\ol{\bm p}_{\infty}$ of $\ms T$, and $E$ and $\ol E$ the exceptional curves over the singularities respectively. 
The latter are smooth rational curves
satisfying $E^2 = \ol E^2 = -(g+1)$.
The resolution eliminates the indeterminacy of the rational conic bundle structure $\pi\circ\Pi:\ms T\lras\Lmd$, so we obtain a genuine conic bundle for which we denote $\tilde f:\tilde{\ms T}\lras \Lmd$.
The curves $E$ and $\ol E$ are sections of $\tilde f$.
We use the same symbols to mean the strict transforms of the above lines into $\tilde{\ms T}$.
Take any point $\lmd\in K_{g+1}$ and 
let $D_{\lmd}$ be the fiber conic $\tilde f\inv(\lmd)$.

\begin{figure}
\includegraphics{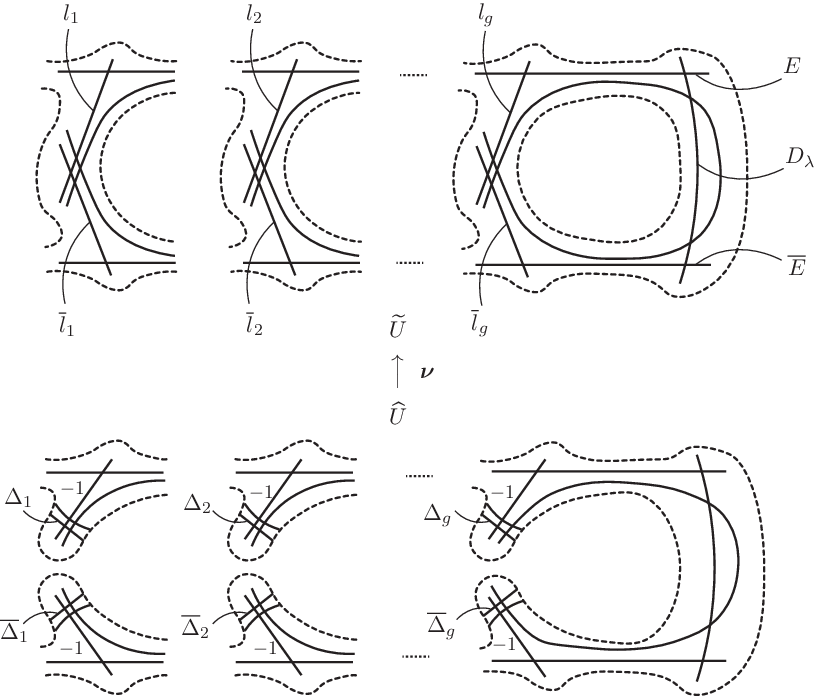}
\caption{The `normalization' of the open subset $\tilde U\subset\ms T$.
}
\label{fig:modif}
\end{figure}

We choose a tubular neighborhood $\tilde U$ in $\tilde{\ms T}$ of the reducible curve
$$
\big(
 l_1\cup \ol l_1\cup\dots\cup
l_g\cup \ol l_g
\big)
\cup
E\cup \ol E \cup D_{\lmd},
$$
chosen in such a way that it contains a subset of the form $\tilde f\inv(Y)$ for an open real subset $Y\subset\Lmd$ that contains the interval $K_{g+1}$.
Making use of a $\sigma$-invariant Riemannian metric on $\tilde {\ms T}$ that is invariant under the $\mathbb S^1$-action, the tubular neighborhood $\tilde U$ can be taken to be real and $\mathbb S^1$-invariant.
Let $\boldsymbol\nu:\hat U\lras \tilde U$ be the ``normalization'' around the intersection points $\rho_i = l_i\cap \ol l_i$,  $1\le i\le g$. See Figure \ref{fig:modif}.
Each $\boldsymbol\nu\inv(l_i)$ and $\boldsymbol\nu\inv(\ol l_i)$, $1\le i\le g$, consists of two irreducible components, one of which is mapped isomorphically onto $l_i$ and $\ol l_i$ respectively.
We use the same letters to denote these components in $\hat U$.
We denote $\DDD_i$ and $\ol\DDD_i$ for
the other components, which are pieces of smooth curves that are mapped onto a portion of the lines $l_i$ and $\ol l_i$ respectively.
These will not be needed in the construction of the modification, but are crucial for proving Proposition \ref{p:le1}.

The curves $l_i$ and $\ol l_i$ in $\hat U$ are disjoint $(-1)$-curves. Let $\boldsymbol{\mu}_1:\hat U\lras \hat U'$ be the blowdown of these $2g$ curves. Then the images of $E$ and $\ol E$ in $\hat U'$ will be $(-1)$-curves. Let $\boldsymbol{\mu}_2:\hat U'\lras \check U$ be the blowdown of these two curves and $e$ and $\ol e\in\check U$ the image points of $E$ and $\ol E$ respectively. The real structure and the $\mathbb S^1$-action are preserved through these operations, so $\check U$ admits them. 

We put $U:=\mu(\tilde U)\subset\ms T$. This is a real $\mathbb S^1$-invariant open subset of $\ms T$.
We regard the surface $\check U$ as a modification of $U$. We then have:

\begin{proposition}\label{p:tmtl}
Every irregular minitwistor line in $U$ is transformed to a smooth rational curve in $\check U$ through $e$ and $\ol e$ and the transformation has self-intersection number two and is $\mathbb S^1$-invariant.
Regular minitwistor lines in $U$ are also transformed to real smooth rational curves with self-intersection number two, but they pass through none of $e$ nor $\ol e$ and are not $\mathbb S^1$-invariant. 
\end{proposition}

Therefore, the modification $\check U$ is a minitwistor space in the original sense and it admits an $\mathbb S^1$-action.
The proposition can be readily shown from the construction
and we omit the proof.
Using the proposition, we can prove:

\begin{theorem}\label{t:extend}
For any $1\le k\le 2^{g-1}$, the EW structure on the punctured 3-manifold $M_k\uc=M_k\minus(\bm I_k\cup\bm I'_k)$ obtained in Proposition \ref{p:W0} extends real analytically across the two intervals $\bm I_k$ and $\bm I'_k$ to give an EW structure on the whole of $M_k\simeq \mathbb S^2\times I$, preserving the $\mathbb S^1$-action.
\end{theorem}

\noindent{\em Proof.}
We still write $l_i = l(\rho_i)$ and $\ol l_i=\ol{l(\rho_i)}$.
By construction, the minitwistor space $\check U$ contains the transformations of all irregular minitwistor lines of the form 
\begin{align}\label{imtl3}
\sum_{i=1}^g\left(l_i + \ol l_i\right) + D_{\lmd},
\quad\lmd\in K_{g+1},
\end{align}
as well as the transformations of all nearby regular minitwistor lines in $U$.
From Proposition \ref{p:tmtl}, this means that the EW structure on $M_k\uc$ extends real analytically across the interval $\bm I_k$.
The same conclusion can be deduced for another interval 
$\bm I'_k$
by just replacing the role of $\sum_{i=1}^g\big(l_i + \ol{l}_i\big)$ with $\sum_{i=1}^g\big(l'_i + \ol l'_i\big)$ in the modification, where
$l'_i=l(\rho'_i)$ and $\ol l'_i=\ol{l(\rho'_i)}$.
The $\mathbb S^1$-symmetry of the EW structure on $M_k\uc$ extends to $M_k$ since the $\mathbb S^1$-action on the open subset $U$ is preserved under the modification procedure.
\proofend

\medskip
On the transformed minitwistor space $\check U$,
the minitwistor lines on it satisfy the following property.
We denote $\check D_{\lmd}$ for the transformation of irregular minitwistor line \eqref{imtl3} into $\check U$. As in Proposition \ref{p:tmtl}, these are $\mathbb S^1$-invariant smooth minitwistor lines.

\begin{proposition} \label{p:le1} 
On the open surface $\check U$, all the transformations $\check D_{\lmd}$ with $\lmd\in K_{g+1}$ are mutually linearly equivalent and satisfies $\dim |\check D_{\lmd}| = 1$.
If a smooth minitwistor line $\check C\subset \check U$ is obtained as the transformation of a regular minitwistor line $C\subset U$,
then $\dim|\check C| = 0$.
\end{proposition}

Again, this can be proved in a similar way to the case of $g=1$ in \cite[Proposition 4.4]{HN22}, using the pieces of curves $\DDD_i$ and $\ol\DDD_i$ that arise in the operation from $\tilde U$ to $\hat U$, as well as the exceptional curves $E$ and $\ol E$. Therefore, we omit the proof.
 
Since any smooth rational curve $C$ satisfying $C^2>0$ on a smooth compact complex surface satisfies $\dim |C|=C^2+1$ (see the proof of \cite[Theorem 4.6]{HN22}),
Proposition \ref{p:le1} means the following:
 
\begin{proposition}\label{p:nc}
The open complex surface $\check U$ does not admit any compactification as a complex surface.
\end{proposition}

\subsection{Spacelike Zoll property of the EW structures.}\label{ss:Zoll}

Let $M_k$ ($1\le k\le 2^{g-1}$) be the Lorentzian EW space obtained in Theorem \ref{t:extend}.
They were the spaces of regular and irregular minitwistor lines on the surface $\ms T$ given in Sections \ref{ss:rmtl} and \ref{ss:imtl}.
For any point $C\in M_k$, regarding $C$ as a minitwistor line, 
let $\hat C$ be the normalization of all the real $g$ nodes of $C$
and $C^{\RR}$ the real circle of $C$.
The normalization does not affect the real circle since the latter necessarily lies on $\ms T\us_{g+1}$.
We denote $\hat C^{\RR}$ the real circle of $\hat C$.
From the general construction \cite{Hi82} of the affine connection of the EW structure, for any point $C\in M_k$, any geodesic on $M_k$ through $C$ is given by
\begin{align}\label{slg}
\mf C_{x,y}:= M_k\cap \ol{xy}^*,
\end{align}
where $x,y\in \hat C$ and $\ol{xy}^*$ means the space hyperplanes in $\CP_{g+3}^*$ containing the line through $x,y$,
and $x,y$ satisfy one of the following conditions.
(i) $x,y\in\hat C^{\RR}$ and $x\neq y$,
(ii) $x,y\in\hat C^{\RR}$ and $x= y$
(in which case $\ol{xy}$ means the tangent line $T_{x}C$), or
(iii) $x,y\not\in \hat C^{\RR}$ and $x= \sigma(y)$
(in which case, $\ol{xy}$ should be understood as $\ol{\nu(x)\nu(y)}$).
In the cases (i) and (ii) we may think $x,y\in C^{\RR}\subset\ms T\us_{g+1}$ and the geodesic \eqref{slg} is spacelike precisely in the case (i). 
Thus every spacelike geodesic on $M_k$ is of the form $\mf C_{x,y}$ for distinct points $x,y\in\ms T\us_{g+1}$.



\begin{proposition}\label{p:slg0}
For any distinct points $x,y\in \ms T\us_{g+1}$, the intersection 
$M_k\cap \ol{xy}^{*\sigma}$ in $\CP_{g+3}^{*\sigma}$ is transveral. 
\end{proposition}

\proof 
This can be proved in a similar way to the case $g=1$ given in \cite[Proposition 4.7]{HN22} using the deformation theory of singular curves, so we only provide a sketch of the proof.

Take any point $C\in M_k\cap \ol{xy}^{*\sigma}$ and the real hyperplane $H\subset\CP_{g+3}$ through $x,y$ such that $C=\ms T\cap H$.
$C$ can be irregular.
Let $p_1,\dots, p_g$ be the real nodes of $C$.
First, suppose that $C$ is regular. 
Let $\nu:\hat C\lras C$ be the normalization of $C$ and $p'_i$ and $\ol p'_i$ the two points of $\hat C$ over the node $p_i$.
Using $\deg N_C = C^2 = \deg \ms T = 2g+2$ where $N_C=[C]|_C$, 
for the equi-singular normal sheaf $N'^{\,p_1,\dots,p_g}_C$ of $C$ in $\ms T$,
we obtain 
\begin{align}\label{Npl11}
\nu^* N'^{\,p_1,\dots,p_g}_C &\simeq 
\nu^*N_C\otimes\ms O_{\hat C}(-p'_1-\ol p'_1-\dots - p'_g - \ol p'_g)\\
&\simeq
\ms O_{\hat C}(2).
\end{align}
From this, using Leray spectral sequence, we obtain
\begin{align}\label{Npl2}
H^0\big(N'^{\,p_1,\dots,p_g}_C\big) &\simeq 
H^0\big(\nu^* N'^{\,p_1,\dots,p_g}_C\big) \simeq 
H^0\big(\ms O_{\hat C}(2)\big)\simeq\CC^3,\\
H^1\big(N'^{\,p_1,\dots,p_g}_C\big) &\simeq
H^1\big(\nu^* N'^{\,p_1,\dots,p_g}_C\big)\simeq
H^1\big(\ms O_{\hat C}(2)\big)=0.
\end{align}
This means that the versal family of equi-singular displacements of $C$ in $\ms T$ is parameterized by a 3-dimensional complex manifold $U$ and from \eqref{Npl11} and \eqref{Npl2} there is a natural isomorphism 
\begin{align}\label{tan1}
T_CU\simeq 
\big\{s\in H^0\big(\nu^*N_C\big)\set 
s(p'_i) =s(\ol p'_i)=0,\,i=1,\dots,g\big\}.
\end{align}
$U$ can be considered as a complexification of $M_k$ around the point $C$.
On the other hand, since $C\subset \ms T$ is non-degenerate and projectively normal from Proposition \ref{p:pn}, 
there is a natural isomorphism $T_H(\CP_{g+3}^*)\simeq H^0(N_C)$.
Further, by the pullback, $H^0(N_C)$ is naturally a subspace of $H^0(\nu^*N_C)$ defined by 
$\{s\in H^0(\nu^*N_C)\set s(p'_i) =s(\ol p'_i),\,i=1,\dots,g\}\simeq\CC^{g+3}$.
Hence, we obtain a natural isomorphism
\begin{align}\label{tan2}
T_C(\ol{xy}^*)\simeq
\big\{s\in H^0\big(\nu^*N_C\big)\set 
s(p'_i) =s(\ol p'_i),\,i=1,\dots,g,\,s(x)=s(y)=0\big\}.
\end{align}
Since $\{x,y\}\cap \{p_1,\ol p'_1,\dots,p_g,\ol p'_g\}=\emptyset$
as $x,y\in \ms T\us_{g+1}$ and $p_i\in\ms T\us_i$ for $1\le i\le g$, 
\eqref{tan1} and \eqref{tan2} mean that 
$\dim (T_CU\cap T_C(\ol{xy}^*)) = 1$.
Hence, $\ol{xy}^*$ and $U$ intersect transversally at $C$.
Therefore, so is the real slices $\ol{xy}^{*\sigma}$ and $M_k$.

Next, suppose that $C$ is irregular.
As in \eqref{imtl1} and \eqref{imtl1'}, it is of the form either
\begin{align}\label{imtl20}
\sum_{i=1}^g \big(l_i + \ol l_i\big) + D_{\lmd}
\,\,{\text{ or }}\,\,
\sum_{i=1}^g \big(l'_i + \ol l'_i\big) + D_{\lmd},
\end{align}
where $l_i= l(\rho_i), l'_i= l(\rho'_i)$ and $D_{\lmd}$ is a fiber conic over a point $\lmd\in K_{g+1}$.
We assume that $C$ is of the former form because the latter form can be argued in the same way.
As before, let $\mu:\tilde{\ms T}\lras\ms T$ be the minimal resolution and $E$ and $\ol E$ the exceptional curves over $\bm p_{\infty}$ and $\ol{\bm p}_{\infty}$ respectively. 
Using the same letters to mean the strict transforms into $\tilde{\ms T}$ for curves in $\ms T$,
we have
\begin{align}\label{imtl4}
\mu^*C=\sum_{i=1}^g \big(l_i + \ol l_i\big) + D_{\lmd} + E + \ol E.
\end{align}
This still has $\rho_1,\dots,\rho_g$ as its all real nodes.
Put $N_{\mu^*C}:=[\mu^*C]|_{\mu^*C}$. Then $\deg N_{\mu^*C}|_{D_{\lmd}} = 2$ as  $D_{\lmd}$ is conic and $\deg N_{\mu^*C}|_{l_i} =\deg N_{\mu^*C}|_{\ol l_i} = 1$ as $l_i$ and $\ol l_i$ are lines.

Let $\nu:\hat C\lras \mu^*C$ be the normalization of $\mu^* C$ at these $g$ nodes.
This makes $l_i$ and $\ol l_i$ disjoint.
Let $N'^{\,\rho_1,\dots,\rho_g}_{\mu^* C}$ be the equi-singular normal sheaf of $\mu^*C$ in $\tilde{\ms T}$ with the equi-singularity being imposed at $\rho_1,\dots,\rho_g$. 
We obtain 
\begin{align}\label{Npl11'}
\nu^* N'^{\,\rho_1,\dots,\rho_g}_{\mu^* C}&\simeq 
(\nu^*N_{\mu^*C})\otimes\ms O_{\hat C}(-p'_1-\ol p'_1-\dots - p'_g - \ol p'_g)
\end{align}
and
\begin{gather}\label{Npl4}
\big(\nu^* N'^{\,\rho_1,\dots,\rho_g}_{\mu^* C}\big)|_{\hat C- D_{\lmd}} \simeq 
\ms O_{\hat C- D_{\lmd}},
\quad
\big(\nu^* N'^{\,\rho_1,\dots,\rho_g}_{\mu^* C}\big)|_{D_{\lmd}}\simeq
\ms O_{D_{\lmd}}(2),\\
H^0\big(N'^{\,\rho_1,\dots,\rho_g}_{\mu^* C}\big)
\simeq
H^0\big(\ms O_{D_{\lmd}}(2)\big)\simeq\CC^3,\\
\label{Npl51}
H^1\big(N'^{\,\rho_1,\dots,\rho_g}_{\mu^* C}\big)\simeq
H^1\big(\ms O_{D_{\lmd}}(2)\big)=0.
\end{gather}
This means that the versal family of displacements of $\mu^* C$ in $\tilde {\ms T}$ with the equi-singularity being imposed on the nodes $\rho_1,\dots,\rho_g$ is  parameterized by a 3-dimensional complex manifold $U$ and  there is a natural isomorphism 
\begin{align}\label{tan3}
T_CU\simeq 
\big\{s\in H^0\big(\nu^*N_{\mu^*C}\big)\set 
s(p'_i) =s(\ol p'_i)=0,\,i=1,\dots,g\big\}.
\end{align}
$U$ can be considered as a complexification of $M_k$ around the point $C$.
Again it is elementary to see that if we further impose the condition $s(x)=s(y)=0$ on \eqref{tan2}, then the space will be 1-dimensional.
Hence, we again obtain $\dim (T_CU\cap T_C(\ol{xy}^*)) = 1$, meaning the required transversality for $M_k\cap \ol{xy}^{*\sigma}$.
\proofend

\medskip
Now we are able to prove the Zoll property of the EW spaces.

\begin{theorem}\label{t:Zoll}
The EW structure on the 3-manifold $M_k \simeq\mathbb S^2\times I$ ($1\le k\le 2^{g-1}$) obtained from the minitwistor space $\ms T$ is spacelike Zoll. Namely, all spacelike geodesics on $M_k$ are simple closed curves.
\end{theorem}

\proof 
We denote $\CP_{g+3}^{*\sigma}$ for the space of real hyperplanes in $\CP_{g+3}$ in which the minitwistor space $\ms T$ is embedded.
For each $1\le i\le g$, recalling that the minitwistor space $\ms T$ is smooth at any point of the sphere $\ms T\us_i$, we put
$$
\ms T_i^{*\sigma}:=
\big\{H\in\CP_{g+3}^{*\sigma}
\set T_p\ms T\subset H {\text{ for some $p\in \ms T_i\us$}}\big\}.
$$
We first show that these are compact subsets of $\CP_{g+3}^{*\sigma}$. 
Let $\ms T^*\subset\CP_{g+3}^*$ be the dual variety of $\ms T$.
Let $I(\ms T)\subset\CP_{g+3}\times\CP_{g+3}^*$ be the incidence variety of $\ms T$,
$I(\ms T)\us$ the real locus of $I(\ms T)$,
 and $\pr_1$ and $\pr_2$ be the projection from 
$I(\ms T)$ to $\CP_{g+3}$ and $\CP_{g+3}^*$ respectively.
The dual variety $\ms T^*$ is exactly $\pr_2(I(\ms T))$, and 
$\ms T_i^{*\sigma} = \pr_2(\pr_1\inv(\ms T\us_i)\cap I(\ms T)\us)$.
Since $\ms T\us_i\simeq \mathbb S^2$ is compact and hence so is $\pr_1\inv(\ms T\us_i)
\cap I(\ms T)\us$
as $\pr_1:I(\ms T)\us\lras\ms T\us$ is an $\RP_g$-bundle,
this means that the image $\ms T_i^{*\sigma}$ is also compact.

We shall show that, for any $x\neq y\in \ms T_{g+1}\us$, 
\begin{align}\label{key1}
(\ms T_1^{*\sigma}\cap \dots\cap \ms T_g^{*\sigma})\cap \ol {xy}^*
= \bigcup_{k=1}^{2^{g-1}} \mf C_{x,y}.
\end{align}
The left-hand side of \eqref{key1} is compact since all $\ms T_i^{*\sigma}$ are compact as above and $\ol{xy}^*$ is also compact.
Hence, once \eqref{key1} is proved, since the right-hand side is a disjoint union by Proposition \ref{prop:disj1}, any $\mf C_{x,y}$ is compact. It then follows from Proposition \ref{p:slg0} that any spacelike geodesic on $M_k$ is a simple closed curve.

To show \eqref{key1}, take any $H\in (\ms T_1^{*\sigma}\cap \dots\cap \ms T_g^{*\sigma})\cap \ol {xy}^*$.
$H$ is a real hyperplane in $\CP_{g+3}$ that is tangent to $\ms T$ at points of $\ms T_1\us,\dots,\ms T_g\us$ and that passes through the two disjoint points $x$ and $y$ of $\ms T\us_{g+1}$.
By Proposition \ref{prop:charW} and Remark \ref{rmk:charW} about the classification of real hyperplanes satisfying these properties, this means that $H\in M_k$ for some $1\le k\le 2^{g-1}$. 
Since $H\in \ol{xy}^*$, from the definition \eqref{slg}, this implies $H\in \mf C_{x,y}$
and hence the inclusion `$\subset$' in \eqref{key1} is proved.

Conversely, take any $1\le k\le 2^{g-1}$ and any $H\in \mf C_{x,y}$.
Then from \eqref{slg}, $H\in \ol{xy}^*$.
Also, since any $H\in M_k$ is tangent to $\ms T$ at points of $\ms T\us_1,\dots,\ms T\us_g$, it follows $H\in \ms T_1^{*\sigma}\cap \dots\cap \ms T_g^{*\sigma}$. Hence $H\in (\ms T_1^{*\sigma}\cap \dots\cap \ms T_g^{*\sigma})\cap \ol {xy}^*$.
Thus, the inclusion `$\supset$' in \eqref{key1} is also proved.
\proofend

\begin{remark}{\em
From Proposition \ref{p:slg0}, the intersection $\mf C_{x,y}= M_k\cap \ol{xy}^{*\sigma}$ is transversal. Hence, all spacelike geodesics on $M_k$ are `reduced'. Namely, they are not multiply covered by nearby geodesics.
}\end{remark}

\section{Automorphisms and the moduli spaces of the Einstein-Weyl structures}\label{s:mod}
In this section, we investigate automorphisms and the moduli space of the EW space $M_k$ ($1\le k\le 2^{g-1}$) discussed so far.
The moduli space can be readily deduced by slightly modifying the argument for the former, so most of this section is devoted to determining the automorphisms.
We note that, as in the case $g=1$, we cannot immediately conclude that any automorphism of the EW space $M_k$ is induced from a holomorphic automorphism of $\ms T$ because the implication ``EW space $\Longrightarrow$ minitwistor'' is not established in the global setting.
For the same reason, we cannot immediately conclude that the EW structure on $M_k$ really deforms even if the complex structure of $\ms T$ is known to deform.

Let $\ms T$ be the minitwistor space from which the EW space $M_k$ is obtained.
Similarly to the case $g=1$ in \cite{HN22}, we will show that there is the following sequence of natural homomorphisms of the automorphism groups:
\begin{align}\label{seqa}
\Aut(M_k) {\lras} \Aut(M_k\uc) {\lras} \Aut(\ms T),%
\end{align}
where $\Aut(M_k)$ and $\Aut(M_k\uc)$ mean the automorphism groups of the EW structure and 
$\Aut(\ms T)$ means the automorphism group of 
the holomorphic structure.

\subsection{The invariance of $\bm I_k\cup\bm I'_k$.}\label{ss:invI}
In this subsection, we show that any automorphism of the EW space $M_k$ preserves the axis $\bm I_k\cup\bm I'_k$. 
This implies that there is a natural homomorphism $\Aut(M_k) {\lras} \Aut(M_k\uc)$ by restriction.
The proof's strategy is similar to that of genus one in that we use a special family of spacelike geodesics in $M_k$. However, we require more effort to understand the behavior of the real circle of the minitwistor lines as a point on the spacelike geodesics moves.


We used the letter $\mf t_1$ to mean the $2^g$-ple covering map from $({\rm J}_{\Sigma})\us_1$ to $({\rm J}_{\Sigma})\us_{\mf o}$, but in the sequel, we use it to mean its restriction to the Seifert surface $\mf S_k$. 
So $\mf t_1:\mf S_k\lras W_{g+1,\,g+1}$ is an unramified double covering.
Thinking points of $\mf S_k\uc=\mf S_k\minus\{\mf o_k,{\mf o}'_k\}$ as hyperplanes in $\CP_{g+2}$ by the natural inclusion $\mf S_k\uc\subset\CP_{g+2}^*$ from Proposition \ref{p:Sei1}, 
for any point $\mf a(\xi+\eta)\in W_{g+1,\,g+1}\minus\{\mf o\}$
so that $\xi,\eta\in \Sigma\us_{g+1}$ and $\eta\neq\tau(\xi)$, we denote $h(\xi,\eta)$ and $h'(\xi,\eta)$ for the hyperplanes corresponding to the two points
$\mf t_1\inv(\mf a(\xi + \eta))$ of $\mf S_k\uc$.
Though the singular points $\mf o_k,\mf o'_k\in\mf S_k$ do not determine the hyperplane uniquely as mentioned at the end of Section \ref{s:Jac}, 
these lifts make sense also at the point $\mf o\in W_{g+1,\,g+1}$ 
by continuity if a continuous path in $W_{g+1,\,g+1}$ through $\mf o$ is given.
We have
$$
h(\xi,\eta)|_{\Sigma} = \sum_{i=1}^g 2p_i + (\xi+\eta)
\qandq
h'(\xi,\eta)|_{\Sigma} = \sum_{i=1}^g 2p'_i + (\xi+\eta)
$$
for some points $p_i,p'_i\in \Sigma\us_i$.
The intersections $h(\xi,\eta)\cap \mathbb D_{g+1}$ and $h'(\xi,\eta)\cap \mathbb D_{g+1}$ are arcs bounded by $\xi$ and $\eta$ unless $\xi=\eta$.
If $\xi=\eta$, then the two hyperplanes $h(\xi,\eta)$ and $h'(\xi,\eta)$ are tangent to $\Sigma$ also at $\xi\in\mathbb D_{g+1}$,
$h(\xi,\eta)\cap \mathbb D_{g+1}=h'(\xi,\eta)\cap \mathbb D_{g+1}=\{\xi\}$ and $h(\xi,\eta),h'(\xi,\eta)\in\ptl \mf S_k$.
By Proposition \ref{p:Sei1} (iii), the tangent points $p_i$ and $p'_i$ ($1\le i\le g$) equal one of the ramification points $r_i,r'_i$ precisely when $\xi=\tau(\eta)$ where $\tau$ is the hyperelliptic involution, and in such a situation, the corresponding point of $\mf S_k$ is one of the singular points $\mf o_k$ or $\mf o'_k$.

\begin{figure}
\includegraphics[height=90mm]{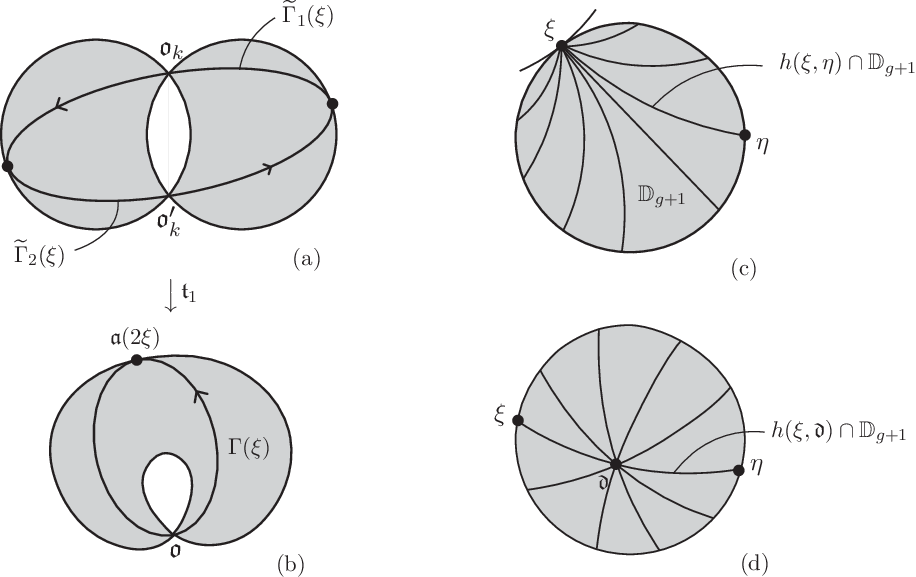}
\caption{The curves in $\mf S_k, W_{g+1,\,g+1}$ and 
the variations of the arcs in $\mathbb D_{g+1}$ }
\label{fig:slg}
\end{figure}

Now take and fix any point $\xi\in\ptl\mathbb D_{g+1}$ and consider the mapping from $\ptl\mathbb D_{g+1}$ to $W_{g+1,\,g+1}$ that sends $\eta\in\ptl\mathbb D_{g+1}$ to $\mf a(\xi+\eta)\in W_{g+1,\,g+1}$. This is injective and from the way of the identification of $W_{g+1,\,g+1}$ as in Figure \ref{fig:Seifert1} given in the proof of Proposition \ref{prop:bdry2}, the image of this mapping is a simple closed curve that satisfies
(i)
it generates the fundamental group $\pi_1(W_{g+1,\,g+1})\simeq\ZZ$,
so that it passes through the origin $\mf o$, and
(ii) it is tangent to the boundary $\ptl W_{g+1,\,g+1}$ at exactly one point $\mf a(2\xi)$.
We denote $\CCC(\xi)$ for this closed curve. 
(See Figure \ref{fig:slg} (b).)
It follows from these that $\tilde\CCC(\xi):=\mf t_1\inv(\CCC(\xi))$, which is of course two-to-one over $\CCC(\xi)$, is a simple closed curve in the Seifert surface $\mf S_k$ which generates $\pi_1(\mf S_k)\simeq\ZZ$ and which is tangent to the boundary $\ptl\mf S_k$ at two points.
Under the inclusion $\mf S_k\uc\subset\CP_{g+2}^*$ as above, these two points correspond to hyperplanes that are tangent to $\Sigma$ not only at a point of each circles $\ptl\mathbb D_1,\dots,\ptl\mathbb D_g$ but also at $\xi\in \ptl\mathbb D_{g+1}$.
We denote the closures of the two connected components of $\tilde\CCC(\xi)\minus\mf t_1\inv(\mf a(2\xi))$ by $\tilde\CCC_1(\xi)$ (the first lap) and $\tilde\CCC_2(\xi)$ (the second lap). So  $\tilde\CCC_1(\xi)\cup\tilde\CCC_2(\xi)= \tilde\CCC(\xi)$
and $\tilde\CCC_1(\xi)\cap\tilde\CCC_2(\xi)= \mf t_1\inv(\mf a(2\xi))$. (See Figure \ref{fig:slg} (a).)

\begin{proposition}\label{p:fol2}
The intersection arcs of the disk $\mathbb D_{g+1}$ with hyperplanes belonging to $\tilde\CCC_1(\xi)$ foliate $\mathbb D_{g+1}$.  Namely, they satisfy 
\begin{align}\label{dis1}
h\cap h'\cap \mathbb D_{g+1} = \{\xi\} 
{\text{ for any $h\neq h'\in \tilde\CCC_1(\xi)$}}
\end{align}
and 
\begin{align}\label{fill1}
\bigcup_{h\in \tilde\CCC_1(\xi)}
h\cap \mathbb D_{g+1} = \mathbb D_{g+1}.
\end{align}
The same thing holds for $\tilde\CCC_2(\xi)$, so the whole $\tilde\CCC(\xi)$ foliates $\mathbb D_{g+1}$ twice.
\end{proposition}

\proof
We keep the above notations.
Since any $h\in \tilde\CCC_1(\xi)$ is tangent to the branch curve $\Sigma$ at points of the circles $\ptl\mathbb D_1,\dots,\ptl\mathbb D_g$, 
for any $h'\in \tilde\CCC_1(\xi)$ which is sufficiently close to $h$, 
the sections $h\cap {\rm C}(\Lmd)$ and $h'\cap {\rm C}(\Lmd)$ meet at a (real) point belonging to a neighborhood of $\ptl\mathbb D_i$ in ${\rm C}(\Lmd)$ for any $1\le i\le g$.
Since $\deg {\rm C}(\Lmd) = \deg \Lmd = g+1$ and $h$ and $h'$ pass through $(g+1)$-th point $\xi\in\ptl\mathbb D_{g+1}$,
this implies  
$$
h\cap h'\cap \mathbb D_{g+1} = \{\xi\}.
$$
If $h\in \tilde\CCC_1(\xi)$ approaches the ends of $\tilde\CCC_1(\xi)$, then the point of $h\cap\ptl\mathbb D_{g+1}$ which is residual to $\xi$ approaches $\xi$ from each side depending on the choice of the end and the arc $h\cap\mathbb D_{g+1}$ shrinks to the single point $\xi$.
Since $\mathbb D_{g+1}$ is a disk, this implies the global disjointness \eqref{dis1} as well as the exhaustion property \eqref{fill1}. (See Figure \ref{fig:slg} (c).)
The case of $\tilde\CCC_2(\xi)$ is obtained in the same way.
\proofend

\medskip
Still fixing any point $\xi\in\ptl\mathbb D_{g+1}$,
from Proposition \ref{p:fol2}, for any point $\mf d\in\mathbb D_{g+1}\minus\ptl\mathbb D_{g+1}$, there uniquely exist $h_1\in\tilde\CCC_1(\xi)$ and $h_2\in\tilde\CCC_2(\xi)$ that pass through $\mf d$.
We denote these hyperplanes $h_1(\xi,\mf d)$ and $h_2(\xi,\mf d)$ respectively.
Moving the point $\xi\in\ptl\mathbb D_{g+1}$ one way around, we obtain two families of hyperplanes in $\CP_{g+2}$ through $\mf d$, which are parameterized by $\ptl\mathbb D_{g+1}$.
Recall that the double covering $\Pi:\ms T\lras {\rm C}(\Lmd)$ restricts to a double covering $\ms T\us_{g+1}\lras\mathbb D_{g+1}$ with branch $\ptl\mathbb D_{g+1}$, so $\Pi\inv(\mf d)$ consists of two points if $\mf d\in \mathbb D_{g+1}\minus\ptl\mathbb D_{g+1}$.
Let $x,y$ be these points and $\mf C_{x,y}$ the spacelike geodesic determined by $\mf d\in D_{g+1}\minus\ptl\mathbb D_{g+1}$, namely the set of minitwistor lines passing through $x$ and $y$.

\begin{proposition}\label{p:slg3}
For $i=1,2$, the assignment $\xi\longmapsto h_i(\xi,\mf d)$ gives a bijection from $\ptl\mathbb D_{g+1}$ to $\mf C_{x,y}.$
\end{proposition}

\proof We keep the above notations.
In a similar way as the proof for the disjointness \eqref{dis1}, using the point $\mf d$ instead of $\xi$, we obtain that, if distinct two points $\xi,\xi'$ on $\ptl\mathbb D_{g+1}$ are sufficiently close, then $h_i(\xi,\mf d)\cap h_i(\xi',\mf d)\cap \mathbb D_{g+1} = \{\mf d\}$.
Denoting $\eta$ for the residual point $(h_i(\xi,\mf d)\cap\ptl\mathbb D_{g+1})\minus\{\xi\}$, this means that if we move a point $\xi$ one way around on $\ptl\mathbb D_{g+1}$, then another point $\eta$ moves one way around in the same direction while $\xi$ and $\eta$ stay `symmetric' with respect to $\mf d$. Namely, if we choose a homeomorphism $\mathbb D_{g+1}\lras \{|z|\le 1\}$ such that $\mf d\longmapsto 0$, then we may assume $\xi(\theta) = e^{i\theta}$ and $\eta(\theta) = e^{i(\theta+\pi)}$.
Hence, the family $\{h_i(\xi,\mf d)\cap \mathbb D_{g+1}\set \xi\in \ptl\mathbb D_{g+1}\}$ of arcs fills the disk $\mathbb D_{g+1}$ twice. Namely, for any point $\ddd\in \mathbb D_{g+1}\minus
(\ptl\mathbb D_{g+1}\cup\{\mf d\})$, there exist precisely two $\xi\in\ptl\mathbb D_{g+1}$ such that $\ddd\in h_i(\xi,\mf d)$, and further, $h_i(\xi,\mf d)\cap h_i(\xi',\mf d)\cap \mathbb D_{g+1} = \{\mf d\}$ for any $\xi\neq\xi'\in\ptl\mathbb D_{g+1}$.
(See Figure \ref{fig:slg} (d).)
%

Again from the identification of the region $W_{g+1,\,g+1}$ as a `contraction' of the M\"obius strip, the above behavior of the two points $\xi$ and $\eta$ implies that the locus of the points $\mf a(\xi+\eta)$ in the above move is a composition of two simple smooth curves in $W_{g+1,\,g+1}$ each of which generates $\pi_1\simeq\ZZ$ and does not meet $\ptl W_{g+1,\,g+1}$. (Namely, the point goes around twice when $\xi\in\ptl\mathbb D_{g+1}$ does once.)
This implies that the locus of the lift $h_i(\xi,\mf d)\in\mf S_k$ as $\xi\in\ptl\mathbb D_{g+1}$ moves one way around is exactly a simple closed smooth curve in $\mf S_k$ which generates $\pi_1\simeq\ZZ$ and which does not meet $\ptl\mf S_k\minus\{\mf o_k,\mf o'_k\}$.
So the family $\{h_i(\xi,\mf d)\set \xi\in \ptl\mathbb D_{g+1}\}$ is effectively parameterized by $\xi\in\ptl\mathbb D_{g+1}\simeq \mathbb S^1$.
The hyperplane section $\ms T\cap \Pi\inv\big(h_i(\xi,\mf d)\big)$
is a minitwisor line passing through $x$ and $y$.
Hence, the assignment $\xi\longmapsto h_i(\xi,\mf d)$ gives an injection $\ptl\mathbb D_{g+1}\lras\mf C_{x,y}$.
If this injection were not surjective, then there would exist a hyperplane $H'\subset\CP_{g+3}$ which passes through $x,y$ but not through the center of $\Pi$. 
But such an $H'$ does not exist because it has to be mapped isomorphically to $\CP_{g+2}$ by $\Pi$ while $x,y$ are mapped to the same point $\mf d$.
\proofend

\medskip
Thus, both for $i=1,2$, the mapping $\xi\longmapsto h_i(\xi,\mf d)$ gives the effective parameterization $\ptl\mathbb D_{g+1}\simeq\mf C_{x,y}$.
Adapting any $i\in\{1,2\}$ and
writing $h(\xi,\mf d):=h_i(\xi,\mf d)$
and $H(\xi,\mf d) := \Pi\inv(h_i(\xi,\mf d))$, we obtain
\begin{align}\label{slg3}
\mf C_{x,y} &= \{H(\xi,\mf d)\set \xi\in \ptl\mathbb D_{g+1}\}\subset\CP_{g+3}^*\\
&\simeq \{h(\xi,\mf d)\set \xi\in \ptl\mathbb D_{g+1}\}
\subset \mf S_k\uc 
\subset \CP_{g+2}^*\label{slg3'}.
\end{align}
As above, this parameterization is effective in the sense that any point of $\mf C_{x,y}$ is written $H(\xi,\mf d)=\Pi\inv (h(\xi,\mf d))$ for a unique $\xi\in\ptl\mathbb D_{g+1}$. 
These geodesics are determined from any point $\mf d\in\mathbb D_{g+1}\minus\ptl\mathbb D_{g+1}$, and further, there is freedom for choosing the center of the projection $\Pi$.
Thus, we have obtained a 3-dimensional family of spacelike geodesics on $M_k$.
All these geodesics pass through $\bm I_k$ and $\bm I'_k$.

We will soon require a description of the (natural) complexification of these geodesics.
For this, let $V$ be a real open tubular neighborhood of $\ptl\mathbb D_{g+1}=\Sigma\us_{g+1}$ in the hyperelliptic curve $\Sigma$ and put $W_{g+1,\,g+1}^{\CC}:=\{\mf a(\xi+\eta)\in {\rm J}_{\Sigma}\set \xi,\eta\in V\}$.
The latter is a $\sigma$-invariant complex surface in the Jacobian ${\rm J}_{\Sigma}$ that contains the region $W_{g+1,\,g+1}$ as its real locus.
Since the double covering map $\mf t_1:\mf S_k\lras W_{g+1,\,g+1}$ was a restriction of the mapping $\mf t:{\rm J}_{\Sigma}\lras {\rm J}_{\Sigma}$ defined by $\mf a(D)\longmapsto\mf a(2D)$ and this is also a covering map, if we choose $V$ sufficiently thin, then $\mf t\inv(W_{g+1,\,g+1}^{\CC})$ has a connected component which is two-to-one over $W_{g+1,\,g+1}^{\CC}$.
We denote $\mf S_k^{\CC}$ for this component and call it a complexification of $\mf S_k$.
Since $\mf C_{x,y}\subset\mf S_k$ from \eqref{slg3'},
the geodesic $\mf C_{x,y}$ can be regarded as included in $\mf S_k^{\CC}$.

In the same way as $\mf S_k\minus\{\mf o_k,\mf o'_k\}$, any point of $\mf S_k^{\CC}\minus\{\mf o_k,\mf o'_k\}$ determines $(g+2)$ points $p_1,\dots, p_g$ and $\xi,\eta$ on the curve $\Sigma$ by the same equation 
$2\mf a(p_1+\dots+p_g) + \mf a(\xi+\eta) = \mf o$ and these points determine the hyperplane $h\subset\CP_{g+2}$ satisfying $h|_{\Sigma} = 2(p_1+\dots+p_g)+ (\xi+\eta)$.
Hence, the inclusion $\mf S_k\minus\{\mf o_k,\mf o'_k\}\subset\CP_{g+2}^*$ extends naturally to $\mf S_k^{\CC}\minus\{\mf o_k,\mf o'_k\}\subset\CP^*_{g+2}$. 
Under this inclusion, the dual hyperplane $\mf d^*\subset\CP_{g+2}^*$ (consisting of hyperplanes passing through $\mf d\in\mathbb D_{g+1}\minus\ptl\mathbb D_{g+1}$) intersects the complex surface $\mf S_k^{\CC}\minus\{\mf o_k,\mf o'_k\}$
in a holomorphic curve and it contains the real curve $\mf C_{x,y}\minus\{\mf o_k,\mf o'_k\}$.
We denote $\mf C_{x,y}^{\CC}$ for the closure of this curve in $\CP_{g+2}^*$.
We use the same symbol to mean the image of this curve under the injection $\Pi^*:\CP_{g+2}^*\lras\CP_{g+3}^*$.
This is the complexification of the geodesic $\mf C_{x,y}$. 
In the same way as its real locus $\mf C_{x,y}$, each element of $\mf C_{x,y}^{\CC}$ may be written $H(\xi,\mf d)=\Pi\inv (h(\xi,\mf d))$ for some point $\xi\in V$ and a hyperplane $h(\xi,\mf d)$ through $\xi,\mf d$ and it defines a complex minitwistor line passing through $\mf d$ as $H(\xi,\mf d)\cap \ms T$.
Since the real geodesic $\mf C_{x,y}$ is smooth, so is the holomorphic curve $\mf C_{x,y}^{\CC}$.

Using these geodesics and their complexification, we shall show:

\begin{proposition}\label{p:C0C1}
Let $C_0$ and $C_1$ be any irregular and regular minitwistor lines on $\ms T$ respectively. Then there exists a real open subset $U$ of $\ms T$ that contains both $C_0$ and $C_1$, such that the modification as in Section \ref{ss:modify} can be applied.
\end{proposition}

\proof
We use the above notations.
Using the $\mathbb S^1$-action as in the proof of Proposition \ref{prop:charW}, 
we may suppose that $C_1=\ms T\cap \Pi\inv(h_1)$ for some $h_1\in \mf S_k\uc$.
Let $h_0\subset\CP_{g+2}$ be the real hyperplane such that $C_0 = \ms T\cap \Pi\inv(h_0)$ and $\lmd\in K_{g+1}$ be the point such that the generating line $l_{\lmd}$ over it is a component of $C_0$.
Put $\mf d:=h_1\cap l_{\lmd}$, which is necessarily a point belonging to the interior of the disk $\mathbb D_{g+1}$.
Then from the parameterization \eqref{slg3'}, 
there uniquely exist points $\xi_1,\xi_0\in\ptl\mathbb D_{g+1}$
satisfying $h_1 = h(\xi_1,\mf d)$ and $h_0 = h(\xi_0,\mf d)$.
These are one of the two points $h_1\cap\ptl\mathbb D_{g+1}$ and $h_0\cap\ptl\mathbb D_{g+1}$ respectively.
Let $\eta_i\in\ptl\mathbb D_{g+1}$ $(i=0,1)$ be the point of $h_i\cap\ptl\mathbb D_{g+1}$
which is residual to $\xi_i$.
Then $\eta_0 = \tau(\xi_0)$ and $\eta_1\neq\tau(\xi_1)$.
By choosing one of the two orientations of the circle $\ptl\mathbb D_{g+1}$, we may suppose that the four points $\xi_0,\xi_1,\eta_0,\eta_1$ are placed in this order.

Let $ A(\xi_0,\xi_1)$ be the closed arc in $\ptl\mathbb D_{g+1}$  whose ends are $\xi_0,\xi_1$ and which does not contain $\eta_0$ and therefore $\eta_1$ also.
For a point $\xi\in A(\xi_0,\xi_1)$, the minitwistor line $\ms T\cap H(\xi,\mf d)$ is irregular only when $\xi=\xi_0$ and if $\xi\neq\xi_0$, then for any $1\le i\le g$, the tangent point $p_i(\xi):=\mathbb D_i\cap h(\xi,\mf d) \,(= \ptl\mathbb D_i\cap h(\xi,\mf d))$ is not a ramification point $r_i$ nor $r'_i$ from Proposition \ref{p:Sei1}. 
This means that for any distinct points $\xi\neq\xi'\in  A(\xi_0,\xi_1)$ and any $1\le i\le g$, the intersection of the two hyperplane sections $h(\xi,\mf d)\cap {\rm C}(\Lmd)$ and $h(\xi',\mf d)\cap {\rm C}(\Lmd)$ has a point belonging to $\pi\inv(\ol K_i)$, where as before $\pi$ is the cone projection.
Since $\deg {\rm C}(\Lmd) = g+1$, this means that 
these points on $\pi\inv(\ol K_i)$ ($1\le i\le g$) and $\mf d$ are all points of $h(\xi,\mf d)\cap h(\xi',\mf d)\cap {\rm C}(\Lmd)$.


For any $1\le i\le g$, the tangent point $p_i(\xi)\in \ptl\mathbb D_i$ constitutes a closed arc as $\xi$ moves on $ A(\xi_0,\xi_1)$ and we denote $ A_i$ for it.
One of the ends of $ A_i$ is a ramification point $r_i$ or $r'_i$ and since the arc $ A(\xi_0,\xi_1)$ does not contain the point $\eta_0$, the arc $ A_i$ does not contain another ramification point.
The image $\pi( A_i)$ is also a closed arc in the closure $\ol K_i$.
Since $h(\xi,\mf d)\cap {\rm C}(\Lmd)$ and $h(\xi',\mf d)\cap {\rm C}(\Lmd)$ intersect at a point over $\ol K_i$ for each $1\le i\le g$ and any $\xi\neq \xi'\in  A(\xi_0,\xi_1)$ as above,
from the continuity of the cone projection $\pi$, for any open tubular neighborhoods 
$T_1,\dots, T_g$ in $\Lmd$ of the arcs $\pi( A_1),\dots, \pi( A_g)$ respectively, there exists an open tubular neighborhood $V(\xi_0,\xi_1)$ of $ A(\xi_0,\xi_1)$ in $\Sigma$, such that $h(\xi,\mf d)\cap h(\xi',\mf d)\cap {\rm C}(\Lmd)$ has a point belonging to $\pi\inv(T_i)$ for any $1\le i\le g$ and any $\xi\neq\xi'\in V(\xi_0,\xi_1)$. 
By letting $T_1,\dots,T_g$ be mutually disjoint, for the same reason to the case $\xi,\xi'\in  A(\xi_0,\xi_1)$, this implies that 
these $g$ points and $\mf d$ are all points of the last intersection.
In the following, we always take $T_1,\dots, T_g$ and $V(\xi_0,\xi_1)$ to be real and the closures satisfy $\ol T_i\cap \ol T_j=\emptyset$ if $i\neq j$. 

In the following we denote $f:=\pi\circ\Pi$ for simplicity, which gives $\ms T$ the (rational) conic bundle structure over $\Lmd$. 
For a point $\xi\in V(\xi_0,\xi_1)$, we denote $C(\xi,\mf d) := \ms T\cap H(\xi,\mf d)$, which is a (real) minitwistor line if $\xi\in A(\xi_0,\xi_1)$.
So its real circle $C(\xi,\mf d)^{\RR}$ is $\ms T\us_{g+1}\cap \Pi\inv(h(\xi,\mf d))$.
Note that $C(\xi_i,\mf d) = C_i$ for $i=0,1$.
As $C(\xi,\mf d)\cap C(\xi',\mf d)=\Pi\inv(h(\xi,\mf d)\cap h(\xi',\mf d)\cap {\rm C}(\Lmd))$, from the above choice of the neighborhoods $V(\xi_0,\xi_1)$ and $T_1,\dots, T_g$, if 
$\xi,\xi'\in V(\xi_0,\xi_1)$, then 
$C(\xi,\mf d)\cap C(\xi',\mf d)$ consists of two points from each of $f\inv(T_1),\dots,f\inv(T_g)$ and $x,y\in\ms T_{g+1}\us$.
The intersections in each $f\inv(T_i)$ are transverse.

Under these setting, taking a real open neighborhood $F$ of $\ol K_{g+1}$ in $\Lmd$ 
such that $(\ol T_1\cup\dots\cup \ol T_g)\cap \ol F = \emptyset$,
we shall define a subset $U$ of $\ms T$ by 
\begin{align}\label{}
U:= \Big(\bigcup_{\xi\in V(\xi_0,\xi_1)} C(\xi,\mf d)\Big)\cup f\inv(F).
\end{align}
This is a real connected open subset of $\ms T$.
For each $1\le i\le g$, let $U_i$ be the subset of $f\inv(T_i)$ which consists of the node of $C(\xi,\mf d)\cap f\inv(T_i)$ itself as $\xi$ moves in $V(\xi_0,\xi_1)$ and the following subset of $f\inv(T_i)$:
\begin{align*}
\bigcup_{\xi,\,\xi'\in V(\xi_0,\xi_1),\,\,\xi\neq\xi'}C(\xi,\mf d)\cap C(\xi',\mf d)\cap f\inv(T_i)
\end{align*}
$U_i$ is also a real connected open subset of $\ms T$ and is clearly included in $U$.
From the transversality for the intersection $C(\xi,\mf d)\cap C(\xi',\mf d)$ in each $f\inv(T_i)$, 
as we did in Section \ref{ss:modify}, by ungluing $U$ over $U_1\cup\dots\cup U_g$, we obtain the ``normalization'' $\hat U$ of $U$ which resolves all nodes of $C(\xi,\mf d)$ in $f\inv(T_1\cup\dots\cup T_g)$ as well as the transversal intersection $C(\xi,\mf d)\cap C(\xi',\mf d)\cap f\inv(T_i)$ for distinct $\xi,\xi'\in V(\xi_0,\xi_1)$.
Since $U$ contains both $C_0$ and $C_1$, this means that $U$ admits the modification as given in Section \ref{ss:modify}.
%
%
\proofend

%

\begin{corollary}\label{c:aut1}
For any two points $p_1\in M_k\uc$ and $p_0\in \bm I_k\cup\bm I'_k$,
there exist no neighborhoods of $p_1$ and $p_0$ in $M_k$ on which the EW structures are isometric to each other.
Therefore,
any automorphism of the EW space $M_k$ preserves the axis $\bm I_k\cup\bm I'_k$, therefore also $M_k\uc = M_k\minus
(\bm I_k\cup\bm I'_k)$.
\end{corollary}

\proof
Let $C_1$ and $C_0$ be the minitwistor lines corresponding to $p_1$ and $p_0$ respectively.
From Proposition \ref{p:C0C1}, we can take a real open subset $U\subset \ms T$ that contains both $C_1$ and $C_0$ and for which the modification in Section \ref{ss:modify} may be applied. Let $\check U$ be the modified minitwistor space and $\check C_1, \check C_0\subset\check U$ the transformations of $C_1$ and $C_0$ respectively. From Proposition \ref{p:tmtl}, both $\check C_1$ and $\check C_0$ are smooth minitwistor lines, and for the same reason to Proposition \ref{p:le1}, they satisfy $\dim |\check C_1|=0$ and  $\dim |\check C_0|=1$. This means that there is no holomorphic isomorphism between any neighborhoods of $\check C_1$ and $\check C_0$ in $\check U$. Since the EW structure on $M_k$ around the points $p_1$ and $p_0$ are uniquely determined from the complex structures on neighborhoods of these points respectively, this means that there is no isomorphism between any neighborhoods of $p_1$ and $p_0$. 
\proofend

\medskip
Thus, we have obtained the homomorphism $\Aut(M_k)\lras \Aut(M_k\uc)$ as the restriction of the domain. In the next two subsections, we show that there is a natural homomorphism $\Aut(M_k\uc) {\lras} \Aut(\ms T)$.

\subsection{The limits of minitwistor lines in future and past infinity}\label{ss:limtl}

As before, suppose that $M_k$ ($1\le k\le 2^{g-1}$) is associated to a minitwistor space $\ms T$ of genus $g$ and let
$\ms H_k$ and $\hat{\ms H_k}$ be the future and past infinity of $M_k$. (See Definition \ref{d:infty}. We do not make distinction between $\ms H_k$ and $\hat{\ms H}_k$ yet.)
So $\ol M_k= M_k \cup \ms H_k\cup\hat{\ms H_k}\simeq \mathbb S^2 \times \ol I$ where $\ol I$ is a closed interval.
Both $\ms H_k$ and $\hat{\ms H_k}$ are $\mathbb S^1$-invariant and 
their intersections with the axis $\bm I_k\cup\bm I'_k$ are exactly the $\mathbb S^1$-fixed points on them.


As before, let $\rho_i\in\{r_i,r'_i\}$ ($1\le i\le g$) with $\rho_1=r_1$ be the choices that correspond to the boundary data $k$
and we denote $l_i=l(\rho_i)$ and $l'_i=l(\rho'_i)$.
The irregular minitwistor lines that correspond to the two ends of the intervals $\bm I_k$ and $\bm I'_k$ are respectively
\begin{align}\label{lines1}
\sum_{1\le i\le g}(l_i+\ol l_i) + \big\{l(r_{g+1}) + \ol{l(r_{g+1})}\big\},\,\,
\sum_{1\le i\le g}(l_i+ \ol{l_i}) +
\big\{l(r'_{g+1}) + \ol{l(r'_{g+1})}\big\}
\end{align}
and 
\begin{align}\label{lines2}
\sum_{1\le i\le g}(l'_i+ \ol{l'_i}) +
\big\{l(r'_{g+1}) + \ol{l(r'_{g+1})}\big\},
\,\,
\sum_{1\le i\le g}(l'_i + \ol l'_i)
+\big\{l(r_{g+1}) + \ol{l(r_{g+1})}\big\}.
\end{align}
Since we had the linear equivalence $\sum_{1\le i\le g} \rho_i + r_{g+1}\sim \sum_{1\le i\le g} \rho'_i + r'_{g+1}$ on the hyperelliptic curve $\Sigma$ as in Proposition \ref{p:123},
from the explicit description of the intersecting pair as in Proposition \ref{p:inter3},
the first divisors in \eqref{lines1} and \eqref{lines2} belong to the same component among $\ms H_k$ and $\hat{\ms H}_k$, and the second divisors belong to another component.
By exchanging $\ms H_k$ and $\hat{\ms H}_k$ if necessary, 
we may suppose that the first divisors belong to $\ms H_k$.
We decompose the first divisor in \eqref{lines1} 
as $L_k + \ol L_k$ where $L_k:=\sum_{1\le i\le g} l_i + l(r_{g+1})$
and the second divisor in \eqref{lines1}
as $\hat L_k+\ol{\hat L_k}$ where $\hat L_k = \sum_{1\le i\le g} l_i + l(r'_{g+1})$.
Similarly, we decompose the first divisor in \eqref{lines2} 
as $L'_k + \ol L'_k$ where $L'_k:=\sum_{1\le i\le g} l'_i + l(r'_{g+1})$ and
the second one in \eqref{lines2} as
$\hat L'_k+\ol{\hat L'_k}$ where $\hat L'_k = \sum_{1\le i\le g} l'_i + l(r_{g+1})$.

All the divisors $L_k$, $\hat L_k$, $L'_k$, $\hat L'_k$ and their conjugations are Caritier divisors on $\ms T$ and they are $\CC^*$-invariant.
It is easy to see that the complete linear systems
generated by these divisors are base point free pencils whose general member is a smooth rational curve, and also 
\begin{align}\label{penc}
|L_k| = \big|\ol {L'_k}\big| \qandq |\hat L_k| = \big|\ol{{\hat L'}}_k\big|.
\end{align}
We denote $$\varpi_k:\ms T\lras Q_k:=\qdr$$ be the morphism which is the product of the two morphisms $\ms T\lras\CP_1$ induced by the pencils $|L_k|$ and $|\ol L_k|$.
Using the pencils $|\hat L_k|$ and $|\ol{\hat L_k}|$, we obtain another morphism $\hat\varpi_k:\ms T\lras \hat Q_k:=\qdr$.
These mappings have the following nice properties, which will be used in the sequel.
We omit the proof since they are readily obtained from the explicit construction.
\begin{proposition}\label{p:cov1}
The mappings $\varpi_k$ and $\hat\varpi_k$ satisfy the following:
\begin{itemize} 
\setlength{\leftskip}{-15pt}
\setlength{\itemsep}{-.5mm} 
\item
they are finite, of degree $(g+1)$ and preserve the real structure and the $\CC^*$-action,
\item in inhomogeneous coordinates $(z,w)$ on $Q_k$ and $\hat Q_k$, 
the induced $\CC^*$-action and real structure are given by $(z,w)\longmapsto (tz,tw)$ and  $(z,w)\longmapsto (\ol w\inv,\ol z\inv)$,
\item the branch divisor consists of $(g+1)$ number of real irreducible $(1,1)$-curves through the two points $(0,0)$ and $(\infty,\infty)$ in the above coordinates and they are disjoint from the real locus $Q_k\us$ and $\hat Q_k\us$. In particular, $\varpi_k$ and $\hat\varpi_k$ restrict to diffeomorphisms from $\ms T_i\us\simeq \mathbb S^2$ onto $Q_k\us\simeq \mathbb S^2$ and $\hat Q_k\us\simeq \mathbb S^2$ respectively for any $1\le i\le g+1$,
\item $\varpi_k^*\ms O(1,1)\simeq \hat\varpi_k^*\ms O(1,1) \simeq \ms O_{\CP_{g+3}}(1)|_{\ms T}$,
\item $\varpi_k(\bm p_{\infty}) = \hat\varpi_k(\bm p_{\infty}) = (0,0)$ and $\varpi_k(\ol{\bm p}_{\infty}) =\hat\varpi_k(\ol{\bm p}_{\infty}) = (\infty,\infty)$, 
\item $\{\rho_1,\dots, \rho_g,r_{g+1}\}\stackrel{\varpi_k}\lras
(0,\infty)$ and $\{\rho'_1,\dots, \rho'_g,r'_{g+1}\}\stackrel{\varpi_k}\lras
(\infty,0)$, 
\item $\{\rho_1,\dots, \rho_g,r'_{g+1}\}\stackrel{\hat\varpi_k}\lras
(0,\infty)$ and $\{\rho'_1,\dots, \rho'_g,r_{g+1}\}\stackrel{\hat\varpi_k}\lras
(\infty,0)$.
%
\end{itemize}
\end{proposition}

Pulling back the pencils of $A$-lines and $B$-lines on the quadrics $Q_k$ and $\hat Q_k$, we obtain four pencils on $\ms T$.
The two pencils obtained from $Q_k$ (resp.\,$\hat Q_k$) are generated by $L_k\sim \ol L'_k$ and $\ol L_k\sim L'_k$
(resp.\,$\hat L_k\sim \ol{\hat L'}_k$ and $\ol {\hat L_k}\sim {\hat L'_k}$) and
all other members of the pencils are smooth rational curves.

\begin{remark}\label{r:g1}{\em
The mappings $\varpi_k$ and $\hat\varpi_k$ can be thought as natural generalizations of 
the double covering maps $\pi_2$ and $\pi_4$ used in \cite{HN22}.
Note that if $g=1$, then the boundary data $k$ can only be 1 as $2^{g-1}=1$.
\proofend
}
\end{remark}

\begin{proposition}\label{p:limtl}
For any point $q\in Q_k\us$, let $\ell_q$ and $\ol \ell_q$ be the pair of lines on the quadric $Q_k$ that intersect at $q$. Then the assignment 
$$
q\longmapsto \varpi_k\inv\big( \ell_q + \ol\ell'_q\big)
$$
gives a bijection from $Q_k\us$ to the boundary $\ms H_k$.
Similarly, there is a natural bijection from $\hat Q_k\us$ to $\hat{\ms H}_k$ using $\hat\varpi_k:\ms T\lras \hat Q_k$.
\end{proposition}

\proof 
From the previous proposition, for any point $q\in Q\us_k$, we can write $\varpi_k\inv(q) = \varpi_k\inv(\ell_q)\cap \varpi_k\inv(\ol \ell_q) = \{p_1,\dots, p_{g+1}\}$, where $p_i\in \ms T\us_i$ for any $i$.
Again from the previous proposition,
the branch divisor of $\varpi_k$ does not hit $Q\us_k$, so 
every $p_i$ is an ordinary node of the curve $\varpi_k\inv(\ell_q)+ \varpi_k\inv(\ol \ell_q)$ and this is a hyperplane section of $\ms T$.
From Proposition \ref{prop:charW} which characterizes the divisors belonging to $\ms H_k$ and $\hat{\ms H}_k$, this means that $\varpi_k\inv(\ell_q)+ \varpi_k\inv(\ol \ell_q)$ belongs to either $\ms H_k$ or $\hat{\ms H}_k$.
But since $L_k + \ol L_k\in \ms H_k$ from the definition of $\ms H_k$, we obtain $\varpi_k\inv(\ell_q)+ \varpi_k\inv(\ol \ell_q)\in \ms H_k$. 
Similarly, $\hat\varpi_k\inv(\ell_q)+ \hat\varpi_k\inv(\ol \ell_q)\in \hat{\ms H}_k$.
These mappings $Q_k\us\lras\ms H_k$ and $\hat Q_k\us\lras\hat{\ms H}_k$ are bijective again by Proposition \ref{prop:charW}.
\proofend

\medskip
If a point of the EW space $M_k$ approaches a point $q$ of the boundary component $\ms H_k$ (resp.\,$\hat{\ms H}_k$),
then the real circle of the corresponding minitwistor line, which lies in $\ms T_{g+1}\us\simeq\mathbb S^2$, shrinks to the node
of the reducible curve
$\varpi_k\inv(\ell_q)+ \varpi_k\inv(\ol\ell_q)$
(resp.\,$\hat\varpi_k\inv(\ell_q)+ \varpi_k\inv(\ol\ell_q)$)
belonging to $\ms T\us_{g+1}$.
So this decomposition can be regarded as a limit of the decomposition of a minitwistor line to two disks divided by the real circle.

\subsection{The families of disks in $\ms T$}\label{ss:disk}
We still fix any boundary data $1\le k\le 2^{g-1}$ so that the ramification points $\rho_1=r_1$ and $\rho_2,\dots,\rho_{g}$ as well.
Let $\ms C_k\uc\subset \ms T\times M_k\uc$ be the tautological family over $M_k\uc$ whose fiber over a point is the corresponding hyperplane sections of $\ms T$.
As before, we write $C^{\RR}=C\cap \ms T\us_{g+1}$ for the real circle of $C\in M_k\uc$. These circles constitute an $\mathbb S^1$-subbundle $(\ms C\uc_k)^{\RR}\subset\ms C\uc_k$.
Let $(\ms C\uc_k)'\lras \ms C\uc_k$ be a simultaneous normalization for all $C\in M_k\uc$.  
This desingularizes all nodes of any $C\in M_k\uc$ and the resulting fiber is a smooth $\CP_1$. 
The $\mathbb S^1$-subbundle can be naturally regarded as a subset of $(\ms C\uc_k)'$ and if we remove it from $(\ms C\uc_k)'$, then we obtain a fiber bundle over $M_k\uc$ whose fibers are two open disks.
By considering the natural extension of this bundle to the completion $M_k$, it is easy to see that this consists of two disk bundles.
We denote $(\ms C\uc_k)^+$ and $(\ms C\uc_k)^-$ for the closed disk bundles which are obtained from the two open disk bundles over $M_k\uc$ by attaching back the circle bundle $(\ms C\uc_k)^{\RR}$ to each.
To distinguish these two bundles, we need the following lemma.
It is a generalization of \cite[Lemma 5.1]{HN22} and can be proved in a similar way using the projections $\varpi_k$ and $\hat\varpi_k$ and Proposition \ref{p:limtl}. So we omit a proof. As before we denote $l_i= l(\rho_i)$ and so on.
These lines are determined from the boundary data $k$.

\begin{lemma}\label{l:pm}
Exactly one of the following two situations occurs:
\begin{itemize}
\setlength{\itemsep}{-.5mm}
\item[\em (i)] every disk $C^+\subset ({\ms C}_k\uc)^+$ intersects 
$l_i$ and $\ol l'_i$ for all $1\le i\le g$ and intersects none of 
$\ol l_i$ and $l'_i$ for $1\le i\le g$,
\item[\em (ii)] every disk $C^+\subset ({\ms C}_k\uc)^+$ intersects $\ol l_i$ and $l'_i$ for all $1\le i\le g$ and intersects none of 
$l_i$ and $\ol l'_i$ for $1\le i\le g$.
\end{itemize}
\end{lemma}

\begin{figure}
\includegraphics{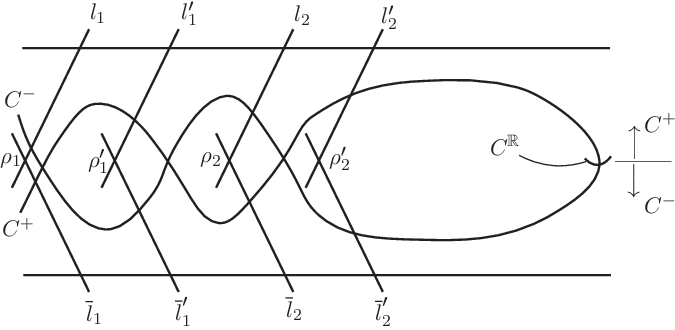}
\caption{
The disks $C^+$ and $C^-$ ($g=2$)}
\label{f:Cpm}
\end{figure}

Let $C\in M_k\uc$ be a regular minitwistor line and $C= C^+\cup C^-$ the decomposition into closed disks by the real circle $C^{\RR}$,
where $C^+\subset (\ms C\uc_k)^+$ and $C^-\subset (\ms C\uc_k)^-$.
From the previous lemma, we may make the following definition.
(See Figure \ref{f:Cpm}.)

\begin{definition}\label{d:Tplus}
{\em In the following,
we make a distinction between the disks $C^+$ and $C^-$ by the following rule.
$C^+$ is the disk that intersects all of $\ol l_1,\ol l_2,\dots,\ol l_g$
(and hence all of $l'_1,l'_2,\dots,l'_g$ as well).
So, $C^-$ is the disk that intersects all of $l_1,l_2,\dots, l_g$
(and hence all of $\ol l'_1,\ol l'_2,\dots,\ol l'_g$ as well).
Further, we define two open complex surfaces by
$$
\ms T^{k+}:= \ms T\setminus \bigcup_{1\le i\le g} \big(l_i\cup \ol l'_i\big)
\qandq 
\ms T^{k-} := \sigma\big(\ms T^{k+}\big).
$$
}
\end{definition}

\begin{proposition}\label{p:disk1}
For any point $z$ of the surface $\ms T^{k+}$, there exists a point of $M_k\uc$ whose corresponding regular minitwistor line $C$ satisfies $z\in C^+$.
In other words, every point of $\ms T^{k+}$ is passed by some closed disk belonging to $(\ms C\uc_k)^+$.
\end{proposition}

This is a generalization of the case $g=1$ given in \cite[Lemma 5.6]{HN22}, and as we did there, we prove the case $z\in \ms T\us$ and the case $z\not\in \ms T\us$ separately. The latter case can be proved in a similar way to the case $g=1$ based on the description of degenerations of minitwistor lines given in Section \ref{ss:limtl}, but in the former case we have to argue in a different way.

\medskip\noindent
{\em Proof of Proposition \ref{p:disk1} in the case $z\in \ms T\us$.}
First, suppose $z\in \ms T\us_{g+1}$.
As before, let $\Pi:\ms T\lras {\rm C}(\Lmd)$ be the double covering map with branch $\Sigma$ and put $\mf d=\Pi(z)\in \mathbb D_{g+1}$.
If $\mf d\not\in\ptl\mathbb D_{g+1}$, then from Proposition \ref{p:fol2}, there exists a real hyperplane $h\subset\CP_{g+2}$ such that $h\in \mf S_k\uc$ and $\mf d\in h$.
If $\mf d\in\ptl\mathbb D_{g+1}$, then taking another point $\eta\in\ptl\mathbb D_{g+1}\minus\{\mf d\}$, the hyperplane $h$ that corresponds to any one of the two points $\mf t_1\inv(\mf a(\mf d + \eta))$ satisfies $h\in \mf S_k\uc$ and $\mf d\in h$. 
Then $C=\ms T\cap \Pi\inv(h)$ is a regular minitwistor line 
whose real circle $C^{\RR}$ passes through $z$.
Since $C^{\RR} \subset C^+\cap C^-$, we have $z\in C^+$.

Next, suppose $z\in \ms T\us_i$ for some $i\le g$.
The sphere $\ms T\us_i$ is $\mathbb S^1$-invariant and on it the $\mathbb S^1$-action is a rotation which fixes the points $r_i$ and $r'_i$. Since $z\neq r_i,r'_i$ as $z\in \ms T^{k+}$, by the $\mathbb S^1$-action we may assume that $z\in \Sigma\us_i=\ptl\mathbb D_i$.
As in the proof of Proposition \ref{p:limtl}, put $q:=\varpi_k(z)\in Q\us_k$
and let $\varpi_k\inv(q) = \{p_1,\dots,p_{g+1}\}$ with $p_j\in\ms T\us_j$ for any $j$.
Then $p_i = z$, and since $z\neq r_i,r'_i$ as above,
by Proposition \ref{p:cov1},
$q\not\in\{(0,\infty),(\infty,0)\}$ in the coordinates on the quadric  $Q_k$ in the same proposition.
This means that there exists a unique pair of lines $\ell$ and $\ol\ell$ on $Q_k$ that intersect at $q$.
Then the divisor $\varpi_k\inv(\ell + \ol\ell)$ on $\ms T$ belongs to the hyperplane class by the same proposition and it restricts to the divisor $2\sum_{j=1}^{g+1}p_j$ on the curve $\Sigma$.
Since $z\in \ptl\mathbb D_i$, there exists a hyperplane
$h\subset\CP_{g+2}$ such that $\ms T\cap \Pi\inv(h)=\varpi_k\inv(\ell + \ol\ell)$. 
(The existence of such $h$ can also be proved by using the pencil $\ms P(\rho_1,\dots, \rho_g,r_{g+1})$ or $\ms P(\rho_1,\dots, \rho_g,r'_{g+1})$ and the boundary $\ptl \mf S_k$.) 
$h$ is unique from Proposition \ref{p:pn} and $h\in\ptl\mf S_k\minus\{\mf o_k,\mf o'_k\}$ as $z\neq r_{g+1},r'_{g+1}$.
In the following we show that $h$ can be moved into a hyperplane $h'\in \mf S_k\uc$ such that the restriction $h'|_{\Sigma}$ still contains $2z$ as a subdivisor.
Once this is proved, $C':=\ms T\cap \Pi\inv(h')$ is a regular minitwistor line having the point $z$ as a real node, which means $z\in (C')^+$ as required.

Let $\mf a:\Sigma\lras {\rm J}_{\Sigma}$ be the Abel-Jacobi map as before.
Since the divisor $2\sum_{j=1}^{g+1}p_j$ on the curve $\Sigma$ belongs to the hyperplane class as above, the point $\sum_{j=1}^g \mf a(p_j)$ belongs to the boundary $\ptl\mf S_k$, and it is not equal to the singular points $\mf o_k,\mf o'_k$ of $\mf S_k$ because $p_{g+1}\neq r_{g+1},r'_{g+1}$.
Now, just for $W_{g+1,\,g+1}$ as preceding to Proposition \ref{p:C0C1}, 
upon choosing 
a real tubular neighborhood $V$ of the circle $\Sigma\us_{g+1}$ in $\Sigma$,
we can take a complexification $(\ptl W_{g+1,\,g+1})^{\CC}$ of the boundary curve $\ptl W_{g+1,\,g+1}$ as 
$\{\mf a(2\xi)\in {\rm J}_{\Sigma}\set \xi\in V\}$.
Then we obtain a complexification 
$(\ptl\mf S_k)^{\CC}$ of the boundary curve $\ptl\mf S_k$ as the connected component of $\mf t\inv((\ptl W_{g+1,\,g+1})^{\CC})$ which contains $\ptl \mf S_k$, where $\mf t:{\rm J}_{\Sigma}\lras {\rm J}_{\Sigma}$ is the twice mapping as before.
These complexifications $(\ptl W_{g+1,\,g+1})^{\CC}$ and $(\ptl\mf S_k)^{\CC}$ are smooth except $\mf o$ and $\mf o_k,\mf o'_k$ respectively if we choose the neighborhood $V$ of $\Sigma\us_{g+1}$ sufficiently thin.

Put $\mf x:=\sum_{j=1}^g\mf a(p_j)\in \ptl \mf S_k$ and $\mf y :=2\mf a(p_{g+1})\in \ptl W_{g+1,\,g+1}$ for simplicity, so that $\mf t(\mf x) + \mf y = \mf o$ in ${\rm J}_{\Sigma}$.
From the definition of $\mf a$, if $\phi:\Sigma\lras\CP_{g-1}$ is the canonical map of $\Sigma$, we have
$T_{\mf y} (\ptl W_{g+1,\,g+1})^{\CC}=\langle\phi(p_{g+1})\rangle$,
where we are thinking $\CP_{g-1}$ as the space of tangent directions at any point of ${\rm J}_{\Sigma}$ and $\langle\phi(p_{g+1})\rangle$ means the line in $\CC^g$ in the direction of $\phi(p_{g+1})$.
Since $\mf t:{\rm J}_{\Sigma}\lras {\rm J}_{\Sigma}$ is a covering map,
it follows that 
\begin{align}\label{Tx}
T_{\mf x} (\ptl\mf S_k)^{\CC}=\langle\phi(p_{g+1})\rangle.
\end{align}
On the other hand, as $\pi(p_j)\neq \pi(p_{l})$ for any $j\neq l$,
from Proposition \ref{prop:he1},
we have $\dim |\sum_{j=1,j\neq i}^g p_j| = 0$. 
Hence, the theta divisor $W_{g-1}\subset {\rm J}_{\Sigma}$ is smooth at the point $\sum_{j=1,j\neq i}^g\mf a(p_j)$. 
So the locus 
\begin{align*}
\mf P_i :&=\big\{ \mf a(q_1+\dots+q_g)\in {\rm J}_{\Sigma}\set q_j\in\Sigma {\text{ if }}j\neq i,\,q_i=p_i \,(=z) \big\}\\
&= \mf a(p_i) + W_{g-1}
\end{align*}
is smooth at the point $\mf x=\sum_{j=1}^g\mf a(p_j)$ and 
\begin{align}\label{Tx2}
T_{\mf x} \mf P_i=\big\langle\phi(p_1), \dots, 
\hat{\phi(p_i)},\dots,
\phi(p_g)\big\rangle
\end{align}
in the same sense as \eqref{Tx} where widehat means removal.
Since the image $\phi(\Sigma)$ is a non-degenerate rational normal curve in $\CP_{g-1}$, 
\eqref{Tx2} means that $\mathbb P(T_{\mf x}\mf P_i) \cap \phi(\Sigma) = 
\{\phi(p_1), \dots, 
\hat{\phi(p_i)},\dots,
\phi(p_g)\}$.
Since  $\phi(p_j)\neq \phi(p_{g+1})$ for any $j\le g$ as $\ol K_j\cap \ol K_{g+1}=\emptyset$, this means $\phi(p_{g+1})\not\in \mathbb P(T_{\mf x}\mf P_i)$.
From \eqref{Tx} and \eqref{Tx2}, this means that the intersection of the curve $(\ptl\mf S_k)^{\CC}$ and the hypersurface $\mf P_i$ is transversal at the point $\mf x$.
Hence, so is the intersection of the real locus $\mf P_i\us$ and $\ptl\mf S_k$ at the same point.
It follows that the intersection $\mf P_i\us\cap \mf S_k$ is a (real) smooth curve in the surface $\mf S_k$.
The hyperplane $h'\subset\CP_{g+2}$ determined by a point of this curve gives the required one.
\proofend

\medskip

 \noindent{\em Proof of Proposition \ref{p:disk1} in the case $z\not\in \ms T\us$.}
We again make use of the map $\varpi_k:\ms T\lras Q_k$.
Since $z\not\in\ms T\us$, the point $q=\varpi_k(z)\in Q_k$ is not a real point, and moreover, since $z\neq \bm p_{\infty},\ol{\bm p}_{\infty}$ as $z\in\ms T^{k+}$, $q\not\in\{(0,0),(\infty,\infty)\}$ in the coordinates of Proposition \ref{p:cov1}. 
From this, it is elementary to see that 
there always exists at least one pair of lines $\ell$ and $\ol \ell$ on the quadric $Q_k$ such that $q\in \ell$, $q\not\in \ol\ell$ and $(\ell\cup\ol \ell)\cap \{(0,0),(0,\infty),(\infty,0),(\infty,\infty)\} = \emptyset$.
Put $D:=\varpi_k\inv(\ell)$. Then $z\in D$ and $z\not\in\ol D$.
These are smooth rational curves on $\ms T$ and $D+\ol D\in\ms H_k$.
Either the pullbacks of $A$-lines or $B$-lines on $Q_k$ correspond to $ + $ -disks under the degeneration of minitwistor lines when a point of $M_k$ approaches $\ms H_k$ as explained at the end of the previous subsection, so by replacing $D$ or $\ol D$ if necessary, we may suppose that $ D $ is a degeneration of $+$ -disks.

As in the previous proof, write $D\cap \ol D= \{p_1,\dots,p_{g+1}\}$ with $p_i\in\ms T\us_i$ and all these are nodes of $D+\ol D$, which  are all singularities of $D+\ol D$. 
This time, the point $z$ is a smooth point of $D+\ol D$.

The rest is the same as the proof of \cite[Lemma 5.3]{HN22} in that we consider displacements of the curve $D+\ol D$ in $\ms T$, which are equisingular at $p_1,\dots, p_g$ and which remain to pass through $z$ and $\ol z$.
Using $D^2 = \ol D^2 = 0$ and $D\ol D = g+1$, we can show that the versal family of such displacements is parameterized by a smooth holomorphic curve and they smooth out the node $p_{g+1}$.
Taking the real locus of this curve and restricting it to the half which is included in $M_k$, we obtain a real 1-dimensional family of regular minitwitor lines through $z$ (and $\ol z$), and $z$ belongs to $+$-disks as $z\in D$ initially and $D$ is a degeneration of $+$ - disks as above.
\proofend

\medskip From Lemma \ref{l:pm} and Proposition \ref{p:disk1}, we obtain the following corollary.

\begin{corollary}\label{c:disk}
For any boundary data $1\le k\le 2^{g-1}$, 
$$\bigcup_{C\in M_k\uc} C^+ = \ms T^{k+}.$$
\end{corollary}

Thus, we are in the same situation as in the case of $g=1$ that we have reached at the end of \cite[Section 5.1]{HN22}, and 
from the inclusion $(\ms C_k\uc)^+\subset\ms T^{k+}\times M_k\uc$, we obtain the double fibration
\begin{align}\label{d:double3}
 \xymatrix{ 
&(\ms C_k\uc)^+ \ar[dl]_{\pr_1} \ar[dr]^{\pr_2} &\\
\ms T^{k+} && M_k\uc
 }
\end{align}
The point is that the projection $\pr_1$ is surjective by Corollary \ref{c:disk}.
Thus, the open complex surface $\ms T^{k+}$ can be regarded as the minitwistor space whose twistor lines are disks \cite{LM09,Nkt09}.
In particular, the image $\pr_2(\pr_1\inv(q))\subset M_k\uc$ is a null surface if $q\in \ms T_{g+1}\us$ and is a timelike geodesic if $q\not\in \ms T_{g+1}\us$.
All maximal null surfaces and maximal timelike geodesics are in this form.
The null surfaces are totally geodesic with respect to the affine connection of the EW structure.

\begin{corollary}\label{c:aut2}
Any automorphism of the punctured EW space $M_k\uc$ induces a holomorphic automorphism of $\ms T^{k+}$. 
\end{corollary}

\proof
Any automorphism of a Lorentzian EW space maps timelike geodesics to timelike geodesics and null surfaces to null surfaces.
Hence, since $\ms T^{k+}$ is exactly the space of these objects as above, any automorphism of $M\uc_k$ induces a diffeomorphism of $\ms T^{k+}$. Every diffeomorphism obtained this way (from an isomorphism of a Lorentzian EW space) is holomorphic in the general setting as proved in \cite[Corollary 5.8]{HN22}.
\proofend

\medskip
Thus, we have obtained a homomorphism $\Aut(M_k\uc) {\lras} \Aut(\ms T^{k+})$.

\subsection{The automorphism groups and the moduli spaces}\label{ss:autmod}
As a final step for completing the sequence \eqref{seqa}, we show the following. Recall that $\ms T$ admits a (rational) conic bundle structure over $\Lmd$ by the (rational) projection $f=\pi\circ\Pi$.


\begin{proposition}\label{p:aut1}
Every holomorphic automorphism of the open minitwistor space $\ms T^{k+}$ 
extends to that of $\ms T$, and the extension maps fibers of the projection $f$ to fibers of the projection $f$.
If the automorphism of $\ms T^{k+}$ is induced from an automorphism of the EW space $M_k\uc$, then the extension commutes with the real structure and it preserves the real minitwistor space $\ms T\us_{g+1}$.
\end{proposition}

\proof
The idea is similar to the case of $g=1$ given in \cite[Theorem 5.11]{HN22}, so we just sketch the proof. Write $l_i = l_i(\rho_i)$ and $l'_i = l_i(\rho'_i)$ for $1\le i\le g$ as before and recall $\ms T^{k+} = \ms T\minus
(l_1\cup\dots\cup l_g)\cup (\ol l'_1\cup\dots\cup\ol l'_g)$. 
The transformations of all lines on $\ms T$ into the minimal resolution $\tilde{\ms T}$ of $\ms T$ are $(-1)$-curves.
We blow down the transformations of the lines $l_1,\dots, l_g$ and $\ol l'_1,\dots, \ol l'_g$. Let $\ms T'$ be the resulting smooth surface.
Then the image of the exceptional curves $E$ and $\ol E$ in $\tilde{\ms T}$ into $\ms T'$ are $(-1)$-curves. Blowing down these two, we obtain another smooth surface $\ms T''$, which is isomorphic to $\qdr$. Let $e$ and $\ol e$ be the images of $E$ and $\ol E$ into $\ms T''$ respectively.
Then there is an obvious isomorphism $\ms T^{k+}\simeq \ms T''\minus\{e,\ol e\}$, so any holomorphic automorphism $\phi^+$ of $\ms T^{k+}$ induces a holomorphic automorphism of $\ms T''\minus\{e,\ol e\}$.
By Hartogs' theorem, this extends to an automorphism $\phi''$ of $\ms T''$ holomorphically.
This automorphism induces an automorphism $\phi$ of $\ms T$ through $\ms T'$ and $\tilde {\ms T}$.

Since fibers of the projection $f$ are transformed to $(1,1)$-curves in $\ms T''\simeq\qdr$ passing through the points $e$ and $\ol e$ and these two points do not belong to the same $(1,0)$-curve nor the same $(0,1)$-curve, 
$\phi$ on $\ms T$ maps fibers of $f$ to fibers of $f$.
Moreover, for the same reason, each of the reducible fibers over the $2g$ points $\lmd_1,\lmd'_1,\dots,\lmd_g,\lmd'_g\in\Lmd$ is mapped to one of these fibers, and each of the remaining reducible fibers over $\lmd_{g+1}$ and $\lmd'_{g+1}$ are mapped to one of these two fibers.

If the automorphism $\phi$ of $\ms T$ is induced from an automorphism of $M_k\uc$, then the commutativity of $\phi$ with the real structure $\sigma$ on $\ms T$ is automatic if we think of each point of $\ms T$ as a complex null direction $[v]$ at a point of $M_k\uc$ and that $\sigma$ assigns the complex conjugation $[\ol v]$ to $[v]$.
Further, $\phi$ preserves $\ms T\us_{g+1}$ from the above information about the automorphism induced on $\Lmd$. (This also from the fact that 
the real minitwistor space $\ms T\us_{g+1}$ is the space of maximal null surfaces in $M_k$ as above.)
\proofend

\medskip
Combining the results obtained in this section so far, we can show:
\begin{proposition}\label{p:aut}
Any automorphism of the EW space $M_k$ naturally induces a holomorphic automorphism of $\ms T$ that commutes with the real structure and that maps fibers of the projection $f$ to fibers of the projection $f$ and
that preserves the real minitwitor space $\ms T\us_{g+1}$. 
\end{proposition}

\proof
Let $\phi:M_k\lras M_k$ be an automorphism preserving the EW structure. By Corollary \ref{c:aut1}, $\phi$ preserves $M_k\uc$. Let $\phi\uc:=\phi|_{M_k\uc}$. From Corollary \ref{c:aut2}, $\phi\uc$ induces a holomorphic automorphism $\phi^+$ of $\ms T^{k+}$.
From Proposition \ref{p:aut1}, $\phi^+$ extends to a holomorphic automorphism of $\ms T$ which satisfies the properties in the theorem.
\proofend

\medskip
Let $\Phi:\ms T\lras \ms T$ be the covering transformation of the double covering $\Pi:\ms T\lras {\rm C}(\Lmd)$. Of course, this covers the identity map on $\Lmd$ and from the equation \eqref{T} of $\ms T$ and the concrete form of $\sigma$, the involution $\Phi$ commutes with the real structure.
Further, $\Phi$ obviously preserves every hyperplane belonging to 
$\Pi^*(\mf S_k)=\{\Pi\inv(h)\set h\in \mf S_k\}$ for each $k$.
Hence, $\Phi$ induces an involution of the EW space $M_k$ for each $k$ preserving the EW structure.
If we use the same letter $\Phi$ to mean this involution on $M_k$, then 
$\Phi$ fixes any point of $\Pi^*(\tilde{\mf S}_k)$, where $\tilde{\mf S}_k$ is the `resolution' of $\mf S_k$ given in Definition \ref{d:ann}. 
Furthermore, since the restriction of $\Pi$ to a generic minitwistor line is of degree one over its image, $\Phi$ does not preserve such minitwistor lines. 
Therefore, $\Phi$ on $M_k\simeq\mathbb S^2\times I$ is non-trivial and is a reflection with respect to the annulus $\tilde{\mf S}_k$.
In particular, $\Phi$ is not included in the $\mathbb S^1$-action \eqref{S1act}.
Topologically, this annulus can be identified with $\mathbb S^1\times I\subset\mathbb S^2\times I$ (see Figure \ref{fig:Seifert4}).

The hyperelliptic involution $\tau$ of $\Sigma$ also induces a holomorphic involution of $\ms T$ that commutes with $\sigma$ and we can see that it induces an involution of the EW space $M_k$ for each $k$. 
But from the explicit equation \eqref{T}, this involution is the composition of the reflection $\Phi$ with the involution $(-1)\in\mathbb S^1$.


\begin{theorem}\label{t:aut}
For any boundary data $1\le k\le 2^{g-1}$, the identity component of $\Aut(M_k)$ is $\mathbb S^1$ given by \eqref{S1act}.
If the hyperelliptic curve $\Sigma$ satisfies $\Aut(\Sigma) = \{\id,\tau\}$, then $\Aut(M_k)$ is generated by the $\mathbb S^1$-action and the involution $\Phi$.
\end{theorem}

\proof
Let $\ms T$ be the minitwistor space from which $M_k$ is obtained and $\Psi$  any automorphism of the EW space $M_k$.
By Proposition \ref{p:aut}, $\Psi$ induces a holomorphic automorphism 
of $\ms T$ which maps fibers of $f$ to fibers of $f$.
Hence it covers a holomorphic automorphism of $\Lmd$ and it preserves the set of the branch points of $\pi:\Sigma\lras\Lmd$.
Thus, we obtain a sequence $\Aut(M_k)\lras\Aut(\ms T)\lras\Aut(\Sigma)$. But the identity component $\Aut_0(\Sigma)$ is $\{\id\}$, so if $\Psi$ belongs to $\Aut_0(M_k)$, then the induced automorphism $\phi$ of $\ms T$ preserves each fiber of $f$.
Further, such an automorphism $\phi$ of $\ms T$ cannot exchange $\bm p_{\infty}$ and $\ol{\bm p}_{\infty}$ because it belongs to the identity component.
This means that $\phi$ descends to an automorphism of $\ms T''\simeq\qdr$ through the blowing down $\ms T\lras\ms T''\simeq\qdr$ in the proof of Proposition \ref{p:aut1} and the resulting automorphism preserves each of the two points $e$ and $\ol e$.
With reality, this means that $\phi\in\mathbb S^1$.
Hence, $\Aut_0(M_k) = \mathbb S^1$.

If $\Aut(\Sigma) = \{\id,\tau\}$, then since $\pi = \pi\circ\tau$, the automorphis $\phi$ of $\ms T$ induced by any automorphism of $M_k$ preserves each fiber of $f$. Therefore, using the covering transformation $\Phi$ of $\Pi:\ms T\lras {\rm C}(\Lmd)$ if necessary, either $\phi$ or $\Phi\circ\phi$ preserves each of $\bm p_{\infty}$ and $\ol{\bm p}_{\infty}$. If $\phi$ does, then $\phi\in\mathbb S^1$ as above.
Otherwise, $\Phi\circ\phi\in \mathbb S^1$. This proves the latter assertion of the theorem.%
\proofend

\medskip

In a similar way, we can show:

\begin{proposition}\label{p:aut3}
If $M_k$ and $M_l$ $(1\le k,l\le 2^{g-1})$ are the EW spaces obtained from the same minitwistor space $\ms T$, and if $\Aut (\Sigma) = \{\id, \tau\}$, then
$M_k\simeq M_l$ as EW spaces if and only if $k=l$.
\end{proposition}

\proof
Recall $M_k = M_k\uc\cup(\bm I_k\cup\bm I'_k)$ and $M_l = M_l\uc\cup(\bm I_l\cup\bm I'_l)$.
From the proof of Corollary \ref{c:aut1}, in terms of the local minitwistor correspondence,
the EW structure on a neighborhood of any point $p$ of $\bm I_k\cup\bm I'_k$ is characterized by the property $\dim |C_p| = 1$ where $C_p$ is the smooth, transformed minitwistor line corresponding to $p$.
Hence, any isomorphism $M_k\lras M_l$ as EW spaces maps $\bm I_k\cup\bm I'_k$ to $\bm I_l\cup\bm I'_l$ and the complement $M_k\uc$ to $M_l\uc$.
Then in the same way as Corollary \ref{c:aut2} and Proposition \ref{p:aut1}, via the open minitwistor spaces $\ms T^{k+}$ and $\ms T^{l+}$, we obtain a holomorphic isomorphism $\phi:\ms T\lras \ms T$ preserving the real structure which covers a holomorphic isomorphism of $\Lmd$. So we obtain an automorphism of $\Sigma$, but by the assumption $\Aut (\Sigma) = \{\id, \tau\}$, the isomorphism of $\Lmd$ has to be the identity map.  
But if $k\neq l$, this contradicts that the automorphism $\phi$ of $\ms T$ is an extension of an isomorphism $\ms T^{k+}\lras\ms T^{l+}$. 
\proofend

\medskip

Finally, we discuss the moduli space of the present EW spaces. 
Let $B:=\{\lmd_i,\lmd'_i\set 1\le i\le g+1\}$ and $\hat B:=\{\hat\lmd_i,\hat\lmd'_i\set 1\le i\le g+1\}$ be two sets of the branch points on $\Lmd\us=\RP_1$ and
$\ms T$ and $\hat{\ms T}$ the minitwistor spaces determined by $B$ and $\hat B$ respectively.
Let $M_k$ and $\hat M_{\hat k}$ be the EW spaces associated to $\ms T$ and $\hat{\ms T}$ with boundary data being $k$ and $\hat k$ respectively,
and $\Sigma$ and $\hat\Sigma$ the hyperelliptic curves determined by $B$ and $\hat B$ respectively. 

\begin{theorem}\label{t:mod}
If the EW spaces $M_k$ and $\hat M_{\hat k}$ are isomorphic, then the two hyperelliptic curves $\Sigma$ and $\hat\Sigma$ are isomorphic.
%
\end{theorem}

\proof
The proof is almost the same as Theorem \ref{t:aut} and Proposition \ref{p:aut3}.
Suppose that $\phi:M_k\lras \hat M_{\hat k}$ is an isomorphism as EW spaces.
In the same way to Corollary \ref{c:aut1} using Proposition \ref{p:C0C1}, $\phi$ has to map the rotation axis $\bm I_k\cup{\bm I'}_k\subset M_k$ to the axis $\hat{\bm I}_{\hat k}\cup\hat{\bm I}'_{\hat k}\subset\hat M_{\hat k}$.
So $\phi$ maps the complement $M_k\uc$ to $\hat M\uc_{\hat k}$.
By Corollary \ref{c:aut2} and Proposition \ref{p:aut1}, this induces a holomorphic isomorphism $\ms T\lras \hat{\ms T}$ preserving the real strcutures which covers a holomorphic isomorphism of $\Lmd$ that maps $B$ to $\hat B$.
This means $\Sigma$ and $\hat\Sigma$ are isomorphic.
\proofend

%

\medskip
Since the moduli space of hyperelliptic curves of genus $g$ is $(2g-1)$-dimensional, Theorem \ref{t:mod} implies the following.

\begin{corollary}\label{c:mod}
The moduli spaces of the EW structures arising from the minitwistor spaces obtained from hyperelliptic curves of genus $g$ are of $(2g-1)$-dimensional.
\end{corollary}

%
%
%
%
%

\section{Appendix}
In \cite{Hi24}, Hitchin obtained the minitwistor spaces whose minitwistor lines are nodal curves from the twistor spaces of ALE gravitational instantons. 
In this section, we show that they are isomorphic to our minitwistor spaces including the real structure, in the case of A$_k$-instantons, and also that the linear systems in which minitwistor lines belong are the same.

Take real numbers $a_1<a_2<\dots<a_{k+1}$ and
let $Z$ be a hypersurface in the total space of the vector bundle $\ms O(k+1)\oplus\ms O(k+1)\oplus\ms O(2)$ over $\CP_1$ defined by the equation
\begin{align}\label{tw1}
xy = (z-a_1u)(z-a_2u)\dots(z-a_{k+1}u),
\end{align}
where $u$ is an affine coordinate on $\CP_1$, $x,y\in\ms O(k+1)$ 
and $z\in \ms O(2)$.
This is smooth away from two singular points lying over $u=0,\infty$, and is invariant under the real structure $\sigma$ induced from the quoternionic structure on $\CC^2$ defined by 
$(z_1,z_2)\longmapsto (-\ol z_2,\ol z_1)$.
Particular resolution of the singularities of $Z$ is the twistor space of an ALE  hyperK\"ahler metric on the minimal resolution of $\CC^2/\CCC$, $\CCC\subset\SU(2)$ with $\CCC\simeq\ZZ_{k+1}$.

The scalar multiplication on $\CC^2/\CCC$ induces a $\CC^*$-action on the hypersurface $Z$ which is given by 
\begin{align*}
(x,y,z,u)\stackrel{t}\longmapsto
\big(t^{k+1}x, t^{k+1}y, t^2 z, t^2u\big),\quad t\in\CC^*,
\end{align*}
which commutes with $\sigma$ iff $|t|=1$. When $k$ is odd, putting $k+1=2l$, $l\ge 1$ and $s=t^2$, we obtain an effective $\CC^*$-action
\begin{align}\label{s}
(x,y,z,u)\stackrel{s}\longmapsto
\big(s^l x, s^ly,s z, su\big),\quad s\in\CC^*.
\end{align}

Hitchin \cite{Hi24} obtained compact minitwistor spaces from these twistor spaces in the case $k$ odd as follows.
Putting $k+1=2l$ as above, consider the fiber of the projection $Z\lras\CP_1$ over the point $u=1$ of $\CP_1$, which is
\begin{align}\label{Hi1}
xy = (z-a_1)(z-a_2)\dots(z-a_{2l})
\end{align}
in $\CC^3$.
This surface can be naturally compactified by thinking $z$ as an affine coordinate on $\CP_1$ and $(x,y)$ as affine fiber coordinates on the $\CP_2$-bundle $\mathbb P(\ms O(l)\oplus\ms O(l)\oplus\ms O)$ over $\CP_1$ this time.
Let $S$ be the compact projective surface obtained this way.
This is non-singular, having a conic bundle structure over $\CP_1$.
Letting $w$ be a fiber coordinate on the $\ms O$-factor,
the divisor added in the compactification is two sections $x=w=0,\, y=w=0$ and the fiber over $z=\infty$.
The two sections are $(-l)$-curves on $S$.
In \cite{Hi24} these are written $D_1,D_2$.
As a conic bundle $S\lras\CP_1$, the reducible fibers are exactly over the $2g$ points $z=a_1,\dots, a_{2l}$.

While the affine surface \eqref{Hi1} or its compactification $S$ is not $\sigma$-invariant, they admit a natural real structure using the $\CC^*$-action \eqref{s}.
Concretely, it is given by the composition of $\sigma$, which maps the fiber over $u=1$ to the fiber over $u=-1$, and the isomorphism between these fibers obtained by letting $s=-1$ in \eqref{s}.
As in \cite{Hi24} we denote $\tau$ for this real structure on $S$. 

\begin{proposition}
The minitwistor space $(S,\tau)$ associated to an ALE hyper-K\"ahler metric on the minimal resolution of A$_{2l-1}$-singularity as above is isomorphic to the minimal resolution of $\ms T$ equipped with the real structure induced from $\sigma$ on $\ms T$, where $\ms T$ is a minitwistor space associated to a hyperelliptic curve of genus $(l-1)$.
\end{proposition}

\proof
Using that $x,y\in\ms O(k+1)$ and $z\in \ms O(2)$, 
the real structure $\sigma$ on $Z$ is given, in the above coordinates, by 
\begin{align}\label{rsZ}
(x,y,z,u)\stackrel{\sigma}\longmapsto \Big( (-1)^{k+1}\frac{\ol y}{\ol u^{k+1}},
\frac{\ol x}{\ol u^{k+1}}, -\frac{\ol z}{\ol u^{2}}, 
-\frac1{\ol u}\Big).
\end{align}
If $k+1=2l$ as above, restricting to $u=1$ and taking the composition  as above, the real structure $\tau$ on $S$ is given by
\begin{align}\label{rsS}
(x,y,z)\stackrel{\tau}{\longmapsto} \big(
(-1)^l \ol y, (-1)^l \ol x, \ol z\big).
\end{align}
From \eqref{Hi1}, it follows that the real locus $S^{\tau}$ consists of $l$ smooth spheres which lie over the intervals 
\begin{itemize}
\setlength{\parskip}{-2mm}
\item
$[a_1,a_2],\, [a_3,a_4],\dots, [a_{2l-1},a_{2l}]$
if $l$ is odd,
\item
$[a_2,a_3],\, [a_4,a_5],\dots, [a_{2l-2},a_{2l-1}]$ and $[a_{2l},a_1]$
if $l$ is even,
\end{itemize}
where the last interval $[a_{2l},a_1]$ means the one in $\RP_1$ bounded by $a_{2l}$ and $a_1$ (where $a_1<a_{2l}$ from the choice) and which includes $\infty=-\infty$.

On the other hand, from \eqref{T1} and \eqref{rs}, $(\ms T,\sigma)$ is defined by
\begin{gather}\label{ts}
uv = - \prod_{i=1}^{l}(z-\lmd_i)(z-\lmd'_i)
\qandq
(u,v,z)\stackrel{\sigma}\longmapsto (\ol v,\ol u,\ol z).
\end{gather}
If $l$ is odd, then we put $\lmd_i = a_{2i-1}$ and 
$\lmd'_i=a_{2i}$ for any $1\le i\le l$ and  
consider the mapping $(x,y,z)\longmapsto (u,v,z) = (x,-y,z)$. 
This transforms the equation \eqref{Hi1} and $\tau$, the real structure \eqref{rsS}, to \eqref{ts} respectively.
Further, since $l=g+1$ and the rational normal curve $\Lmd$ is in $\CP_{g+1}$, from the definition of the coordinates in \eqref{T1}, $u,v$ can be considered as fiber coordinates on the line bundle $\ms O(l)$.
So the above mapping naturally extends holomorphically to the mapping from $S$ to $\ms T$ and it is exactly the contraction of the two sections added in the compactification.
Hence we have seen that $(S,\sigma)$ is isomorphic to the minimal resolution of $(\ms T,\tau)$ in the case $l$ odd.

If $l$ is even, because of the discrepancy of the above $l$ intervals, we cannot use $a_j$-s for $\lmd_i$ and $\lmd'_i$.
Instead, we take any real number $\bbb\in (a_1,a_2)$ and let $\psi\in\PGL(2,\RR)$ be any M\"obius transformation which preserves the orientation and which satisfies $\psi(\bbb) = \infty$. Put $\lmd_i = \psi(a_{2i})$ for $1\le i\le l$, $\lmd'_i = \psi(a_{2i+1})$ for $1\le i\le l-1$ and $\lmd'_l = \psi(a_1)$. Then $\lmd_1<\lmd'_1<\dots<\lmd_l<\lmd'_l$. 
Using the letter $z'$ instead of $z$ for a coordinate on $\ms T$, it is possible to see that there uniquely exists a positive real number $c$ such that 
the mapping $(x,y,z)\longmapsto (u,v,z') = (x,-cy,\psi(z))$
transforms the equations of $S$ and $\tau$ into those of $\ms T$ and $\sigma$ respectively.
Then for the same reason as the case $l$ odd, this gives the identification of $(\ms T,\sigma)$ as the contraction of $D_1$ and $D_2$ in $(S,\tau)$. 
\proofend

\begin{proposition}
Under the isomorphism in the previous proposition, the minitwistor lines in $S$ and the pullbacks of minitwistor lines in $\ms T$ under the minimal resolution belong to the same linear system.
\end{proposition}

\proof
Let $C$ be a fiber of the projection $S\lras\CP_1$. 
Then by \cite[Theorem 1]{Hi24} the minitwistor lines in $S$ belong to the linear system $|lC + D_1 + D_2|$, where $D_1$ and $D_2$ are the two sections of $S\lras\CP_1$ added in the compactification. Note that $D_2 = \tau(D_1)$. 
On the other hand, an irregular minitwistor line in $\ms T$ is the inverse image of $g+1=l$ generating lines of the cone under the double covering map $\Pi:\ms T\lras {\rm C}(\Lmd)$ and its pullback to the minimal resolution of $\ms T$ is exactly the sum of $l$ fibers of $\tilde{\ms T}\lras\CP_1$ plus the exceptional curves of the minimal resolution. Therefore the class is also $lC + D_1 + D_2$.
\proofend

\medskip
Using $D_1^2 = D_2^2 = -l$, we readily see that the complete linear system $|lC + D_1 + D_2|$ is base point free and induces a birational mapping into $\CP_{l+2}$ which is exactly the contraction of $D_1$ and $D_2$.
Since $(lC+D_1+D_2)^2 = -2l + 4l = 2l$, the degree of the image surface is $2l$.
Of course, this agrees with the degree of $\ms T$ in $\CP_{g+3}$
which was $2\deg {\rm C}(\Lmd) = 2(g+1)$ as $l=g+1$.

\begin{remark}
{\em
In \cite{Hi21, Hi24}, Hitchin obtained a $\tau$-conjugate pair of pencils on the minitwistor space $S$ from the central sphere in the fibers of $Z\lras\CP_1$ over $u=0$ and $u=\infty$ in the above coordinates. If $D+\tau(D)$ is the sum of members of the two pencils, then the image curve of it under the contraction $S\lras\ms T$ belongs to future or past infinity of the EW space $M_k$ from the description of the curves corresponding to future and past infinity $\ptl M_k=\ms H_k\cup{\hat{\ms H}}_k$ given in Proposition \ref{prop:charW}.
From this, it seems likely that ALE gravitational instantons of type A$_{2l-1}$ and the present EW spaces obtained from hyperelliptic curve of genus $l+1$ are adjacent in the sense that they share a common boundary in the dual projective space $\RP_{g+3}^*$ (= the space of real hyperplanes in $\CP_{g+3}$), which would be an analogue of the fact that the hyperbolic ball and the standard deSitter space (or more precisely its quotient by a fixed-point free involution) are adjacent in the sense that they share a common boundary in the dual projective space $\RP_3^*$, where we are thinking the two spaces as the EW spaces associated to a real smooth quadric in $\CP_3$ that has real points.
}
\end{remark}

\end{document}